\documentclass[]{article}
\makeatletter\if@twocolumn\PassOptionsToPackage{switch}{lineno}\else\fi\makeatother

\usepackage{amsfonts,amssymb,amsbsy,latexsym,amsmath,tabulary,graphicx,times,caption,fancyhdr}
\usepackage[utf8]{inputenc}

\usepackage{url,multirow,morefloats,floatflt,cancel,tfrupee}
\makeatletter

\AtBeginDocument{\@ifpackageloaded{textcomp}{}{\usepackage{textcomp}}}
\makeatother
\usepackage{colortbl}
\usepackage{xcolor}
\usepackage{pifont}
\usepackage[nointegrals]{wasysym}
\usepackage{amsmath}
\usepackage{graphicx}
\usepackage{amssymb}
\usepackage{epstopdf}
\usepackage{amsmath,amssymb,amsthm}
\usepackage{bm}
\usepackage{bbm}
\usepackage{mathtools}
\usepackage{gensymb}
\usepackage{setspace}
\usepackage{float}
\usepackage{comment}
\usepackage{caption}
\usepackage{indentfirst}
\newcommand{\bfn}{\textbf{n}}
\usepackage{epsfig}
\makeatletter

\def\mcWidth#1{\csname TY@F#1\endcsname+\tabcolsep}

\def\cAlignHack{\rightskip\@flushglue\leftskip\@flushglue\parindent\z@\parfillskip\z@skip}
\def\rAlignHack{\rightskip\z@skip\leftskip\@flushglue \parindent\z@\parfillskip\z@skip}

\@ifundefined{etal}{}{}

\usepackage{ifxetex}
\ifxetex\else\if@twocolumn\@ifpackageloaded{stfloats}{}{\usepackage{dblfloatfix}}\fi\fi

\AtBeginDocument{
\expandafter\ifx\csname eqalign\endcsname\relax
\def\eqalign#1{\null\vcenter{\def\\{\cr}\openup\jot\m@th
  \ialign{\strut$\displaystyle{##}$\hfil&$\displaystyle{{}##}$\hfil
      \crcr#1\crcr}}\,}
\fi
}

\AtBeginDocument{%
  \@ifpackageloaded{endfloat}%
   {\renewcommand\efloat@iwrite[1]{\immediate\expandafter\protected@write\csname efloat@post#1\endcsname{}}}{\newif\ifefloat@tables}%
}%

\def\BreakURLText#1{\@tfor\brk@tempa:=#1\do{\brk@tempa\hskip0pt}}
\let\lt=<
\let\gt=>
\def\processVert{\ifmmode|\else\textbar\fi}

\@ifundefined{subparagraph}{
\def\subparagraph{\@startsection{paragraph}{5}{2\parindent}{0ex plus 0.1ex minus 0.1ex}%
{0ex}{\normalfont\small\itshape}}%
}{}

\newcommand\role[1]{\unskip}
\newcommand\aucollab[1]{\unskip}
  
\@ifundefined{tsGraphicsScaleX}{\gdef\tsGraphicsScaleX{1}}{}
\@ifundefined{tsGraphicsScaleY}{\gdef\tsGraphicsScaleY{.9}}{}
\def\checkGraphicsWidth{\ifdim\Gin@nat@width>\linewidth
	\tsGraphicsScaleX\linewidth\else\Gin@nat@width\fi}

\def\checkGraphicsHeight{\ifdim\Gin@nat@height>.9\textheight
	\tsGraphicsScaleY\textheight\else\Gin@nat@height\fi}

\def\fixFloatSize#1{}
\let\ts@includegraphics\includegraphics

\def\inlinegraphic[#1]#2{{\edef\@tempa{#1}\edef\baseline@shift{\ifx\@tempa\@empty0\else#1\fi}\edef\tempZ{\the\numexpr(\numexpr(\baseline@shift*\f@size/100))}\protect\raisebox{\tempZ pt}{\ts@includegraphics{#2}}}}

\AtBeginDocument{\def\includegraphics{\@ifnextchar[{\ts@includegraphics}{\ts@includegraphics[width=\checkGraphicsWidth,height=\checkGraphicsHeight,keepaspectratio]}}}

\DeclareMathAlphabet{\mathpzc}{OT1}{pzc}{m}{it}

\def\URL#1#2{\@ifundefined{href}{#2}{\href{#1}{#2}}}

\def\UrlOrds{\do\*\do\-\do\~\do\'\do\"\do\-}%
\g@addto@macro{\UrlBreaks}{\UrlOrds}

\edef\fntEncoding{\f@encoding}

\makeatother

\newif\ifmultipleabstract\multipleabstractfalse%
\newtheorem{problem}{\n{Problem}}[section]

\newtheorem{remark}{\n{Remark}}

\newtheorem{proposition}{\n{Proposition}}[section]

\font\n=cmcsc10

\newcommand{\rvline}{\hspace*{-\arraycolsep}\vline\hspace*{-\arraycolsep}}

   \newcommand{\norm}[1]{\left\lVert #1\right\rVert}
\newcommand{\abs}[1]{\left\lvert #1 \right\rvert}

\makeatletter

\def\wileyIndent{1pt}
\usepackage[paperheight=10in,paperwidth=6.5in,margin=2cm,headsep=.5cm,top=2.5cm,headheight=1cm]{geometry}

\renewenvironment{abstract}
{\vspace*{-1pc}\trivlist\item[]\leftskip\wileyIndent\hrulefill\par\vskip4pt\noindent\textbf{\abstractname}\mbox{\null}\\}{\par\noindent\hrulefill\endtrivlist}

\usepackage[]{footmisc}

\def\author#1{\gdef\@author{\hskip-\dimexpr(\tabcolsep)\hskip\wileyIndent\parbox{\dimexpr\textwidth-\wileyIndent}{\centering\bfseries#1}}}

\def\title#1{\linespread{1}\gdef\@title{\centering\bfseries\ifx\@articleType\@empty\else\@articleType\\\fi#1}}

\let\@articleType\@empty \def\articletype#1{\gdef\@articleType{{\normalfont\itshape#1}}}

\AtBeginDocument{\fancypagestyle{headings}{\fancyhf{}\fancyhead[C]{\RunningHead}\fancyhead[R]{\thepage}}\pagestyle{headings}}

 \def\audegree#1{}

\date{}

\makeatother

\usepackage[T1]{fontenc}
\makeatother
\usepackage[numbers,sort&compress]{natbib}

\def\thanksspace{{\phantom{\textsuperscript{\thefootnote}}}}
 \makeatletter
\def\@fnsymbol#1{\ensuremath{\ifcase#1\or $1$\or
   $2$\or $3$\else\@ctrerr\fi}}
\begin{document}

\title{Unfitted mixed finite element methods for elliptic interface problems}
\author{Najwa~Alshehri\textsuperscript{A}\thanks{E-mail:                     
                    najwa.alshehri@kaust.edu.sa}{\thanksspace}, Daniele~Boffi\textsuperscript{A,B}\thanks{E-mail:                     
                    daniele.boffi@kaust.edu.sa}{\thanksspace}\space and Lucia~Gastaldi\textsuperscript{C}\thanks{E-mail:                     
                    lucia.gastaldi@unibs.it}{\thanksspace}~\\[-3pt]\normalsize\normalfont  \itshape ~\\
\textsuperscript{A}{Computer, electrical and mathematical sciences and engineering division\unskip, King Abdullah University of Science and Technology (KAUST)\unskip, Thuwal\unskip, 23955\unskip, Saudi Arabia}~\\
\textsuperscript{B}{Dipartimento di Matematica `F. Casorati'\unskip,
University of Pavia\unskip, I-27100 Pavia\unskip, Italy}~\\
\textsuperscript{C}{Dipartimento di Ingegneria Civile, Architettura,
Territorio, Ambiente e di Matematica\unskip, Universit\`{a} degli Studi di
Brescia\unskip, I-25123 Brescia\unskip, Italy}}

\def\RunningHead{Unfitted FEM for interface problems}
\def\RunningAuthor{N.~Alshehri, D.~Boffi, and L.~Gastaldi}

\maketitle

\begin{abstract}
In this paper, new unfitted mixed finite elements are presented for elliptic interface problems with jump coefficients. Our model is based on a fictitious domain formulation with distributed Lagrange multiplier. The relevance of our investigations is better seen when applied to the framework of fluid-structure interaction problems. Two finite element schemes with piecewise constant Lagrange multiplier are proposed and their stability is proved theoretically. Numerical results compare the performance of those elements, confirming the theoretical proofs and verifying that the schemes converge with optimal rates. \def\keywordstitle{Keywords}

\end{abstract}
    
\section{Introduction}
\label{se:intro} 
Elliptic interface problems with jump coefficients are important and widely used  in applications including  bio-science and fluid-dynamics applications. 

We consider a problem where the coefficients in the governing partial
differential equation may jump across the interface that separates two or more
sub-domains. There are various possibilities for the decomposition of the
domain. In this paper, we consider the case where one subdomain is immersed
into another one. Other cases might be considered as
well~\cite{hansbo2005lagrange} and limiting to two subdomains doesn't affect
the generality of our discussion.

\begin{figure}[H]
\begin{minipage}[c]{.4445\linewidth}
\begin{center}
	\includegraphics[width=0.7\textwidth]{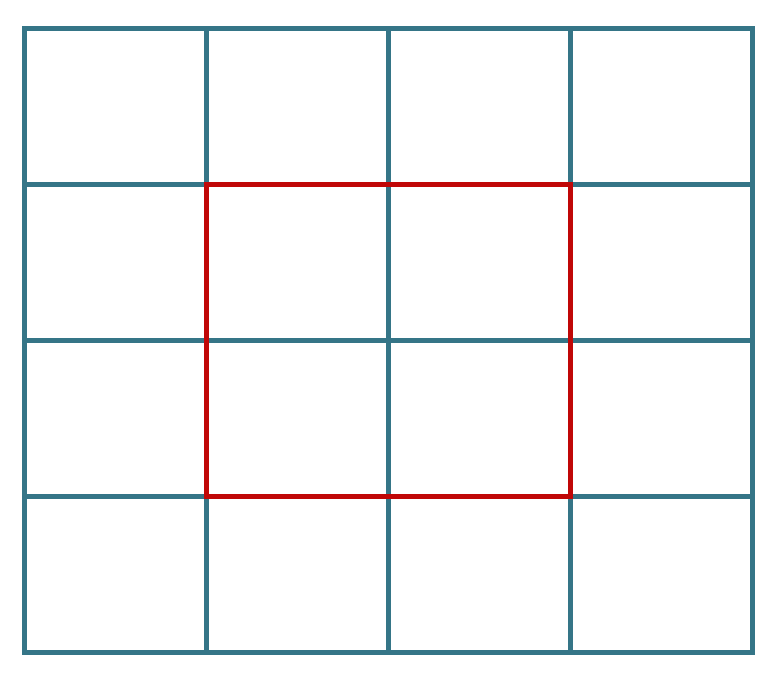} \\
	\small \text{A: Fitted boundaries.}
\end{center}
\end{minipage}
\begin{minipage}[c]{.4445\linewidth}
\begin{center}
	\includegraphics[width=0.7\textwidth]{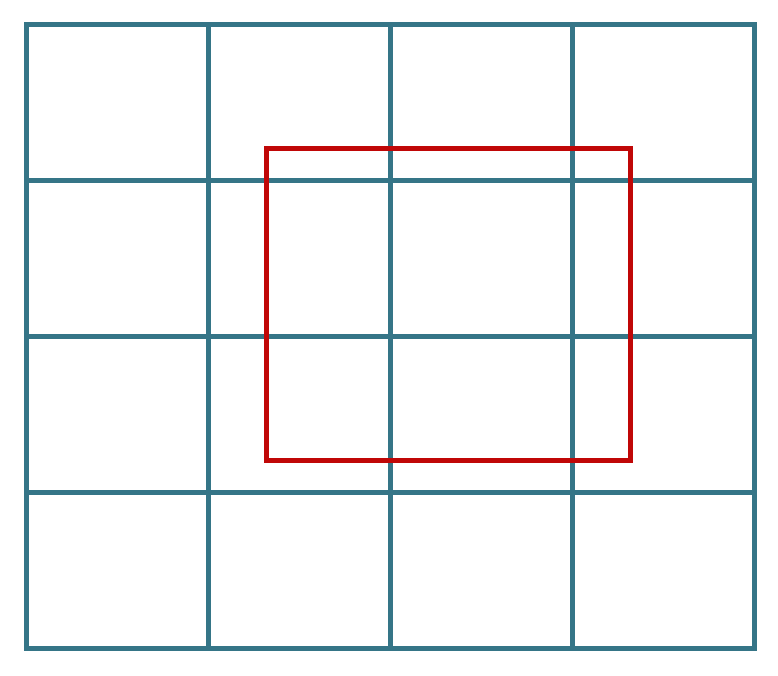} \\
	\small \text{B: Unfitted boundaries.}
\end{center}
\end{minipage}
	\caption{Fitted vs. unfitted boundaries.}
	\label{fig:fitted_unfitted}
\end{figure}

One way to address interface problems is to use fitted meshes. This approach involves generating a mesh that conforms to the interface between different materials or regions within the domain, as illustrated in Figure~\ref{fig:fitted_unfitted}.A. Through the years, several methods have been developed within this approach, including the extended finite element method using fixed-grids which is based on an 
extended Finite Element Method strategy~\cite{gerstenberger2008extended} and
multiscale finite element method~\cite{hou1999convergence}, as well as
others~\cite{barrett1987fitted},
\cite{chen2017interface}, \cite{huynh2013high}, and \cite{mu2016new}. One
common technique used to deal with interface problems within this category
involves utilizing an arbitrary Lagrangian–Eulerian (ALE) coordinate system
\cite{donea1977lagrangian}, \cite{hirt1974arbitrary}, and
\cite{hughes1981lagrangian}. While this method can provide
accurate solutions to such problems, it can be challenging,
particularly for time-dependent problems like fluid-structure interaction
problems. In these cases, the mesh must track the system evolution and
interface explicitly, which can result in ill-shaped meshes around the
interface due to large displacements and deformations. One way to address this
issue is to update the mesh at each time step, although this approach can be
computationally expensive and difficult to use.

An alternative approach is to use unfitted meshes. In this
approach, the meshes are independent of the interface, which is allowed to cut
through the interior of elements, as shown in
Figure~\ref{fig:fitted_unfitted}.B. A popular way to enforce the jump
conditions in the unfitted mesh approach is to modify the finite element basis
near the interface. Many numerical methods have been proposed in this
direction, such as the immersed boundary method, immersed interface
method,
immersed finite element methods, multiscale finite element methods, extended
finite element methods, and many others, see, e.g. \cite{leveque1994immersed},
\cite{Li-Ito}, \cite{li2003new}, \cite{chu2010new}, \cite{kwak2010analysis}, and
\cite{melenk1996partition}. The jump condition can also be imposed
by introducing penalty terms across interfaces~\cite{ babuvska1970finite}
or applying the Nitsche's method \cite{ nitsche1971variationsprinzip} as in the
cut finite element methods (CUT-FEM), see, e.g., \cite{massjung2012unfitted},
\cite{hansbo2002unfitted}, and  \cite{hansbo2014cut}.
This paper introduces the Fictitious Domain with Distributed Lagrange
Multiplier method (FD-DLM), which adapts the unfitted mesh approach. This
method is based on the Immersed Boundary Method (IBM), which was introduced by
Peskin in~\cite{peskin1972flow} for the simulation of blood flow in the heart
in the early seventies using finite differences. Later, in 2003, IBM was
established in the framework of finite elements in~\cite{boffi2003finite}. This
research led to the FD-DLM method, which avoids mesh regeneration by using
fixed meshes. It makes use of the fictitious domain approach introduced
in~\cite{glowinski1994fictitious}, \cite{yu2005dlm} to fictitiously extend one
mesh into the other. Then, the two meshes are considered independent of each
other and constructed only once. To enforce that the solution in the extended
domain coincides with the solution in the immersed domain, we add a coupling
term, and in our model, a Lagrange multiplier term is responsible for that. The
advantage of the FD-DLM method is that it avoids the need for re-meshing and
simplifies the mesh generation process, especially when the interface is moving
in time.

Our focus on studying the elliptic interface problem in this formulation goes beyond the development of an advanced method for approximating interface problems and finding a novel approach with optimal solutions for them. Rather, our research represents a theoretical exploration of this issue, which serves as a foundation for modeling fluid-structure interaction (FSI) problems where the interface is dynamic and possibly deforming over time. The FD-DLM formulation is advantageous in these cases. As a result, finding alternative combinations of finite element spaces that can provide a stable discretization of our problem is valuable in the context of using these same spaces for more intricate scenarios.

A crucial aspect of this method is how to deal with the coupling term.  Such a term is represented by a bilinear form defined on suitable Hilbert spaces. For instance, in~\cite{auricchio2015fictitious}, the discretization of this term was represented as the scalar product of $L^2$. In~\cite{boffi2014mixed} the
coupling term was represented as the scalar product either in $L^2$ or $H^1$.
Furthermore, the computation of this term involves the evaluation of the
integral on the immersed domain of shape functions that are supported on two
different meshes. In~\cite{boffi2022interface} it is discussed how to
implement this term in practice. It is shown that, in order to achieve the optimal convergence rate of the method, one has to perform the integral exactly by examining the intersection of the two meshes. Moreover, it is observed that finding the geometric intersection cannot be avoided even if the precision of the used quadrature rule is increased. Hence, in our numerical tests we will follow the intersection approach. More details on the coupling terms are given
in Section~\ref{se:Numerical_results}.

In~\cite{auricchio2015fictitious} continuous piecewise linear finite element
spaces are considered for the discretization of the problem on triangular
meshes. In~\cite{boffi2014mixed} continuous piecewise bilinear finite element
spaces are considered on quadrilateral meshes. Recently, \cite{boffiexistence}
showed, in the framework of FSI, the stability of a linearization of the
continuous problem and introduced a unified setting for the choice of the
finite element spaces.  This setting allows for more general choices of
spaces. However, so far only continuous finite element spaces are considered
for the multiplier responsible for the coupling term.

Our work is an extension to those papers where we are interested in finding stable elements with more flexibility in the choice of the multiplier. In particular, we addressed for the first time the question of whether piecewise discontinuous elements can be used for the approximation of the Lagrange multiplier. The motivation for our choice originates from FSI problems, where having a discontinuous Lagrange multiplier could improve the local mass conservation properties.

In Section~\ref{se:Model_problem} we introduce in detail the problem in the continuous setting. This problem is well-posed; we then consider its discretization and propose possible choices of elements that are the main object of our study in Section~\ref{se:discretization}.  Next, we prove the well-posedness of our discrete schemes in Section~\ref{se:A_Priori}.  The stability of our schemes is based on the presence of interior degrees of freedom (so-called bubble functions) in the space approximating the solution
where the Lagrange multiplier is distributed. We are showing numerically in Section~\ref{sse:weakinfsup} that the presence of the bubble functions is necessary for the discrete inf-sup condition.
Section~\ref{sse:rateofconv} reports a series of numerical tests which confirm the theoretical results.

\section{Formulation of the method}
\subsection{Model problem}
\label{se:Model_problem}
\begin{figure}[htp]
\centering
\begin{minipage}[c]{.2225\linewidth}
\begin{center}
	\includegraphics[width=0.9\textwidth]{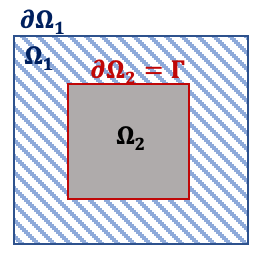} 
\end{center}
\end{minipage}
\begin{minipage}[c]{.1\linewidth}
\begin{center}
=
\end{center}
\end{minipage}
\begin{minipage}[c]{.2225\linewidth}
\begin{center}
	\includegraphics[width=0.9\textwidth]{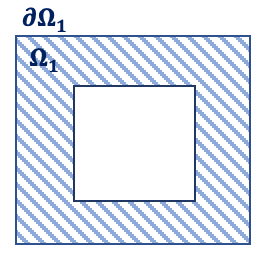} 
\end{center}
\end{minipage}
\begin{minipage}[c]{.1\linewidth}
\begin{center}
+
\end{center}
\end{minipage}
\begin{minipage}[c]{.2225\linewidth}
\begin{center}
	\includegraphics[width=0.9\textwidth]{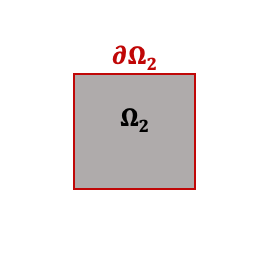} 
\end{center}
\end{minipage}
	\caption{ Domain decomposition in 2D.}
	\label{fig:domain_decomposition}
\end{figure}

Let $\Omega$ be a domain in $\mathbb{R}^d, d=1,2,3$ with a bounded Lipschitz boundary $\partial \Omega$. We assume that $\Omega$ is subdivided into two subdomains $\Omega_{i},i=1,2$ so that
$\overline\Omega=\overline\Omega_1\cup\overline\Omega_2$. The subdomains are separated by a Lipschitz continuous interface
$\Gamma=\overline\Omega_1\cap\overline\Omega_2$.
In order to simplify the presentation, we assume that $\Omega_2$ is immersed in
$\Omega$ so that $\overline\Gamma\cap\partial\Omega=\emptyset$. 
Figure~\ref{fig:domain_decomposition} reports a sketch of the situation in 2D.
Then, we consider the following problem:
\begin{problem}
\label{eq:cont_interface_pbm}
Given $f_1:\Omega_1 \to \mathbb{R}$, $f_2:\Omega_2 \to \mathbb{R}$, find
$u_1:\Omega_1 \to \mathbb{R}$ and $u_2:\Omega_2 \to \mathbb{R}$ such that:
\begin{subequations}
\begin{align}
-\nabla \cdot (\beta_i \nabla u_i)&=f_i
&& \mathrm{in }\; \Omega_i,\ i=1,2
\label{eq:cont_interface_pbm1}\\
u_1 &= u_2 && \mathrm{on }\;  \Gamma
\label{eq:cont_interface_pbm2}\\
\beta_1 \nabla(u_1) \cdot \textbf{n}_1 &=-\beta_2 \nabla(u_2) \cdot
\textbf{n}_2 
&& \mathrm{on }\;  \Gamma
\label{eq:cont_interface_pbm3}\\
u_1&=0 && \mathrm{ on } \; \partial \Omega_1
\label{eq:cont_interface_pbm4}
\end{align}
\end{subequations}
where $\textbf{n}_i$, $i=1,2$, is the unit vector pointing out of $ \Omega_i$ 
and normal to $\Gamma$. 
\end{problem}
This is an elliptic interface problem with a jump in the coefficients and
homogeneous Dirichlet boundary condition. 
We assume that the coefficients $\beta_i$ belong to $L^\infty(\Omega_i)$
($i=1,2$) and that they are bounded by below as follows:
\begin{align}
\label{eq:coeff}
\beta_1&>\underline\beta_1>0\\
\beta_2&>\underline\beta_2>0.
\end{align}
In 2D, this model typically describes the displacement of a membrane
made of two materials. The coefficients $\beta_i$ stand for the
stiffness of the materials, $f_i$ $i=1,2$ for the loads applied to the
membrane, and $u_i$ for the vertical displacement in $\Omega_i$ ($i=1,2$), respectively. Equation~\eqref{eq:cont_interface_pbm2} guarantees the continuity of the solutions $u_1$ and $u_2$ on the interface $\Gamma$. This means that we are considering connected materials that do not break. Moreover, Equation~\eqref{eq:cont_interface_pbm3} prescribes a jump of the
normal derivatives of $u_1$ and $u_2$ at the interface that is inversely proportional to the ratio of the coefficients. 

\begin{figure}[htp]
\centering
\begin{minipage}[c]{.2225\linewidth}
\begin{center}
	\includegraphics[width=0.9\textwidth]{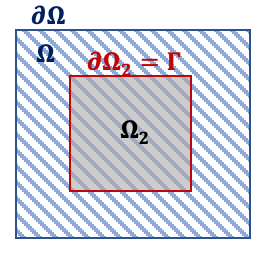} 
\end{center}
\end{minipage}
\begin{minipage}[c]{.1\linewidth}
\begin{center}
=
\end{center}
\end{minipage}
\begin{minipage}[c]{.2225\linewidth}
\begin{center}
	\includegraphics[width=0.9\textwidth]{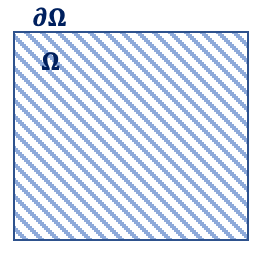} 
\end{center}
\end{minipage}
\begin{minipage}[c]{.1\linewidth}
\begin{center}
+
\end{center}
\end{minipage}
\begin{minipage}[c]{.2225\linewidth}
\begin{center}
	\includegraphics[width=0.9\textwidth]{omega2.png} 
\end{center}
\end{minipage}
	\caption{ $\Omega_1$ fictitiously extended in $\Omega_2$  in 2D.}
	\label{fig:ficticiousdomain}
\end{figure}

In this paper, we consider a fictitious domain approach, therefore
we reformulate Problem~\ref{eq:cont_interface_pbm}
following~\cite{auricchio2015fictitious}~and~\cite{boffi2014mixed}.
More precisely, we extend $u_1$, $\beta_1$, and $f_1$ to $\Omega$ and denote such extensions by $u$, $\beta$, and $f$, respectively, so that $u|_{\Omega_1}=u_1$, $f|_{\Omega_1}=f_1$, and
$\beta|_{\Omega_1}=\beta_1$. Moreover, we enforce the extended solution $u$ to coincide with $u_2$ in $\Omega_2$, i.e. $u|_{\Omega_2}=u_2$ by introducing a 
Lagrange multiplier. The resulting formulation will be called FD-DLM.

In view of the introduction of the variational formulation of the problem, we recall some notations. For any open connected domain, or part of a domain, $\omega\subset\mathbb{R}^d$ for $d=2,3$, we denote the standard Lebesgue and Sobolev spaces by $L^2(\omega)$ and $H^1(\omega)$, respectively. Those spaces are endowed with their norms;
$\norm{\cdot}_{L^2(\omega)}=\norm{\cdot}_{0,\omega}$ and
$\norm{\cdot}_{H^1(\omega)}=\norm{\cdot}_{1,\omega}$. Moreover, 
$(\cdot , \cdot)_{\omega}$ stands for the scalar product in $L^2(\omega)$.

Let us consider the spaces: 
\begin{align*}
V&=H^1_0(\Omega)=\left\lbrace v \in H^1(\Omega) |~v=0 ~\mathrm{on~} \partial \Omega \right\rbrace \\
V_2&=H^1(\Omega_2)
\end{align*}
endowed with their natural norms $\norm{v}_{V}=|v|_{1,\Omega}=\|\nabla v\|_{0,\Omega}$ and $\norm{v_2}_{V_2} =\norm{v_2}_{1,\Omega_2}$, respectively.

We denote by $\Lambda$ the dual space of $V_2$, i.e.
$\Lambda=\left[ H^1(\Omega_2)\right] ^*$, endowed with the following dual norm:

\[
\norm{\mu}_\Lambda = \sup_{v_2 \in V_2} \dfrac{\left\langle \mu,v_2\right\rangle }{\norm{v_2}_{V_2}}
\]
where  $\left\langle \cdot  , \cdot \right\rangle$  is the duality paring
between $V_2$ and its dual space $\Lambda$. 

In \cite{auricchio2015fictitious} it has been proved that
Problem~\ref{eq:cont_interface_pbm} is equivalent to the following fictitious
domain formulation with distributed Lagrange multiplier.
\begin{problem}
\label{pbm:1}
Given $f \in L^2(\Omega)$, $f_2 \in L^2(\Omega_2)$, $\beta\in L^\infty(\Omega)$
and $\beta_2\in L^\infty(\Omega_2)$ with $f|_{\Omega_1}=f_1$ and 
$\beta|_{\Omega_1}=\beta_1$, find
$(u,u_2,\lambda)\in V \times V_2 \times \Lambda$ such that
\begin{align*}
(\beta \nabla u, \nabla v)_{\Omega} + \left\langle \lambda ,v|_{\Omega_2} \right\rangle  
&=(f,v)_{\Omega} && \forall v \in V\\
( (\beta_2-\beta )\nabla u_2, \nabla v_2)_{\Omega_2} - \left\langle \lambda ,v_2 \right\rangle 
&=(f_2-f,v_2)_{\Omega_2} && \forall v_2 \in V_2\\
\left\langle \mu , u|_{\Omega_2} -u_2 \right\rangle
&=0 && \forall \mu \in \Lambda.
\end{align*}
\end{problem}
After some standard calculations, one can obtain the following characterization of
$\lambda$, that will be useful to estimate the approximation error,
see~\cite{auricchio2015fictitious}: 
\begin{equation}
\label{eq:dual_lambda}
\left\langle \lambda, v_2 \right\rangle = -
\int_{\Omega_2}\left(\frac{\beta}{\beta_2} f_2-f \right)\; v_2\;dx
+\int_{\Gamma} (\beta_2 -\beta) \nabla u_2 \cdot \textbf{n}_2\; v_2\; d \gamma.
\end{equation}
Clearly, Problem~\ref{pbm:1} is a saddle point problem that can be written in operator form as follow:
\begin{equation}
\label{eq:full_matrix}
\begin{pmatrix}
A_1 & 0&\rvline &C_1^T\\
0&A_2&\rvline &-C_2^T\\
\hline
C_1&-C_2&\rvline &0\\
\end{pmatrix}\\
\begin{pmatrix}
u\\
u_2\\
\hline
\lambda\\
\end{pmatrix}\\
 = 
\begin{pmatrix}
F_1\\
F_2\\
\hline
0\\
\end{pmatrix},
\end{equation}
where $A_1$ and $A_2$ are the operators associated with the bilinear forms
$( \beta \nabla u, \nabla v)_{\Omega}$ and 
$( (\beta_2-\beta )\nabla u_2, \nabla v_2)_{\Omega_2}$, respectively.
Moreover, $(C_1,\, -C_2)$ is the operator pair that is associated with the
bilinear form $\left\langle \mu , u|_{\Omega_2} -u_2 \right\rangle$
with kernel:
\[
\mathbb{K}=\left\lbrace (u,u_2)\in V \times V_2 :
\left\langle \mu , u|_{\Omega_2}-u_2\right\rangle =0 , 
\forall\mu\in\Lambda\right\rbrace.
\]
Notice that due to the definition of $\Lambda$ we have that $u|_{\Omega_2}=u_2$
for all $(u,u_2)\in\mathbb{K}$.

Typically, in order to prove the well-posedness of a continuous saddle
point problem like Problem~\ref{pbm:1}, one needs to verify the following two
sufficient conditions (see~\cite{boffi2013mixed}, \cite{babuvska1973finite}, \cite{brezzi1974existence}).
\begin{itemize}
\item \textbf{Continuous elker condition:} 
There exists a constant $ \overline{\gamma_1}>0$ such that for all $(v,v_2)$ in
$ \mathbb{K} $ the following inequality holds true.
\[
(\beta \; \nabla v, \nabla v)_{\Omega} +((\beta_2 - \beta)\; \nabla
v_{2},\nabla v_{2})_{\Omega_2} \geq \overline{\gamma_1} \left( \norm{v}_{V}^2
+\norm{v_{2}}_{V_2}^2\right).
\]
\item \textbf{Continuous inf-sup condition: } 
There exists a constant $ \overline{\gamma_2}>0$, such that
for all $ \mu\in\Lambda$ the following bound holds true.
\[
\sup_{(v,v_2)\in V \times V_2} \dfrac{\left\langle \mu,
v|_{\Omega_2}-v_{2}\right\rangle }{\left( \norm{v}_{V}^2+\norm{v_{2}}_{V_2}^2
\right)^\frac{1}{2}} \geq\overline{\gamma_2}  \norm{\mu}_{\Lambda}.
\]
\end{itemize}
Both conditions were shown to hold in~\cite{auricchio2015fictitious}.
Therefore, a unique stable solution exists for Problem~\ref{pbm:1} and we can
state the following proposition (see~\cite{boffi2013mixed}, \cite{babuvska1973finite}, \cite{brezzi1974existence}).\\

\begin{proposition}[\textbf{Stability}]  \label{prop:a_priori_estimate}
Given $f \in L^2(\Omega)$ and $f_2 \in L^2(\Omega_2)$, there exists a unique
solution $(u,u_{2},\lambda)$ in $ V \times V_{2} \times \Lambda$ for
Problem~\ref{pbm:1}  which satisfies the following a priori estimate:
\[
\norm{u}_{V}+\norm{u_{2}}_{V_2}+\norm{\lambda}_{\Lambda} \leq C \left(
\norm{f}_{0,\Omega}+\norm{f_2}_{0,\Omega_2}\right).
\]
\end{proposition}

\begin{remark}
\label{rmk:regularity}
It is known that the regularity of the solution $(u_1,u_2)$ of an elliptic interface problem with discontinuous coefficients and a Lipschitz interface $\Gamma$, such as Problem~\ref{eq:cont_interface_pbm}, is limited by the presence of re-entrant corners of the interface and of the external boundary. Hence, we have that there exists $s$ with $\frac32<s\le2$ such that
$u_i\in H^{s}(\Omega_i)$ for $i=1,2$, (see~\cite{nicaise}).
Moreover, since the solution $u\in H^1_0(\Omega)$ of Problem~\ref{pbm:1} exhibits jump in the derivative normal to the interface, it belongs to
$H^r(\Omega)$, with $1 <r<\frac{3}{2}$. 
\end{remark}
\subsection{Finite element discretization}
\label{se:discretization}
Let $\mathcal{T}$ and $\mathcal{T}_2$ be two shape-regular meshes of the fictitiously extended domain $\Omega$ and the immersed domain $\Omega_2$,
respectively. Here, we are considering quadrilateral meshes in 2D and
hexahedral in 3D, with $h$ and $h_2$ denoting the maximum mesh size of
$\mathcal{T}$ and $\mathcal{T}_2$, respectively. We introduce the 
finite element spaces $V_h \subset V$, $V_{2h} \subset V_2$, and $\Lambda_h \subset \Lambda$. $V_h$ and $V_{2h}$ contain piecewise polynomials continuous across the interelement boundaries, while for $\Lambda_h$ we choose discontinuous finite elements. However, since $\Lambda_h\subset L^2(\Omega_2)$,
the duality paring in Problem~\ref{pbm:1} can be evaluated using the $L^2$ scalar product in $\Omega_2$. Then, the discrete version of Problem~\ref{pbm:1} reads:
\begin{problem}
\label{pbm:2}
Given $f \in L^2(\Omega)$ and $f_2 \in L^2(\Omega_2)$, find $(u_h,u_{2h},\lambda_h) \in V_h \times V_{2h} \times \Lambda_h$ such that
\begin{align*}
(\beta \nabla u_h, \nabla v_h)_{\Omega} + (\lambda_h ,v_h |_{\Omega_2} )_{\Omega_2}  &=(f,v_h)_{\Omega}  && \forall v_h \in V_h\nonumber \\
((\beta_2-\beta )\nabla u_{2h}, \nabla v_{2h})_{\Omega_2}  - ( \lambda_h ,v_{2h} )_{\Omega_2} & =(f_2-f,v_{2h})_{\Omega_2}  && \forall v_{2h} \in V_{2h}\\
( \mu_h,u_h |_{\Omega_2} -u_{2h} )_{\Omega_2}&=0&&\forall \mu_h \in\Lambda_h.
\nonumber
\end{align*}
\end{problem}
For $K$ an element of $\mathcal{T}$ or $\mathcal{T}_2$, we define $Q_k(K)$, $k\geq 1$, to be the space of finite elements that are polynomials of degree at most $k$, separately in each variable on $K$. If $k=0$ then $Q_0(K)$ is the space of constant polynomials and will be denoted by $P_0(K)$. Moreover, let $B(K) \in Q_2(K)$ be a bubble function defined on $K$ and vanishing at the boundary of the element $K$. This function will be used in this work to enrich the space $Q_1(K)$.

In the following, we are going to discretize the Lagrange multiplier by piecewise constants, hence we introduce the following natural choices for $V_h\times V_{2h}\times\Lambda_h$.

\begin{figure}[H]
\centering
	\includegraphics[width=0.25\textwidth]{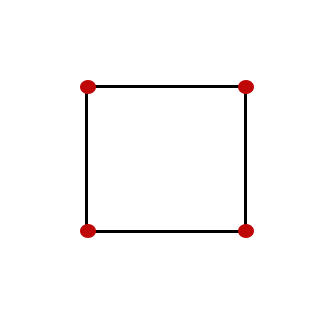} 
	\includegraphics[width=0.25\textwidth]{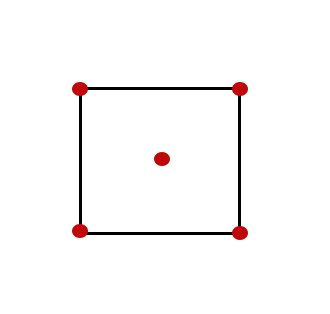} 
	\includegraphics[width=0.25\textwidth]{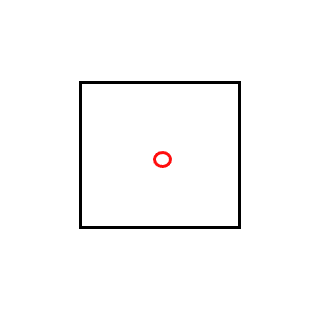} 
	\caption{Element 1 in 2D , $Q_1-(Q_1+B)-P_0$.}
	\label{fig:Q1Q1BP0}
\end{figure}

\begin{figure}[H]
\centering
	\includegraphics[width=0.25\textwidth]{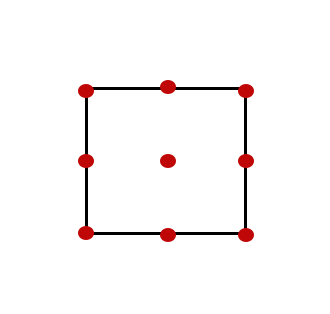} 
	\includegraphics[width=0.25\textwidth]{q2.png} 
	\includegraphics[width=0.25\textwidth]{p0.png} 
	\caption{Element 2 in 2D, $Q_2-Q_2-P_0$.}
	\label{fig:Q2Q2P0}
\end{figure}

\begin{itemize}
\item \textbf{Element 1:} $Q_1-\left(Q_1+B\right) -P_0$\\
In this element we define the discrete subspaces $V_h,V_{2h},\Lambda_h$  to be 
\begin{align}
\label{eq:elm1}
V_h&=\left\lbrace v_h \in V: v_h |_{K} \in \mathit{Q}_1(K),
\forall K \in \mathcal{T}\right\rbrace\nonumber\\
V_{2h}&=\left\lbrace v_{2h} \in V_2: v_{2h} |_{K}\in\mathit{Q}_1(K)+B(K),
\forall K \in \mathcal{T}_2 \right\rbrace\\
\Lambda_h&=\left\lbrace \mu_h \in \Lambda :  \mu_h |_{K} \in
\mathit{P}_0(K) , \forall K\in \mathcal{T}_2 \right\rbrace.
\nonumber
\end{align}
Hence, the solution $u_{h}$ is approximated by piecewise bilinear polynomials and $u_{2h}$ by piecewise bilinear polynomials enriched by bubble functions, so that, five degrees of freedom are used per element in the space $V_{2h}$ (see Figure~\ref{fig:Q1Q1BP0}). 
\item \textbf{Element 2: $Q_2 - Q_2 - P_0$}\\
Here, we define the discrete subspaces $V_h,V_{2h},\Lambda_h$  to be 
\begin{align}\label{eq:elm2}
V_h&=\left\lbrace v_h \in V: v_h |_{K} \in \mathit{Q}_2(K), \forall K\in
\mathcal{T}\right\rbrace 
\nonumber\\
V_{2h}&=\left\lbrace v_{2h} \in V_2: v_{2h} |_{K} \in \mathit{Q}_2(K),
\forall K \in \mathcal{T}_2 \right\rbrace\\
\Lambda_h&=\left\lbrace \mu_h \in \Lambda :  \mu_h |_{K} \in
\mathit{P}_0(K) , \forall K \in \mathcal{T}_2 \right\rbrace .
\nonumber
\end{align}
where we approximate both $u_h,u_{2h}$ by piecewise biquadratic polynomials.
Hence, we use nine nodes for each element in the spaces $V_h$ and $V_{2h}$ (see
Figure~\ref{fig:Q2Q2P0}). It is clear that compared to Element 1, $V_{2h}$
uses four extra nodes at the middle of each edge, while $V_{h}$ requires five
extra degrees of freedom, so this element is computationally more
expensive than Element 1.
\end{itemize}

\section{Error estimate}
\label{se:A_Priori}
In this section, we study the existence of a unique stable solution of
Problem~\ref{pbm:2}. Since it is
a discrete saddle point problem, sufficient conditions for its well-posedness 
are the discrete elker and inf-sup conditions (see~\cite{boffi2013mixed}, \cite{babuvska1973finite}, \cite{brezzi1974existence}). \\
Let the discrete kernel $\mathbb{K}_h$ associated with the bilinear form
$( \mu_h , u_h|_{\Omega_2} -u_{2h} )_{\Omega_2}$ be defined as follow:
\[
\mathbb{K}_h=\left\lbrace (u_h,u_{2h})\in V_h \times V_{2h} : ( \mu_h , u_h|_{\Omega_2}-u_{2h})_{\Omega_2}  =0 , \forall \mu_h \in \Lambda_h\right\rbrace.
\]
Then, the following proposition states the discrete elker condition.
\begin{proposition}[\textbf{Discrete elker condition}] \label{prop:discreteelker}
Let us consider $V_h,V_{2h},\Lambda_h$ defined in \eqref{eq:elm1} and
\eqref{eq:elm2} and assume that the coefficients
satisfy $\beta_2 > \beta >\underline\beta>0$ in $\Omega_2$, then there exists a
constant $\gamma_1>0$ independent of the discretization parameters $h,h_2$,
such that for all $ (v_h ,v_{2h}) \in \mathbb{K}_h$ the following inequality
holds true.
\[
\left(\beta\nabla v_h, \nabla v_h\right)_{\Omega}
+\left((\beta_2 - \beta)\nabla v_{2h},\nabla v_{2h}\right)_{\Omega_2} 
\geq \gamma_1 \Bigg[ \norm{v_h}_{V}^2+\norm{v_{2h}}_{V_2}^2\Bigg].
\]
\end{proposition}
\begin{proof}
For all $ (v_h ,v_{2h}) \in \mathbb{K}_h$ we have:
\begin{align*}
(\beta\nabla v_h, \nabla v_h)_{\Omega}
+((\beta_2 - \beta)\nabla v_{2h},\nabla v_{2h})_{\Omega_2}
&=\beta \norm{\nabla v_h}_{0,\Omega}^2
+(\beta_2 - \beta) \norm{\nabla v_{2h}}^2_{0,\Omega_2}\\
&\geq \underline\beta \norm{\nabla v_h}_{0,\Omega}^2
+\eta_0 \norm{\nabla v_{2h}}^2_{0,\Omega_2},
\end{align*}
where $\eta_0$ is such that $\beta_2-\beta\ge\eta_0>0$.
\\
Since $V=H^1_0(\Omega)$, we can apply the Poincaré inequality
$\|v_h\|_{0,\Omega}\le C_\Omega\|\nabla v_h\|_{0,\Omega}$, and we have that:
\[
(\beta\nabla v_h,\nabla
v_h)_{\Omega}\ge\frac{\underline\beta}2\left(1+\frac1{C_\Omega^2}\right) 
\norm{ v_h}_{1,\Omega}^2.
\]
It remains to bound by below $\norm{\nabla v_{2h}}_{0,\Omega_2}$ by means of $\norm{v_{2h}}_{1,\Omega_2}$. In order to use the Poincaré–Wirtinger inequality, we split $v_{2h}$ as follows
\[
v_{2h}=E(v_{2h}) + \bigg[ v_{2h}-E(v_{2h})\bigg]
\]
where $E(v_{2h})\in\mathbb{R}$ is the mean value of $v_{2h}$ . Then,
\begin{align*}
\norm{v_{2h}}_{0,\Omega_2} &\leq \norm{E(v_{2h}) }_{0,\Omega_2}
+\norm{v_{2h}-E(v_{2h})}_{0,\Omega_2} \\
&\leq \norm{E(v_{2h}) }_{0,\Omega_2} +C\norm{\nabla v_{2h}}_{0,\Omega_2}. 
\end{align*}
Now, since $(v_{h}, v_{2h})\in \mathbb{K}_h$, and $\Lambda_h$ contains
constant functions on $\Omega_2$, we have 
\begin{equation}
\label{eq:mean_value_bound}
\aligned
( E(v_{2h}), E(v_{2h}))_{\Omega_2}&
=(E(v_{2h}),v_{2h})_{\Omega_2}-(E(v_{2h}), v_{2h}-E(v_{2h}))_{\Omega_2}\\ 
&=( E(v_{2h}) , v_h|_{\Omega_2})_{\Omega_2}-
(E(v_{2h}), v_{2h}-E(v_{2h}))_{\Omega_2}\\
&=( E(v_{2h}) , v_h|_{\Omega_2})_{\Omega_2}.
\endaligned
\end{equation}
Then, \eqref{eq:mean_value_bound} gives the bound
\begin{align*}
\norm{E(v_{2h})}_{0,\Omega_2}^2 &\leq \norm{E(v_{2h})}_{0,\Omega_2} \norm{v_h}_{0,\Omega_2}.
\end{align*}
Hence,
$$\norm{E(v_{2h})}_{0,\Omega_2} \leq \norm{v_h}_{0,\Omega_2}$$
which concludes the proof.

\end{proof}

\begin{remark}
In a particular case and under the additional condition that $h_2/h^d$ is
sufficiently small, in~\cite{boffi2014mixed} it has been shown that the
discrete elker condition holds true even if the constraint $\beta_2>\beta$ is
not verified. In that case triangular linear elements were used for the
approximation of both $V_2$ and $\Lambda$.
\end{remark}

In the next proposition, we prove the discrete inf-sup condition. 

\begin{proposition}[\textbf{Discrete inf-sup condition}] Let
$V_h,V_{2h},\Lambda_h$ be the spaces defined in~\eqref{eq:elm1} and~\eqref{eq:elm2} for Element~1 and Element~2, respectively. Then there
exists a positive constant $\gamma_2>0$
independent of the discretization parameters $h,h_2$, such that the following
discrete inf-sup condition holds.\\
\begin{equation}
\label{eq:Discreteinfsup}
\sup_{(v_h,v_{2h})\in V_h \times V_{2h}}
\dfrac{( \mu_h, v_h|_{\Omega_2}-v_{2h})_{\Omega_2}  }{\left(
\norm{v_h}_{V}^2+\norm{v_{2h}}_{V_2}^2 \right)^\frac{1}{2}} \geq \gamma_2
\norm{\mu_h}_{\Lambda}.
\end{equation}
\label{prop:discreteinfsup}
\end{proposition}

\begin{proof}
We prove this proposition for Element~1. The same proof carries on for
Element~2.

Since $v_h=0$ is a possible choice in $V_h$, we have that
\[
\sup_{(v_h,v_{2h})\in V_h \times V_{2h}} 
\dfrac{( \mu_h, v_h|_{\Omega_2}-v_{2h})_{\Omega_2}  }
{\left(\norm{v_h}_{V}^2+\norm{v_{2h}}_{V_2}^2 \right)^\frac{1}{2}} \geq
\sup_{v_{2h}\in V_{2h}} \dfrac{( \mu_h, v_{2h})_{\Omega_2}
}{\norm{v_{2h}}_{V_2}}.
\]
Hence, it is enough to prove the following bound:
\begin{equation}
\label{eq:simplierDiscreteinfsup}
\sup_{v_{2h}\in V_{2h}} \dfrac{( \mu_h, v_{2h})_{\Omega_2} }
{\norm{v_{2h}}_{V_2}} \geq \gamma_2 \norm{\mu_h}_{\Lambda}.
\end{equation}

In order to show that \eqref{eq:simplierDiscreteinfsup} is satisfied, we use a
Fortin trick (see~=\cite[Prop.~5.4.2]{boffi2013mixed}).
Due to the fact that the continuous inf-sup condition is satisfied, the aim of this proof is to find a linear Fortin operator $\Pi_h:V_2 \to V_{2h}$ that satisfies the following relations for all $v_2 \in V_{2}$:
\begin{equation}
\label{eq:Fortin}
\aligned
&\text{ i) } ||\Pi_h v_2||_{V_2} \le C ||v_2||_{V_2}.\\
&\text{ii) } ( \mu_h, v_{2})_{\Omega_2} =( \mu_h, \Pi_h v_{2})_{\Omega_2}& \forall \mu_h\in\Lambda_h.
\endaligned
\end{equation}
We introduce the following subspaces of $V_{2h}$:\\
$$
\overline{\overline{V_{2h}}}=\left\lbrace v_{2h} \in V_2: v_{2h} |_{K} \in \mathit{Q}_1(K), \forall K \in \mathcal{T}_2 \right\rbrace \subset V_{2h}
$$
$$
\overline{V_{2h}}=\left\lbrace v_{2h} \in V_2: v_{2h} |_{K} \in B(K) , \forall K \in \mathcal{T}_2 \right\rbrace \subset V_{2h}.
$$
Let $\Pi_h=\Pi_1+\Pi_2(I-\Pi_1)$ where $\Pi_1:V_2 \to\overline{\overline{V_{2h}}} \subset V_{2h}$ is the Cl\'ement's
operator, such that for all $v_2\in V_2$\\
\begin{subequations}
\begin{align}
&\sum_{K \in \mathcal{T}_2} 
h_{K}^{-2}\norm{v_2-\Pi_1 v_2}^2_{0,K}
\leq C \norm{v_2}_{V_2}^2\label{eq:boundspi1}\\
&\sum_{K \in \mathcal{T}_2} 
\norm{v_2-\Pi_1 v_2}^2_{1,K}
\leq C \norm{v_2}_{V_2}^2,\label{eq:boundspi2}
\end{align}
\end{subequations}
and $\Pi_2:V_2 \to \overline{V_{2h}} \subset V_{2h}$ such that $\Pi_2(p_i)=0 $ for all  nodes $p_i$ at vertices of the element $K \in \mathcal{T}_2$ and  \begin{equation}\label{eq:pi2}
\int_{K} \Pi_2 v_2 =\int_{K} v_2 \hspace{1in} \forall v_2 \in V_2 \text{ and }
K \in \mathcal{T}_2.
\end{equation}
Let $\mathcal{F}_{K}^{-1}$ be the affine mapping that maps objects from the current element, $K\in \mathcal{T}_2$, to the reference element denoted by $\hat{K}$, that is the unit square in 2D and the unit cube in 3D. Then, $\hat{v}_2=v_2\circ\mathcal{F}_{K}$ where symbols with hat refer to quantities evaluated in the reference domain $\hat{K}$.\\ 

Let $\abs{K}$ be the measure of the element $K$. Then,
\[
\int_{\hat{K}} \widehat{\Pi_2 v_2}=\abs{K}^{-1} \int_{K} \Pi_2 v_2
=\abs{K}^{-1} \int_{K} v_2 
=\abs{K}^{-1} \abs{K} \int_{\hat{K}}\hat{v}_2 =\int_{\hat{K}} \hat{v}_2
\]
which means that $\widehat{\Pi_2}$ also satisfies equation \eqref{eq:pi2}. Moreover, let us define the following:
\begin{align*}
	\abs{\norm{\widehat{\Pi_2 v_2}}}=\abs{\int_{\hat{K}}\widehat{\Pi_2 v_2}}.
\end{align*}
It is clear that this is a norm in this specific case, since the function $\widehat{\Pi_2 v_2} $ by definition is the bubble in the reference element.
Hence, $ \abs{\int_{\hat{K}}\widehat{\Pi_2 v_2}}=0$ if and only if
$\widehat{\Pi_2 v_2}=0$.\\
Using the definition of $\Pi_h$ and the equality in \eqref{eq:pi2}, we get\\
\[
( \mu_h,v_2-\Pi_h v_2)_{\Omega_2}
=( \mu_h,(v_2-\Pi_1 v_2)-\Pi_2(v_2-\Pi_1 v_2))_{\Omega_2} =0
\quad \forall\mu_h\in\Lambda_h
\]
which proves that the proposed operator $\Pi_h$ satisfies condition ii)
in~\eqref{eq:Fortin}.

Now, we need to show condition i) for all $v_2 \in V_2$, that is:
\begin{equation}
\label{eq:condition1}
\norm{\Pi_h v_2}_{V_2}=\norm{\Pi_1 v_2 + \Pi_2 (v_2-\Pi_1 v_2)}_{V_2}
\leq C \norm{ v_2}_{V_2}.
\end{equation}
It is clear that \eqref{eq:boundspi2} implies
\begin{equation}
\label{eq:satbility1}
 \norm{\Pi_1v_2}_{V_2} \leq C \norm{v_2}_{V_2}.
 \end{equation}
For the sake of simplicity, let $w=v_2-\Pi_1 v_2$ and write the $H^1$-norm as follows:
\begin{equation}\label{eq:norm1}
\norm{\Pi_2 w}^2_{1,\Omega_2}=\norm{\Pi_2 w}^2_{0,\Omega_2} +\abs{\Pi_2 w}^2_{1,\Omega_2}.
\end{equation}

For all $K\in\mathcal{T}_2$, we estimate the first term in~\eqref{eq:norm1} as
follows:
\begin{align*}
	\norm{\Pi_2 w}_{0,K} 
			&\leq C \;  \abs{K}^{\frac{1}{2}} \; 
				 \norm{ \widehat{\Pi_2 w}}_{0,\hat{K}}
				&& \textit{mapping to the reference element $\hat{K}
					$}\\
			&\leq C \; \abs{K}^{\frac{1}{2}} \;  \abs{ \int_{\hat{K}}
				 \widehat{\Pi_2 w}}
				&& \textit{by equivalence of norms in finite dimensions}\\
			&= C \;  \abs{K}^{\frac{1}{2}} \; \abs{\int_{\hat{K}} \hat{ w}}
				&& \textit{since $\widehat{\Pi_2}$ satisfies
					 \eqref{eq:pi2}}\\
			&\leq  C \;  \abs{K}^{\frac{1}{2}}\;  \norm{\hat{w}}_{0,\hat{K}}
				&& \textit{by Cauchy-Schwarz inequality}\\
		&\leq C\;\abs{K}^{\frac{1}{2}}\;\abs{K}^{-\frac{1}{2}}\;\norm{w}_{0,K}
				&&\textit{mapping back  to the physical element $K$}\\
			&= C \;  \norm{w}_{0,K}.
\end{align*}
Similarly, using the same argument as before we have 
\begin{align*}
	\norm{\nabla\Pi_2 w}_{0,K}
			&\leq C \;  h_{K}^{-1} \;  \abs{K}^{\frac{1}{2}} \; 
				\abs{\widehat{\Pi_2 w}}_{1,\hat{K}}
				\leq C \;  h_{K}^{-1} \;  \abs{K}^{\frac{1}{2}} \;   
				\abs{ \int_{\hat{K}} \widehat{\Pi_2 w}}\\
			&\leq C \;  h_{K}^{-1} \;  \abs{K}^{\frac{1}{2}} \; 
				 \abs{ \int_{\hat{K}} \hat{ w}}
			\leq  C \;  h_{K}^{-1} \;  \abs{K}^{\frac{1}{2}} \; 
				\norm{\hat{w}}_{0,\hat{K}}.
\end{align*}
Mapping back to the current element $K$, we obtain:\\
\begin{equation}
\label{eq:seminorm1T}
\norm{\nabla\Pi_2 w}_{0,K} 
\leq C\; h_{K}^{-1}\;\abs{K}^{\frac{1}{2}}\;\abs{K}^{-\frac{1}{2}}\;
\norm{w}_{0,K} 
= C  \;  h_{K}^{-1} \;\norm{w}_{0,K}.
\end{equation}
By substituting \eqref{eq:seminorm1T} in \eqref{eq:norm1} and using the bounds of $\Pi_1$ in \eqref{eq:boundspi1}, we get
\begin{equation}\label{eq:satbility2}
\norm{\Pi_2 (v_2-\Pi_1 v_2)}^2_{V_2}\leq
C \sum_{K \in \mathcal{T}_2} 
h_{K}^{-2} \norm{v_2-\Pi_1 v_2}_{0,K}^2
 \leq C \norm{v_2}_{V_2}^2.
\end{equation}
Using the triangle inequality, and applying the results
of~\eqref{eq:satbility1} and~\eqref{eq:satbility2}, yield that condition i)
holds, in fact
\begin{align*}
\norm{\Pi_h v_2}_{V_2} & = \norm{\Pi_1 v_2 + \Pi_2 (v_2-\Pi_1 v_2)}_{V_2}\\
		   & \leq \norm{\Pi_1 v_2}_{V_2} + \norm{\Pi_2 (v_2-\Pi_1 v_2)}_{V_2}\\
				   &\leq C \norm{ v_2}_{V_2}.
\end{align*}
This concludes the proof and shows that the $H^1-$ stability of the constructed Fortin operator $\Pi_h$  is satisfied. Therefore, the discrete inf-sup condition \eqref{eq:Discreteinfsup} holds.
\end{proof}

\begin{remark}
Proposition~\ref{prop:discreteinfsup} shows rigorously that the inf-sup
condition for Element~1 and Element~2 is satisfied uniformly.
The bubble function in the space $V_{2h}$ is necessary for the stability of the element. In fact, in Section~\ref{sse:weakinfsup} we give numerical evidence that we do not have a uniform inf-sup bound if we modify Element~1 and remove the bubble to become $Q_1-Q_1-P_  0$.
\end{remark}

Thanks to Propositions~\ref{prop:discreteelker} and~\ref{prop:discreteinfsup}
there exists a unique stable solution $(u_h,u_{2h},\lambda_h)$ for
Problem~\ref{pbm:2} in the spaces $ V_h \times V_{2h} \times \Lambda_h$ that are defined in~\eqref{eq:elm1} and~\eqref{eq:elm2}. Moreover, by taking into consideration the regularity results of the solutions $u$ and $u_2$ recalled in Remark~\ref{rmk:regularity}, we can state the following a priori error estimate.
\\

\begin{proposition}[\textbf{Error estimate}]
Given  $(f,f_2) \in L^2(\Omega) \times L^2(\Omega_2)$.  Let $(u,u_{2},\lambda)
\in  V \times V_2 \times \Lambda$ and
$(u_h,u_{2h},\lambda_h) \in V_h \times V_{2h} \times \Lambda_h$ 
be the solutions of Problem~\ref{pbm:1} and~Problem~\ref{pbm:2}, respectively.
We consider Element~1 and Element~2 defined in \eqref{eq:elm1} and \eqref{eq:elm2}. Then,
under the assumption that the mesh $\mathcal{T}_2$ is quasi-uniform,
the following error estimate holds true:
\[ \aligned
\norm{u-u_h}_{V}&+\norm{u_2-u_{2h}}_{V_2}+\norm{\lambda-\lambda_h}_{\Lambda}\\
&\leq C\left(h^{r-1} \norm{u}_{r,\Omega} 
+\max(h_2^{s-1},h_2^{1-t})\norm{u_2}_{s,\Omega_2}
+h_2\norm{\frac{\beta}{\beta_2}f_2-f}_{0,\Omega_2} \right)
\endaligned
\]
with $1<r<3/2$ and $3/2<s\le2$, defined in Remark~\ref{rmk:regularity},
and for $1/2<t<1$.
\label{prop:error_estimate}
\end{proposition}

\begin{proof}
Thanks to the continuous elker and inf-sup conditions and their discrete
version in Propositions~\ref{prop:discreteelker} and~\ref{prop:discreteinfsup}, the theory of the saddle point problem yields the usual quasi-optimal error estimate (see~\cite[Th.~5.2.2]{boffi2013mixed}), that in this case, reads as follows
\[ \aligned
\norm{u-u_h}_{V}&+\norm{u_2-u_{2h}}_{V_2}+\norm{\lambda-\lambda_h}_{\Lambda}\\
&\le C\left(\inf_{v\in V_h}\norm{u-v}_V+\inf_{v_2\in V_{2h}}\norm{u_2-v_2}_{V_2}
+\inf_{\mu\in\Lambda_h}\norm{\lambda-\mu}_\Lambda\right).
\endaligned
\]
As a consequence of the regularity of $u$ and $u_2$ reported in
Remark~\ref{rmk:regularity},  standard arguments on the approximation error for
the spaces $V_h$ and $V_{2h}$ imply that
\[
\aligned
&\inf_{v\in V_h}\norm{u-v}_V\le Ch^{r-1}|u|_{r,\Omega}\\
&\inf_{v_2\in V_{2h}}\norm{u_2-v_2}_{V_2}\le Ch_2^{s-1}|u_2|_{s,\Omega_2}.
\endaligned
\]
The estimate of the best approximation of $\lambda$ requires a more careful
analysis. A similar consideration has been performed
in~\cite[Prop.~6]{auricchio2015fictitious} in a different setting. Actually,
taking into account~\eqref{eq:dual_lambda},
the Lagrange multiplier can be represented as the sum of two pieces
$\lambda=\lambda_1+\lambda_2$ with
\begin{equation}
\label{eq:lambda12}
\aligned
&\langle\lambda_1,v_2\rangle=\int_{\Omega_2}((\beta/\beta_2)f_2-f)v_2\,dx\\
&\langle\lambda_2,v_2\rangle=
\int_\Gamma(\beta_2-\beta)\nabla u_2\cdot\bfn_2 v_2\,d\gamma.
\endaligned
\end{equation}
We treat separately the two pieces and look for two elements in $\Lambda_h$
which approximate $\lambda_1$ and $\lambda_2$.

Since $f\in L^2(\Omega)$, $f_2\in L^2(\Omega_2)$ and since we assumed the
bound~\eqref{eq:coeff}, the first equality
in~\eqref{eq:lambda12} implies that $\lambda_1\in L^2(\Omega_2)$ with 
$\norm{\lambda_1}_{0,\Omega_2}\le\norm{(\beta/\beta_2)f_2-f}_{0,\Omega_2}$.
Let $P_0:L^2(\Omega_2)\to\Lambda_h$ be the $L^2$-projection onto $\Lambda_h$,
we set $\lambda_{1h}=P_0\lambda_1\in\Lambda_h$ with
\[
(\lambda_1-\lambda_{1h},\mu_h)=0\quad\forall \mu_h\in\Lambda_h.
\]
Then we observe that
\[
\norm{\lambda_1-\lambda_{1h}}_{\Lambda}
=\sup_{v_2\in V_2}
\frac{\langle\lambda_1-\lambda_{1h},v_2\rangle}{\norm{v_2}_{V_2}}
=\sup_{v_2\in V_2}\frac{(\lambda_1-\lambda_{1h},v_2)}{\norm{v_2}_{V_2}},
\]
but
\[
\aligned
(\lambda_1-\lambda_{1h},v_2)&=(\lambda_1-\lambda_{1h},v_2-P_0 v_2)
=(\lambda_1,v_2-P_0 v_2)\\
&\le\norm{\lambda_1}_{0,\Omega_2}\norm{v_2-P_0 v_2}_{0,\Omega_2}
\le C\norm{(\beta/\beta_2)f_2-f}_{0,\Omega_2}h_2\norm{v_2}_{V_2}.
\endaligned
\]
Therefore we end up with
\begin{equation}
\label{eq:estlambda1}
\norm{\lambda_1-\lambda_{1h}}_{\Lambda}\le 
C h_2 \norm{(\beta/\beta_2)f_2-f}_{0,\Omega_2}.
\end{equation}
Let us now construct an approximation of $\lambda_2$ and bound the approximation
error. Remark~\ref{rmk:regularity} states that $u_2\in H^s(\Omega_2)$ for
$3/2<s\le2$. Therefore the trace of the normal derivative of $u_2$ belongs
to $H^{s-3/2}(\Gamma)$ and we can infer from the second equality
in~\eqref{eq:lambda12} that $\lambda_2\in H^{-t}(\Omega_2)$ with $1/2<t<1$,
namely by trace inequality we have
\[
\aligned
\norm{\lambda_2}_{H^{-t}(\Omega_2)}&=\sup_{v_2\in H^t(\Omega_2)}
\frac{\langle\lambda_2,v_2\rangle}{\norm{v_2}_{H^t(\Omega_2)}}
=\sup_{v_2\in H^t(\Omega_2)}
\frac{\int_\Gamma (\beta_2-\beta)\nabla u_2\cdot\textbf{n}_2 v_2\,d\gamma}
{\norm{v_2}_{H^t(\Omega_2)}}\\
&\le \abs{\beta_2-\beta} \sup_{v_2\in H^t(\Omega_2)}\frac
{\norm{\nabla u_2\cdot\textbf{n}_2}_{H^{s-3/2}(\Gamma)}
\norm{v_2}_{H^{3/2-s}(\Gamma)}}{\norm{v_2}_{H^t(\Omega_2)}}\\
&\le C\sup_{v_2\in H^t(\Omega_2)}\frac
{\norm{u_2}_{H^{s}(\Omega_2)}
\norm{v_2}_{H^{t-1/2}(\Gamma)}}
{\norm{v_2}_{H^t(\Omega_2)}}\\
&\le C\sup_{v_2\in H^t(\Omega_2)}\frac
{\norm{u_2}_{H^{s}(\Omega_2)}
\norm{v_2}_{H^{t}(\Omega_2)}}{\norm{v_2}_{H^t(\Omega_2)}}
=C\norm{u_2}_{H^{s}(\Omega_2)}.
\endaligned
\]
Let $\lambda_{2h}\in\Lambda_h$ be such that
\[
\langle\lambda_{2h},v_{2h}\rangle
=\int_\Gamma (\beta_2-\beta)\nabla u_2\cdot\textbf{n}_2 v_{2h}\,d\gamma
\quad\forall v_{2h}\in V_{2h}.
\]
Then, it is easily seen that $\langle\lambda_2-\lambda_{2h},v_{2h}\rangle=0$
for all $v_{2h}\in V_{2h}$ so that we have for $v_{2h}\in V_{2h}$
\begin{equation}
\label{eq:firststepl2}
\langle\lambda_2-\lambda_{2h},v_{2}\rangle
=\langle\lambda_2-\lambda_{2h},v_{2}-v_{2h}\rangle
=\langle\lambda_2,v_{2}-v_{2h}\rangle-(\lambda_{2h},v_{2}-v_{2h}).
\end{equation}

We now choose $v_{2h}$ so that $v_{2h}|_\Gamma$ is the $L^2(\Gamma)$
projection of $v_2|_\Gamma$ and at the interior nodes of $\Omega_2$ coincides
with the Cl\'ement interpolant of $v_2$. In particular, we have the estimates
\[
\aligned
&\|v_2-v_{2h}\|_{0,\Gamma}\le Ch_2^{1/2}\|v_2\|_{H^{1/2}(\Gamma)}\\
&\|v_2-v_{2h}\|_{0,\Omega_2}\le Ch_2\|v_2\|_{1,\Omega_2}.
\endaligned
\]
Hence, using~\eqref{eq:lambda12}, we get

\[
\aligned
\langle\lambda_2,v_{2}-v_{2h}\rangle&=
\int_\Gamma (\beta_2-\beta)\nabla u_2\cdot\textbf{n}_2(v_{2}-v_{2h})\,d\gamma\\
&\le \abs{\beta_2-\beta} \norm{\nabla u_2\cdot\textbf{n}_2}_{H^{s-3/2}(\Gamma)}
\norm{v_2-v_{2h}}_{H^{3/2-s}(\Gamma)}.
\endaligned
\]
The last norm can be estimated as follows
\[
\aligned
\norm{v_2-v_{2h}}_{H^{3/2-s}(\Gamma)}&=\sup_{w\in H^{s-3/2}(\Gamma)}
\frac{(v_2-v_{2h},w)_{\Gamma}}{\|w\|_{H^{s-3/2}(\Gamma)}}\\
&=\sup_{w\in H^{s-3/2}(\Gamma)}
\frac{(v_2-v_{2h},w-w^I)_{\Gamma}}{\|w\|_{H^{s-3/2}(\Gamma)}}\\
&\le Ch_2^{s-3/2}\|v_2-v_{2h}\|_{0,\Gamma}\\
&\le Ch_2^{s-1}\|v_2\|_{H^{1/2}(\Gamma)}
\le Ch_2^{s-1}\|v_2\|_{1,\Omega_2},
\endaligned
\]
where $w^I$ is the Cl\'ement interpolant of $w$ satisfying
\[
\norm{w-w^I}_{0,\Gamma}\le Ch_2^{s-3/2}\|w\|_{H^{s-3/2}(\Gamma)}.
\]

It remains to bound the second term in~\eqref{eq:firststepl2}: we have
\[
(\lambda_{2h},v_{2}-v_{2h})\le
\norm{\lambda_{2h}}_{0,\Omega_2}\norm{v_{2}-v_{2h}}_{0,\Omega_2}
\le C h_2 \norm{\lambda_{2h}}_{0,\Omega_2}\norm{v_2}_{1,\Omega_2}.
\]
%
%
In order to conclude the proof, we now estimate
$\norm{\lambda_{2h}}_{0,\Omega_2}$. Recalling the Fortin operator $\Pi_h$
satisfying~\eqref{eq:Fortin}, we have
\[
\aligned
\norm{\lambda_{2h}}_{0,\Omega_2}&=
\sup_{v_2\in V_2}\frac{\langle\lambda_{2h},v_2\rangle}
{\norm{v_2}_{0,\Omega_2}}
=\sup_{v_2\in V_2}\frac{\langle\lambda_{2h},\Pi_h v_2\rangle}
{\norm{v_2}_{0,\Omega_2}}
=\sup_{v_2\in V_2}\frac{\langle\lambda_2,\Pi_h v_2\rangle}
{\norm{v_2}_{0,\Omega_2}}\\
&\le \sup_{v_2\in V_2}\frac{\norm{\lambda_2}_{H^{-t}(\Omega_2)}
\norm{\Pi_h v_2}_{H^t(\Omega_2)}}{\norm{v_2}_{0,\Omega_2}}\\
&
\le Ch_2^{-t}\norm{u_2}_{H^s(\Omega_2)}\sup_{v_2\in V_2}
\frac{\norm{\Pi_h v_2}_{0,\Omega_2}}{\norm{v_2}_{0,\Omega_2}}.
\endaligned
\]
The last supremum can be bounded if we provide a uniform $L^2(\Omega_2)$ estimate of the Fortin operator.
Looking at the construction of the Fortin operator as
$\Pi_h v_2=\Pi_1 v_2+\Pi_2(v_2-\Pi_1v_2)$, and inspecting the proof of its stability, it turns out that this can be done if we replace the Cl\'ement interpolation $\Pi_1$ with an interpolation that is bounded in $L^2(\Omega_2)$. This can be done by adopting the quasi-interpolation operator introduced in~\cite{Ern}, denoted again $\Pi_1$ for simplicity, and which has been proved to have the required stability
\[
\|\Pi_1v_2\|_{0,\Omega_2}\le C\|v_2\|_{0,\Omega_2}.
\]

We finally get
\[
(\lambda_{2h},v_{2}-v_{2h})\le
C h_2\norm{\lambda_{2h}}_{0,\Omega_2}\norm{v_2}_{1,\Omega_2}\le
C h_2^{1-t}\|u_2\|_{H^s(\Omega_2)}\|v_2\|_{H^1(\Omega_2)}.
\]
Putting together all the pieces we see that the best approximations of $u$,
$u_2$, and $\lambda$ converge with order $h^{r-1}$, $h_2^{s-1}$, and
$h_2^{1-t}$, respectively. In the case when $h_2<1$, since
$s-1>1-t$, we get that $h_2^{s-1}$ is dominated by $h_2^{1-t}$. Therefore, the
final rate of convergence is given by the maximum between $h^{r-1}$ and
$h_2^{1-t}$.
\end{proof}

\begin{remark}
This paper only considers quadrilateral and hexahedra mesh. However, all results and proofs might be easily extended to triangular and tetrahedra  mesh with the corresponding element choices $P_1-(P_1 +\overline{B})-P_0$ and $P_2-(P_2+\overline{B})-P_0$ where $\overline{B}$ is a proper bubble function associated with an internal node. 
\end{remark}

\section{Numerical results}
\label{se:Numerical_results}
In this section, we delve into a comprehensive numerical study of the elliptic interface problem in the FDDLM formulation. Our primary objective is to investigate the accuracy, stability, and convergence rate of the two proposed elements. 

We consider a two-dimensional setting, and we fix the right-hand sides of the problem as $f=1$ and $f_2=1$.

Our simulations involve various cases, including different ratios between the mesh sizes of the two domains, as well as different combinations of coefficients between the two domains, including cases where the sufficient condition of the elker, $\beta_2-\beta>0$, is violated. In particular, we investigate the three cases: \\
Case 1: $\beta=1, \beta_2 =10$.\\
Case 2: $\beta=1, \beta_2 =10000$. \\
Case 3: $\beta=10, \beta_2 =1$.

To further assess the performance of the proposed elements, we consider four examples involving various immersed shapes. Specifically, we consider a square, an L-shape, a circle, and a flower, each equipped with suitable Dirichlet boundary conditions. These examples were chosen to demonstrate the robustness of our proposed method in dealing with different geometries. In each example, we choose the domains $\Omega$ and $\Omega_2$ as follows:\\
\textbf{Example 1:} We consider a rectangular mesh with a square immersed shape mesh, where $\Omega=[0,6]^2$ and $\Omega_2=[e, 1+\pi]^2$. \\
\textbf{Example 2:} We investigate a rectangular mesh with an immersed L-shape mesh, where $\Omega=[0,6]^2$ and $\Omega_2=[1,3]^2\setminus[2,3]^2$. \\
\textbf{Example 3:} We study a rectangular mesh with a circular immersed shape mesh, where $\Omega=[-1.4,1.4]^2$ and $\Omega_2 =B((0,0),1)$.\\
\textbf{Example 4:} We examine a rectangular mesh with an immersed flower shape mesh, where $\Omega=[-2,3]^2$ and $\Omega_2$ is a flower with radius $1+0.1\cos(5\theta)$ and center (0,0).

For the first, second, and fourth examples, we obtain reference solutions by solving Laplace problems with jumping coefficients on a very fine mesh (up to 16 million degrees of freedom) using a high-order mixed finite element method. In the case of the circle, we have an exact solution in the different three cases given as follows:
\begin{align*}
\text{Case 1:}
\begin{cases}
U_1 = \dfrac{4-x^2-y^2}{4}& \text{in } \Omega1\\
U_2 = \dfrac{31-x^2-y^2}{40}& \text{in } \Omega2,\\
\end{cases}
\end{align*}
\begin{align*}
\text{Case 2:}
\begin{cases}
U_1 = \dfrac{4-x^2-y^2}{4}& \text{in } \Omega1\\
U_2 = \dfrac{30001-x^2-y^2}{40000}& \text{in } \Omega2,\\
\end{cases}
\end{align*}
\begin{align*}
\text{Case 3:}
\begin{cases}
U_1 = \dfrac{4-x^2-y^2}{40}& \text{in } \Omega1\\
U_2 = \dfrac{13-10x^2-10y^2}{40}& \text{in } \Omega2.\\
\end{cases}
\end{align*}

\begin{figure}
\begin{minipage}[c]{.45\linewidth}
\includegraphics[width=1.1 \linewidth  ] {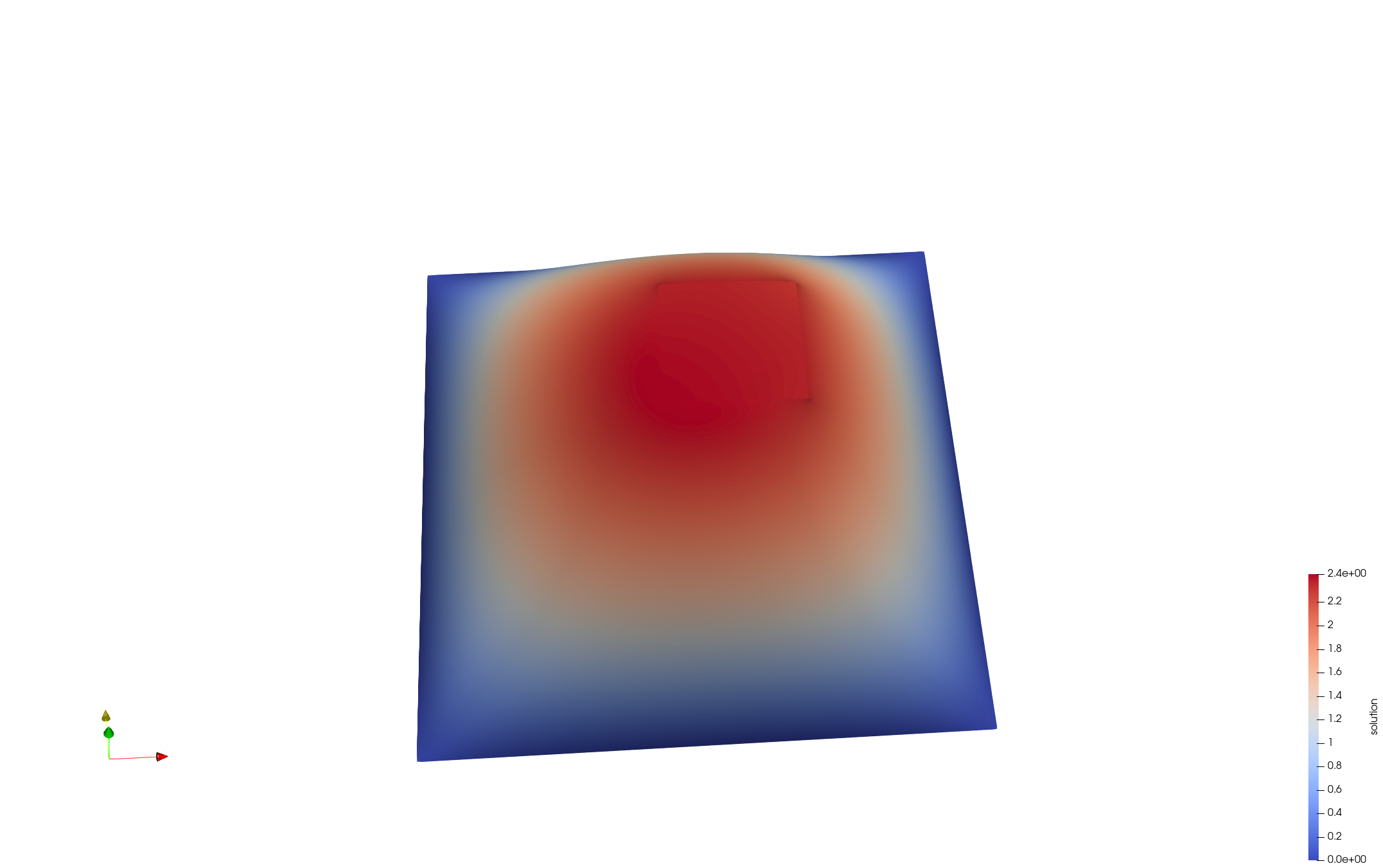}
\caption{ Exact solution of Example 1.}
\label{fig:exact_solution_square}
\end{minipage}
\begin{minipage}[c]{.45\linewidth}
\includegraphics[width=1.1 \linewidth]{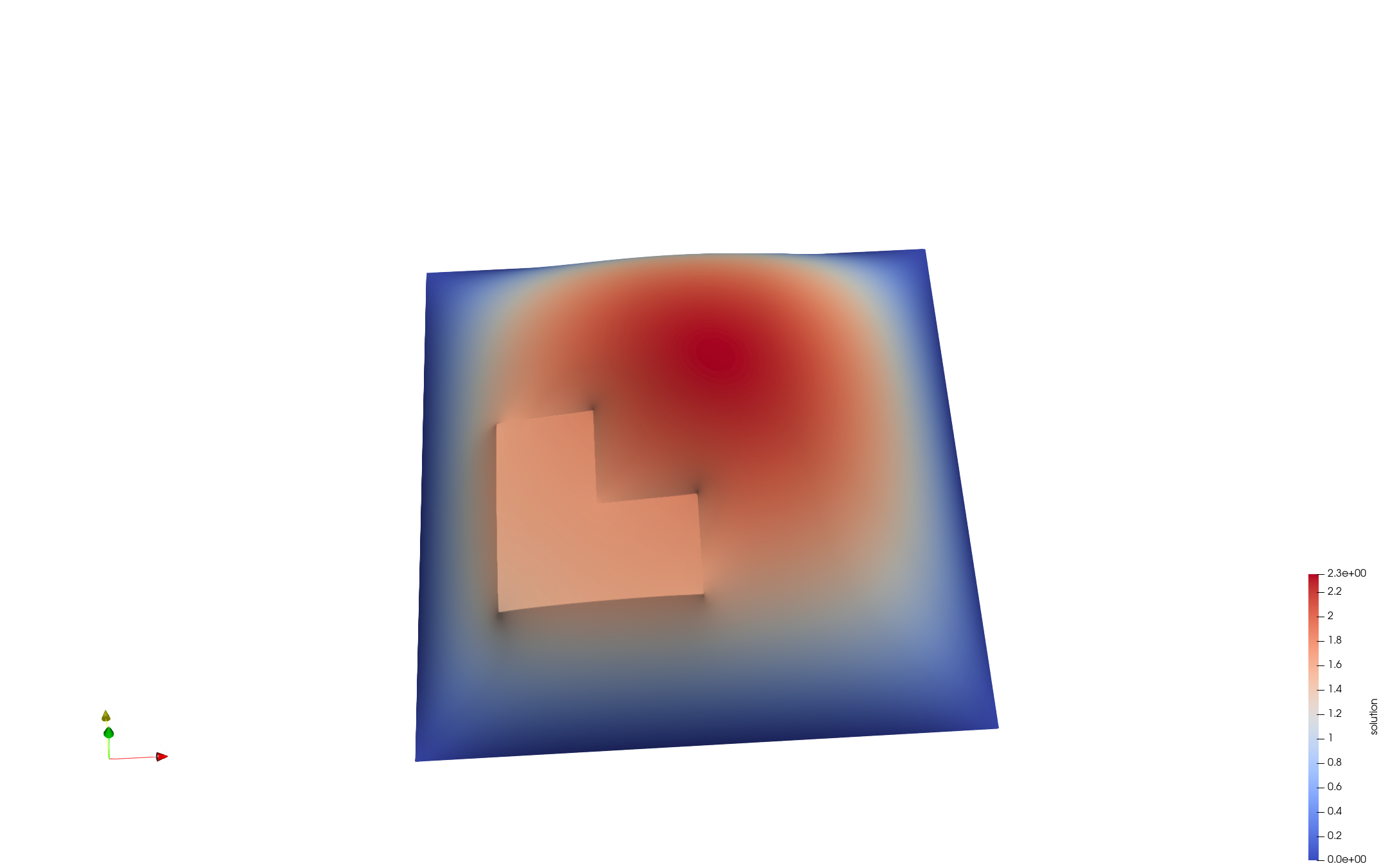}
\caption{ Exact solution of Example 2.}
\label{fig:exact_solution_Lshape}
\end{minipage}

\begin{minipage}[c]{.45\linewidth}
\includegraphics[width=1.1 \linewidth]{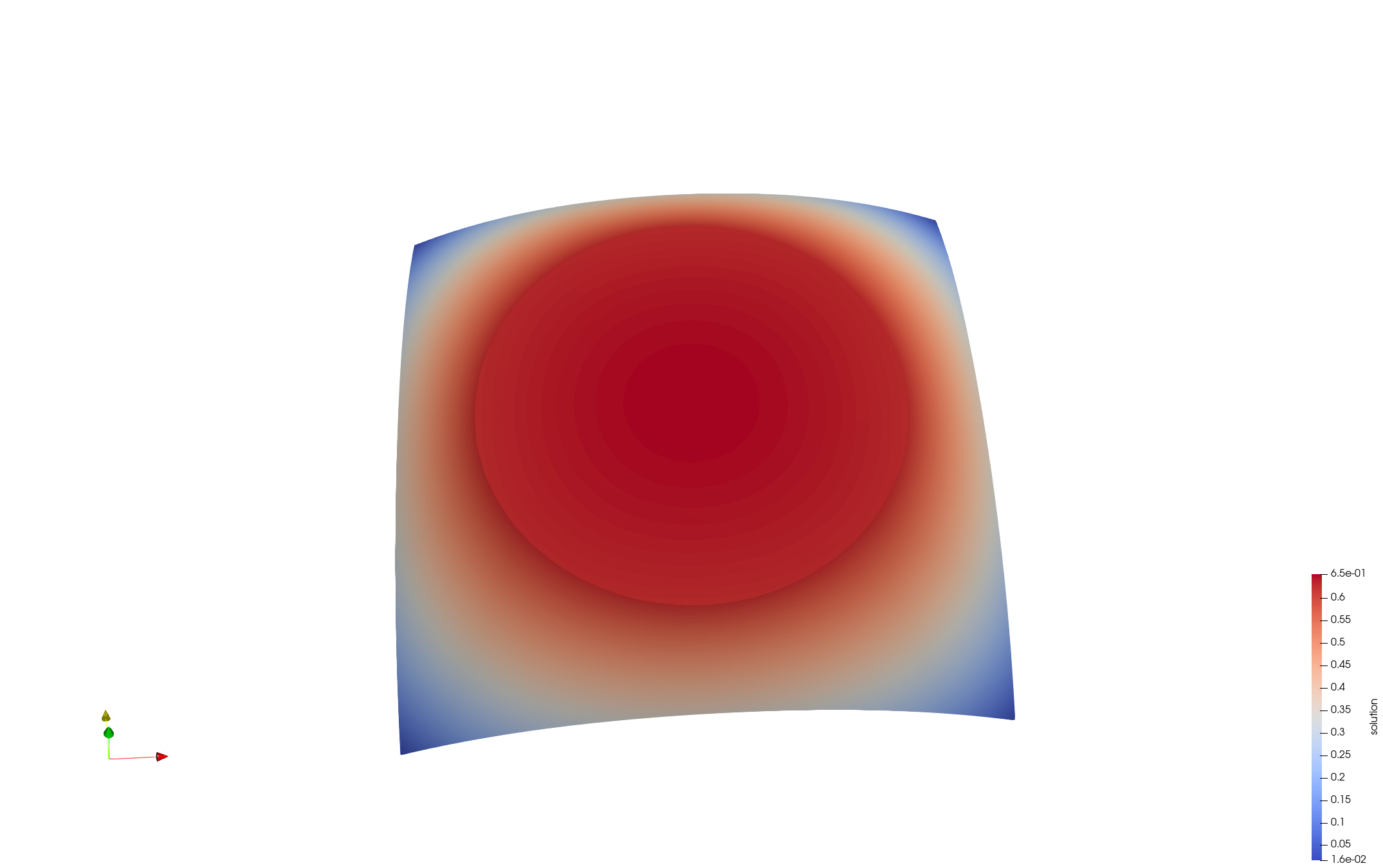}
\caption{ Exact solution of Example 3.}
\label{fig:exact_solution_circle}
\end{minipage}
\begin{minipage}[c]{.45\linewidth}
\includegraphics[width=1.1 \linewidth]{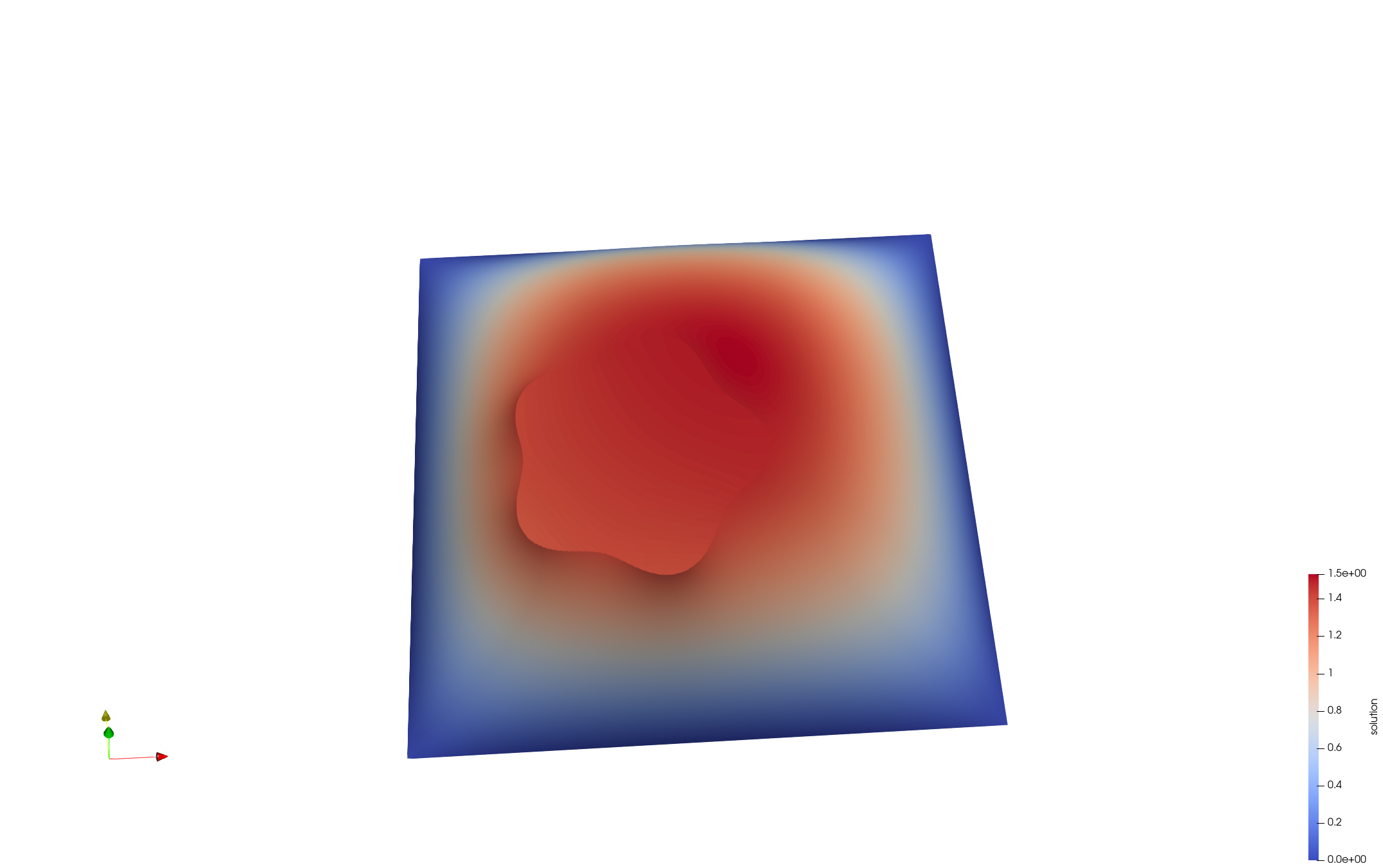}
\caption{ Exact solution of Example 4.}
\label{fig:exact_solution_flower}
\end{minipage}
\end{figure} 
We plot in Figures~\ref{fig:exact_solution_square},\ref{fig:exact_solution_Lshape},\ref{fig:exact_solution_circle},\ref{fig:exact_solution_flower} the exact$\setminus$reference solutions for each of the examples in the first case to provide a visual representation of the problems.

In any case, dealing with the coupling term arising from non-matching meshes between the two domains is one of the main challenges in interface problem simulations. To tackle this issue, we describe the numerical methods and techniques utilized in our study. Additionally, we analyze the need for the bubble function to ensure the stability of the proposed elements, and further investigate the inf-sup condition numerically. Lastly, we present a detailed convergence rate analysis and discuss the performance of the two proposed elements in all cases. We present some numerical tests that confirm the expected convergence rate of our schemes.

In the following, we will describe how we handle the numerical treatment of the coupling term $(\lambda_h, v_h|_{\Omega_2})_{\Omega_2}$ in Problem~\ref{pbm:2} for Example~2 (without loss of generality), where $\lambda_h \in \Lambda_h$ and $v_h \in V_h$. This term is represented as an $L^2$ scalar product in $\Omega_2$. Evaluating
this term requires the evaluation of the following integral:
\begin{equation}
\label{eq:integral}
\int_{\Omega_2} \phi_i \psi_j|_{\Omega_2}\, dx
=\sum_{K\in\mathcal{T}_2}\int_{K} \phi_i \psi_j|_{\Omega_2}\, dx
\end{equation}
where $\phi_i$ $i=1,\dots,\dim(\Lambda_h)$ and $\psi_j$ $j=1,\dots,\dim(V_h)$
are the basis functions that span $\Lambda_h$ and $V_h$, respectively.
Thanks to the choice of $\Lambda_h$, $\phi_i=1$ on the element
$K_i\in\mathcal{T}_2$ and vanishes elsewhere. Hence the integral
in~\eqref{eq:integral} reduces to the integral on the intersection between
$K_i$ with the support of $\psi_j$.

\begin{figure}[htp]
\centering
\begin{minipage}[c]{.22\linewidth}
\centering
	\includegraphics[width=0.9\textwidth]{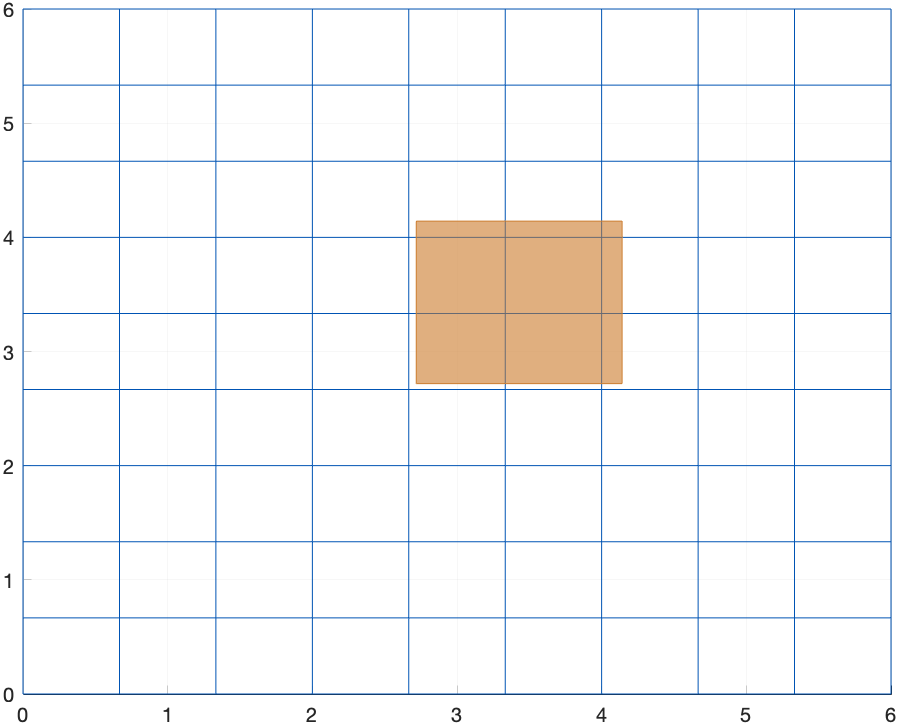} 
	\tiny{\text{A:  $\mathcal{T}_2$ and $\mathcal{T}$.}}
\end{minipage}
\begin{minipage}[c]{.02\linewidth}
\centering
\tiny $\rightarrow$
\end{minipage}
\begin{minipage}[c]{.22\linewidth}
\centering
	\includegraphics[width=0.9\textwidth]{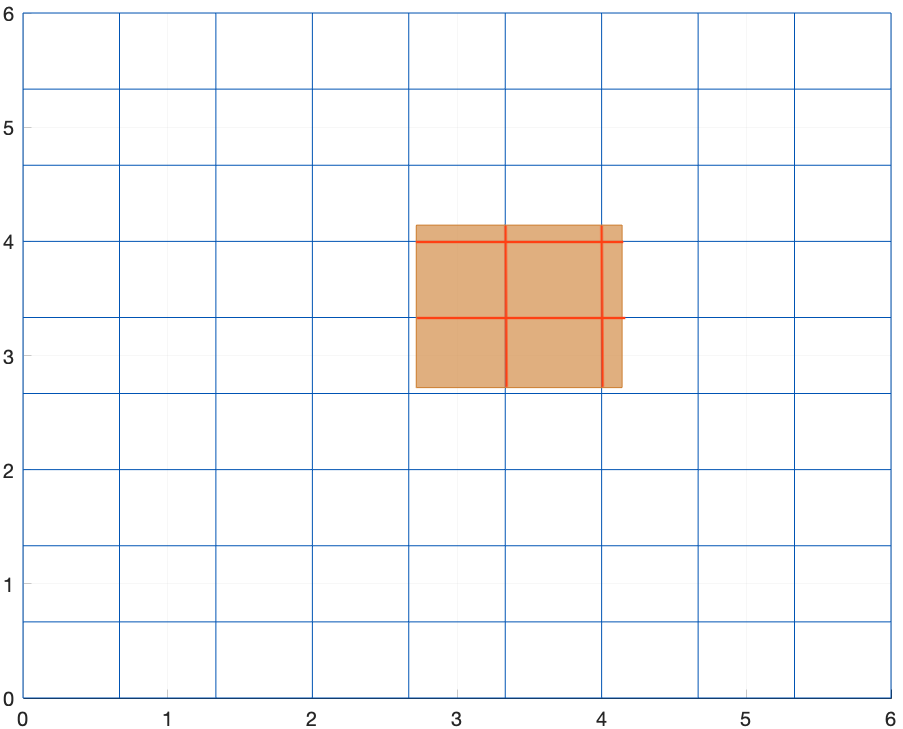} 
	\tiny{\text{B: $\overline{\mathcal{T}_2}$ is a  remesh of $\mathcal{T}_2$.}}
\end{minipage}
\begin{minipage}[c]{.02\linewidth}
\centering
\tiny $\rightarrow$
\end{minipage}
\begin{minipage}[c]{.22\linewidth}
\centering
	\includegraphics[width=0.9\textwidth]{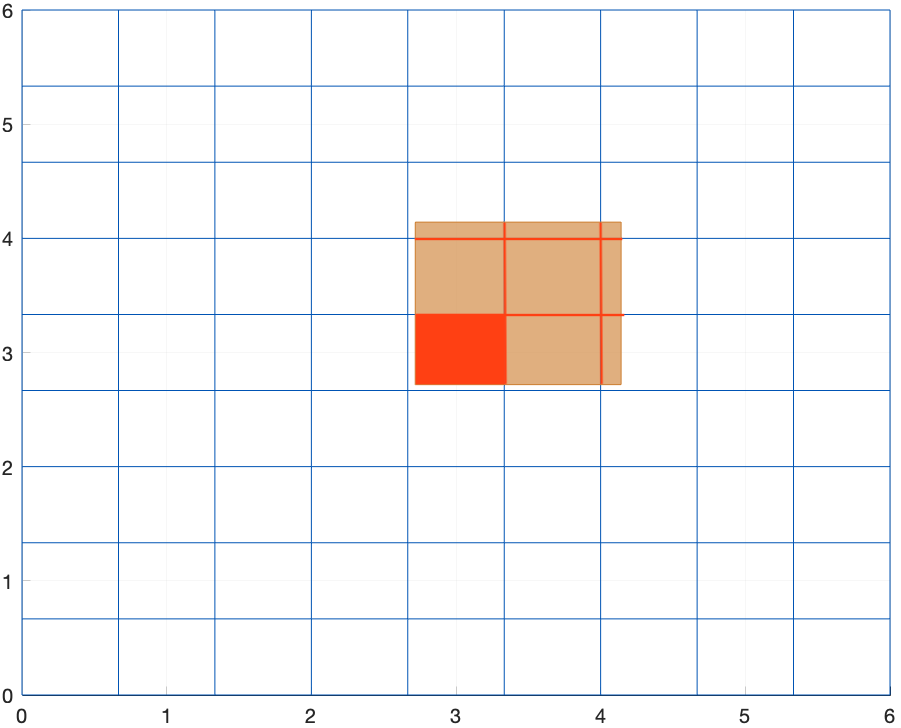} 
		\tiny{\text{C: Support of $\phi_i \psi_j$.}}
\end{minipage}
\begin{minipage}[c]{.02\linewidth}
\centering
\tiny $\rightarrow$
\end{minipage}
\begin{minipage}[c]{.22\linewidth}
\centering
	\includegraphics[width=0.9\textwidth]{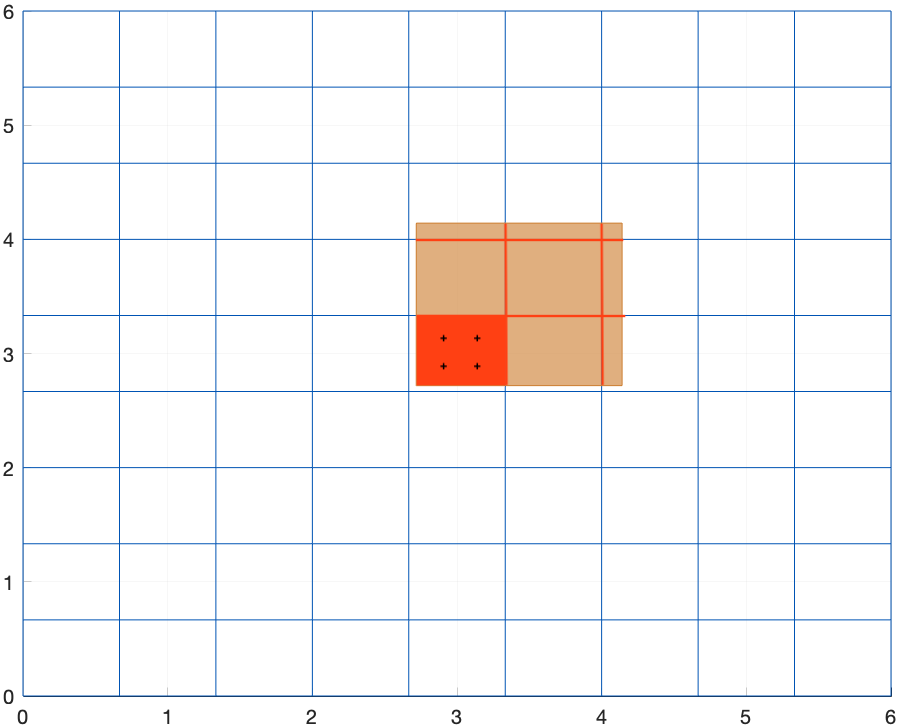} 
		\tiny{\text{D: Mapped quadrature points .}}
\end{minipage}
	\caption{ Intersection approach.}
	\label{fig:Intersection_approach}
\end{figure}

To simplify the idea, consider the 2D case and assume that the support of
$\phi_i$ is the element $K_i\in\mathcal{T}_2$ colored with brown in
Figure~\ref{fig:Intersection_approach}.A. We start by finding the geometric
intersections of elements in $\mathcal{T}$ with $K_i$. This leads us
to the introduction of a new mesh in $\Omega_2$, that we denote by 
$\mathcal{\overline{T}}_2$ with the property that 
each element $\overline{K}\in \mathcal{\overline{T}}_2$ is contained in 
one element in $\mathcal{T}$ as in Figure~\ref{fig:Intersection_approach}.B.
This new mesh is also a quadrilateral mesh but finer than $\mathcal{T}_2$.
Therefore we have 
\[
\int_{\Omega_2} \phi_i \psi_j|_{\Omega_2}\, dx=
\int_{K_i\cap\mathrm{supp}(\psi_j)} \psi_j\, dx
\]
where $K_i\cap\mathrm{supp}(\psi_j)$ is the rectangle colored in red in 
Figure~\ref{fig:Intersection_approach}.C.

This approach will give us an exact evaluation of the scalar product since it
takes into the account the exact region for which the two shape functions
interact~\cite{boffi2022interface}. The procedure for adapting this approach is
briefly summarized bellow:

\begin{itemize}
\item Choose the order of the quadrature rule depending on the degree of
$\psi_j$.
\item Find the corresponding quadrature points and weights in the reference
element $\widehat{K}$. Denote by $\widehat{q_k}$ the $k^{th}$ quadrature
point and by $\widehat{w_k}$ the associated quadrature weight.
\item Map the quadrature points and weights to $\overline{K}$ where
$\overline{K}=K_i\cap\mathrm{supp}(\psi_j)$. So, we have
$w_k=\abs{\overline{K}} \widehat{w_k}$, and
$\psi_j(q_k)=\widehat{\psi_j}(\widehat{q_k})$,
see, for example, Figure~\ref{fig:Intersection_approach}.D.
\item Evaluate the integral as follows:
\[
\int_{\Omega_2} \phi_i \psi_j
=\int_{\overline{K}} \psi_j 
=\abs{\overline{K}} \int_{\widehat{K}} \widehat{\psi_j}
=\abs{\overline{K}} \sum_{k} \widehat{w_k} \widehat{\psi_j}(\widehat{q_k})=
\sum_{k} w_k \psi_j(q_k).
\]
\end{itemize}

\subsection{Numerical discrete inf-sup test}
\label{sse:weakinfsup}
In this section, numerical tests on the discrete inf-sup bound of the proposed
elements is presented to validate the theoretical proofs given in
Section~\ref{se:A_Priori}.

For the same reason explained in Section~\ref{se:A_Priori}, it is enough to show that there exists a constant $\gamma_2 >0$ such that the following inf-sup bound is satisfied:
\begin{equation}
\sup_{v_{2h}\in V_{2h}} \dfrac{( \mu_h, v_{2h})_{\Omega_2} }{\norm{v_{2h}}_{V_2}} \geq \gamma_2 \norm{\mu_h}_{\Lambda}.
\end{equation}
In order to estimate numerically the constant $\gamma_2$, we use the following standard procedure.\\
Let $N_1$ and $N_2$ be the two matrices associated with the following norms:
\begin{align*}
\norm{v_2}_{V_2}^2     &=\abs{v_2}_{V_2}^2+\norm{v_2}_{0,\Omega_2}^2= v_2^T N_2 v_2 
							&& v_2 \in V_{2h}\\
\norm{\mu_h}_{\Lambda} ^2 &\approx h_2^2 \norm{\mu_h}_{0,\Omega_2} ^2 = \mu_h^T  \left(h_2^2 N_1\right)\mu_h
							&& \mu_h \in \Lambda_h.
\end{align*}
Arguing as in~\cite{malkus1981eigenproblems}, the eigenvalue equation associated with this inf-sup condition is
\begin{equation}\label{pb:eig_pbm}
\left(  C_2\right)  \left( h_2^2 N_1\right) ^{-1} \left( C_2^T\right)  v_2=\sigma  N_2 v_2
\end{equation}
where $C_2$ is the operator associated with the bilinear form $(\mu_h, v_{2h})_{\Omega_2}$ .
If the considered finite element satisfies the inf-sup condition, then, with increase refinements, the square root of the smallest eigenvalue $\sigma_1$ is bounded from below away from zero independently from the mesh sizes. This bound is the desired inf-sup  bound. 

\begin{figure}[H]
\centering
		\includegraphics[width=0.6 \linewidth]{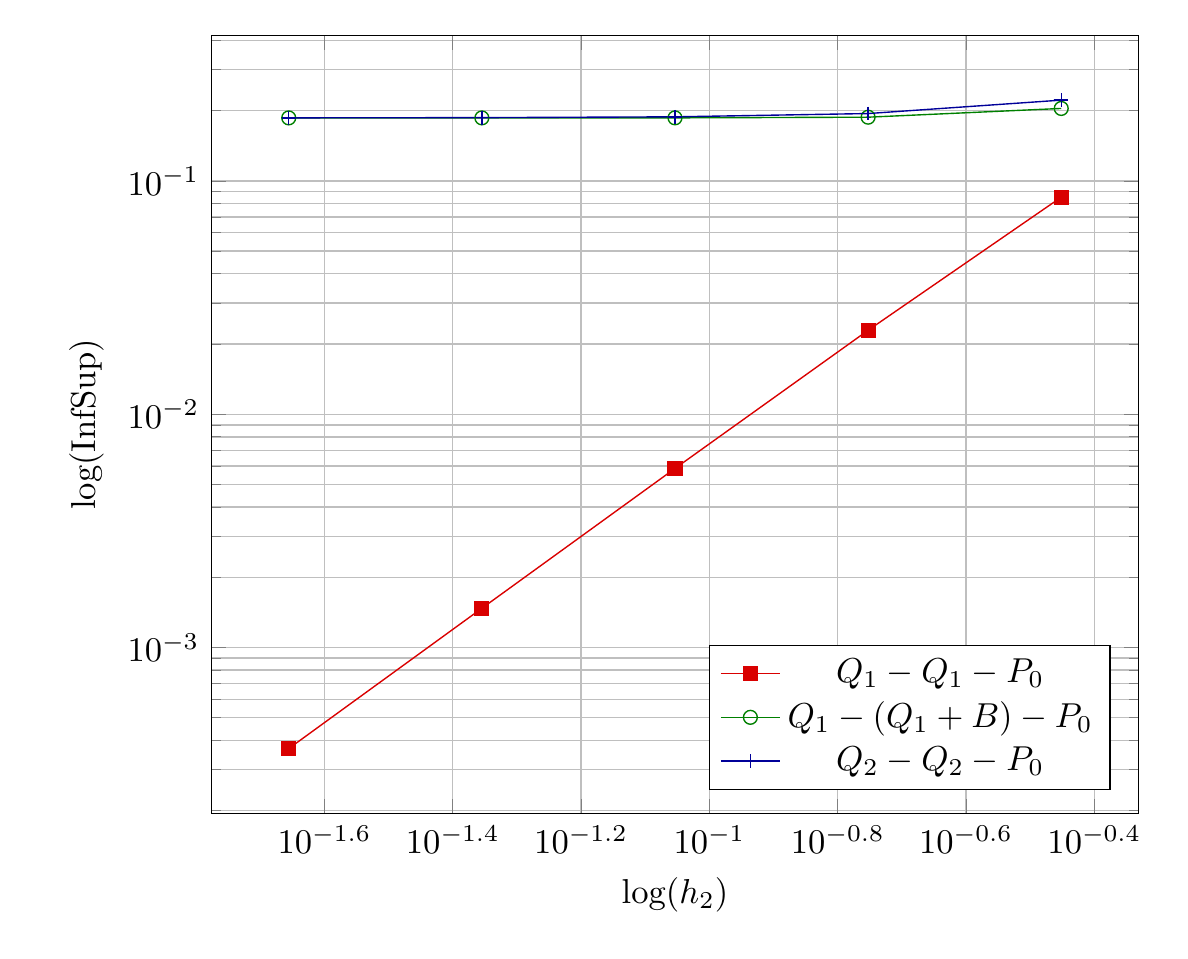} 
  \caption{Numerical discrete Inf-Sup Test.}
  \label{fig:infsup}
\end{figure}

The problem is solved in a sequence of five refinements. The results are
plotted using logarithmic scaling.
In Figure~\ref{fig:infsup} we report the results of our test for the two
elements presented in Section~\ref{se:A_Priori} together with the test for the
element $Q_1-Q_1-P_0$, where we didn't add the bubble to the space $V_{2h}$.

It is clear that Element~1 and Element~2 are stable as the mesh is refined,
that is the inf-sup constant doesn't degenerate, while the element without the
bubble is not, that is, the inf-sup constant tends to zero as $h$ goes to
zero.

\subsection{Rate of convergence}
\label{sse:rateofconv}
The rate of convergence is a critical aspect in numerical methods for solving partial differential equations as it measures how quickly the numerical solution approaches the exact solution with mesh refinement which depends on the smoothness of the solution and the efficiency of the method. In this section, we analyze the rate of convergence of our numerical methods for Examples 1, 2, 3, and 4, comparing the convergence rates with the theoretical rates predicted by the finite element method to assess the effectiveness of our proposed method in dealing with different immersed shapes.

To gain a better understanding of the performance of our numerical method, we conduct extensive testing for each example, examining three different choices of the ratio between mesh sizes. Furthermore, we explore the performance of our method by allowing the diffusion coefficients to vary considering Cases 1, 2, and 3 introduced at the beginning of this section. We analyze the behavior of the solution using the two proposed elements, namely $(Q_1,Q_1+B,P_0)$ and $(Q_2,Q_2,P_0)$. In addition, we investigate the numerical behavior of the $(Q_1,Q_1,P_0)$ element, which we previously demonstrated that it fails to pass the numerical inf-sup test.

\begin{figure}[H]
\centering
\begin{minipage}[c]{0.3\linewidth}
		\includegraphics[width=1\linewidth]{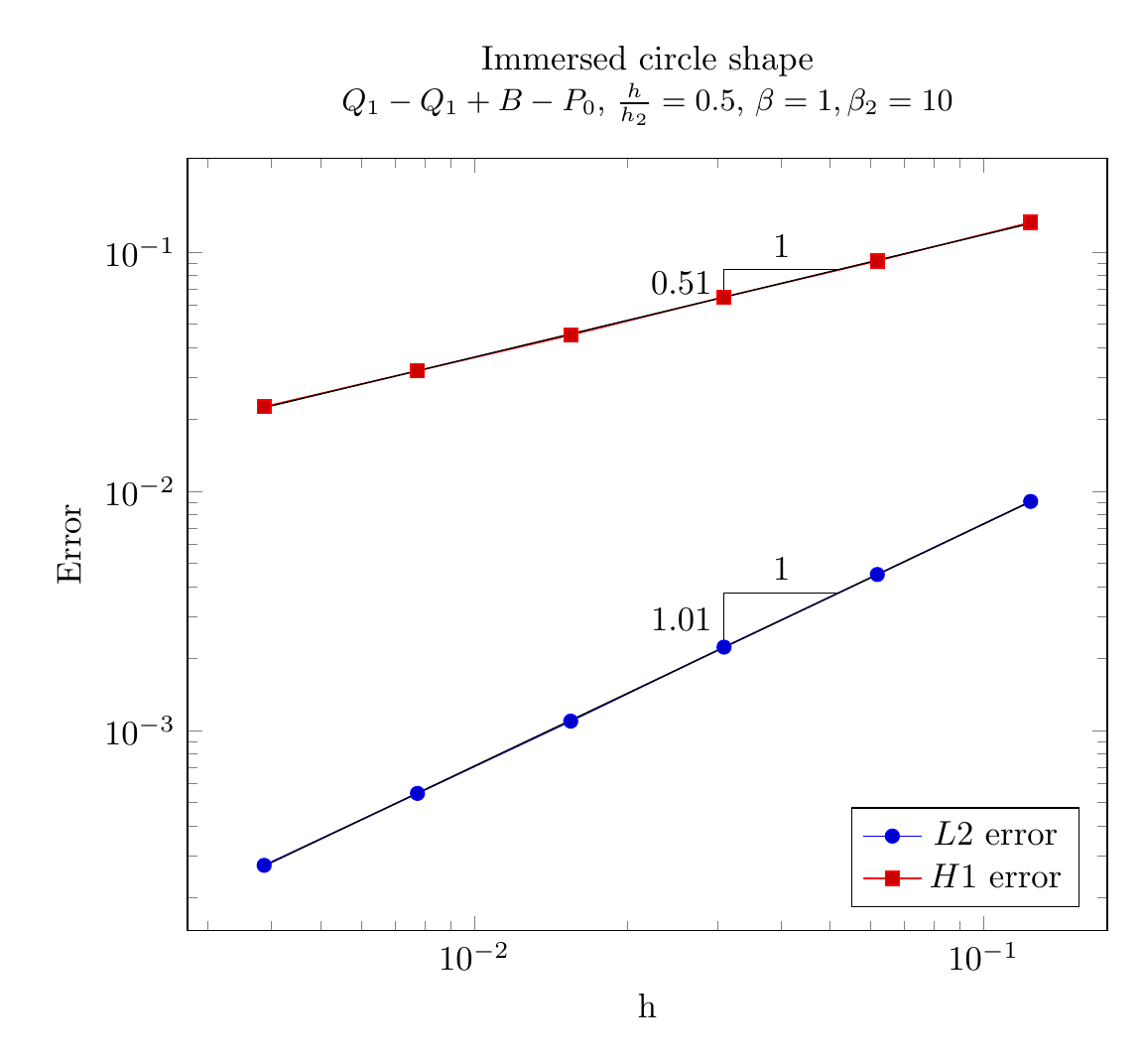} 
\end{minipage}
\begin{minipage}[c]{.3\linewidth}
		\includegraphics[width=1\linewidth]{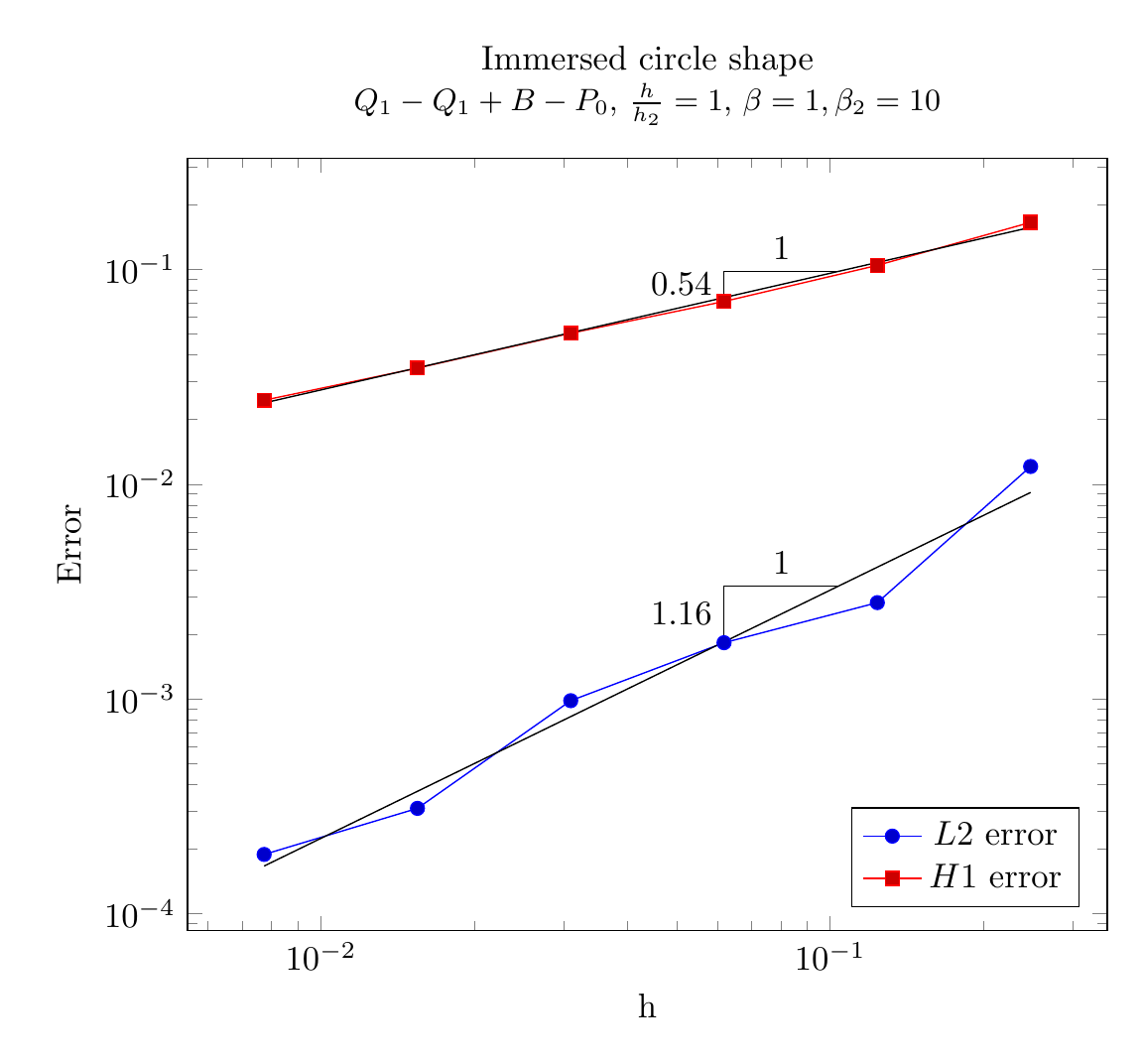} 
\end{minipage}
\begin{minipage}[c]{0.3\linewidth}
		\includegraphics[width=1\linewidth]{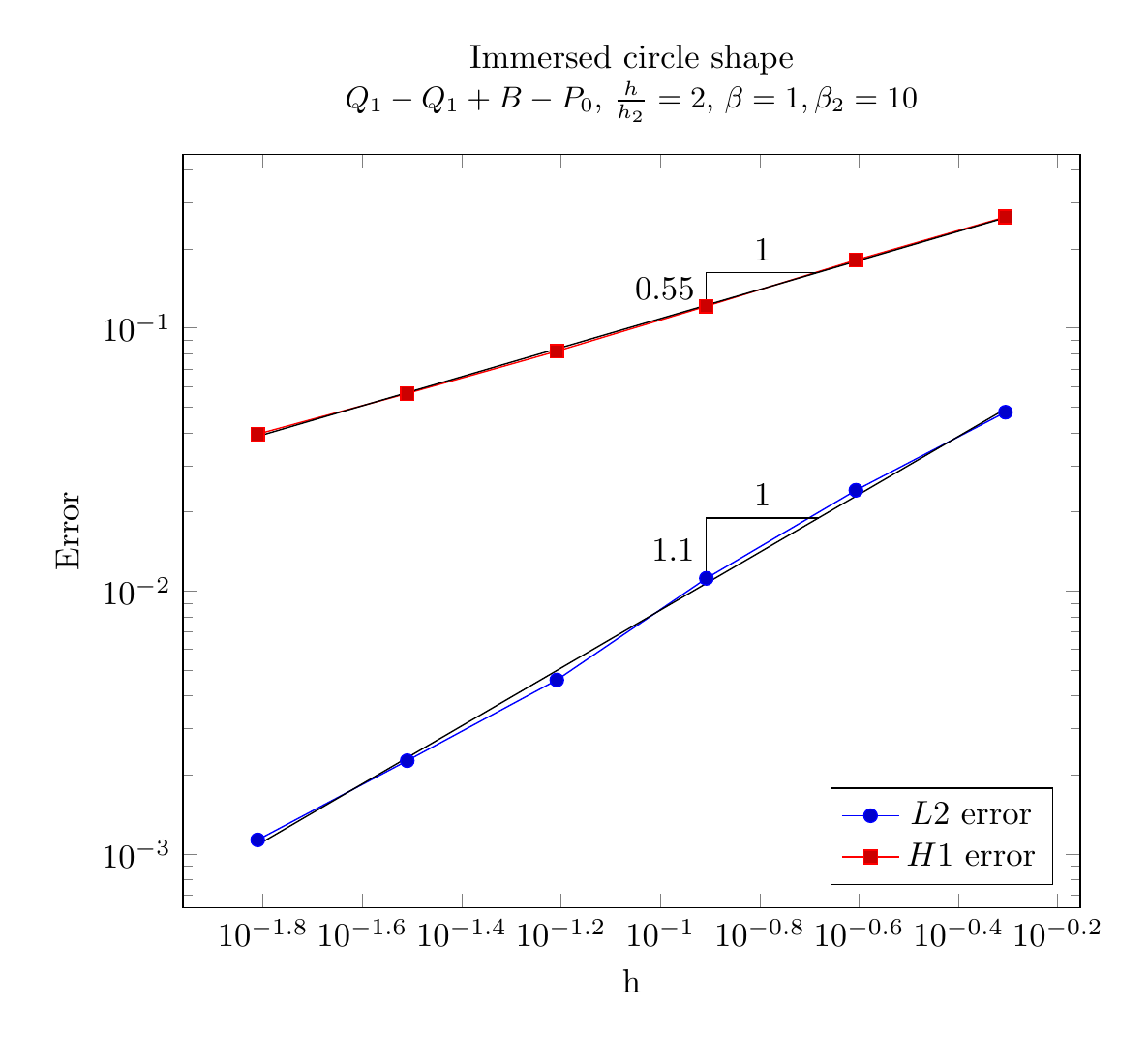} 
\end{minipage}
\captionof{figure}{Convergence of $Q_1-(Q_1+B)-P_0$ FDDLM: $\beta=1~\beta_2=10$.}
  \label{fig:error1_1_0_110_c}
%
\begin{minipage}[c]{0.3\linewidth}
		\includegraphics[width=1\linewidth]{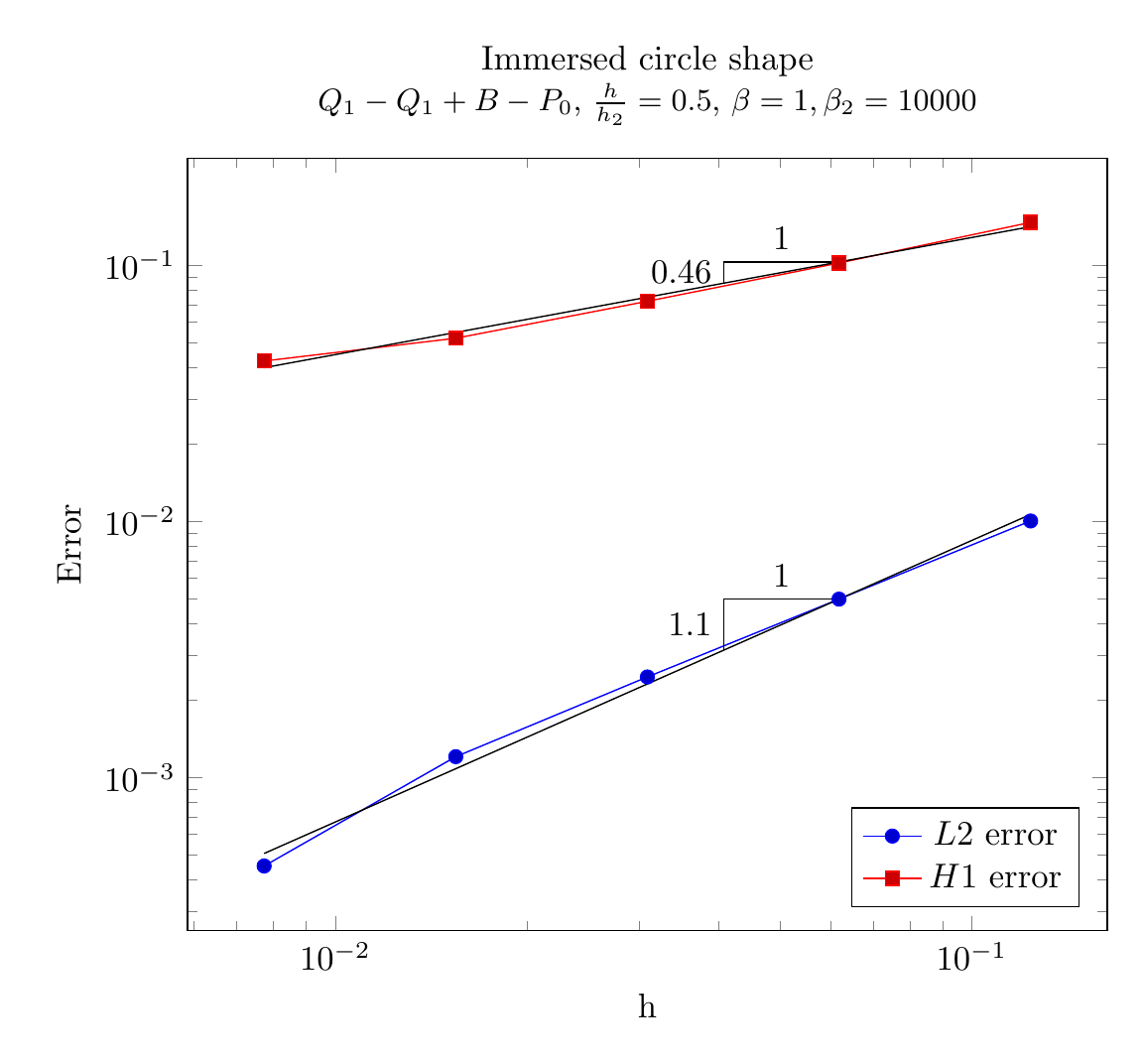} 
\end{minipage}
\begin{minipage}[c]{.3\linewidth}
		\includegraphics[width=1\linewidth]{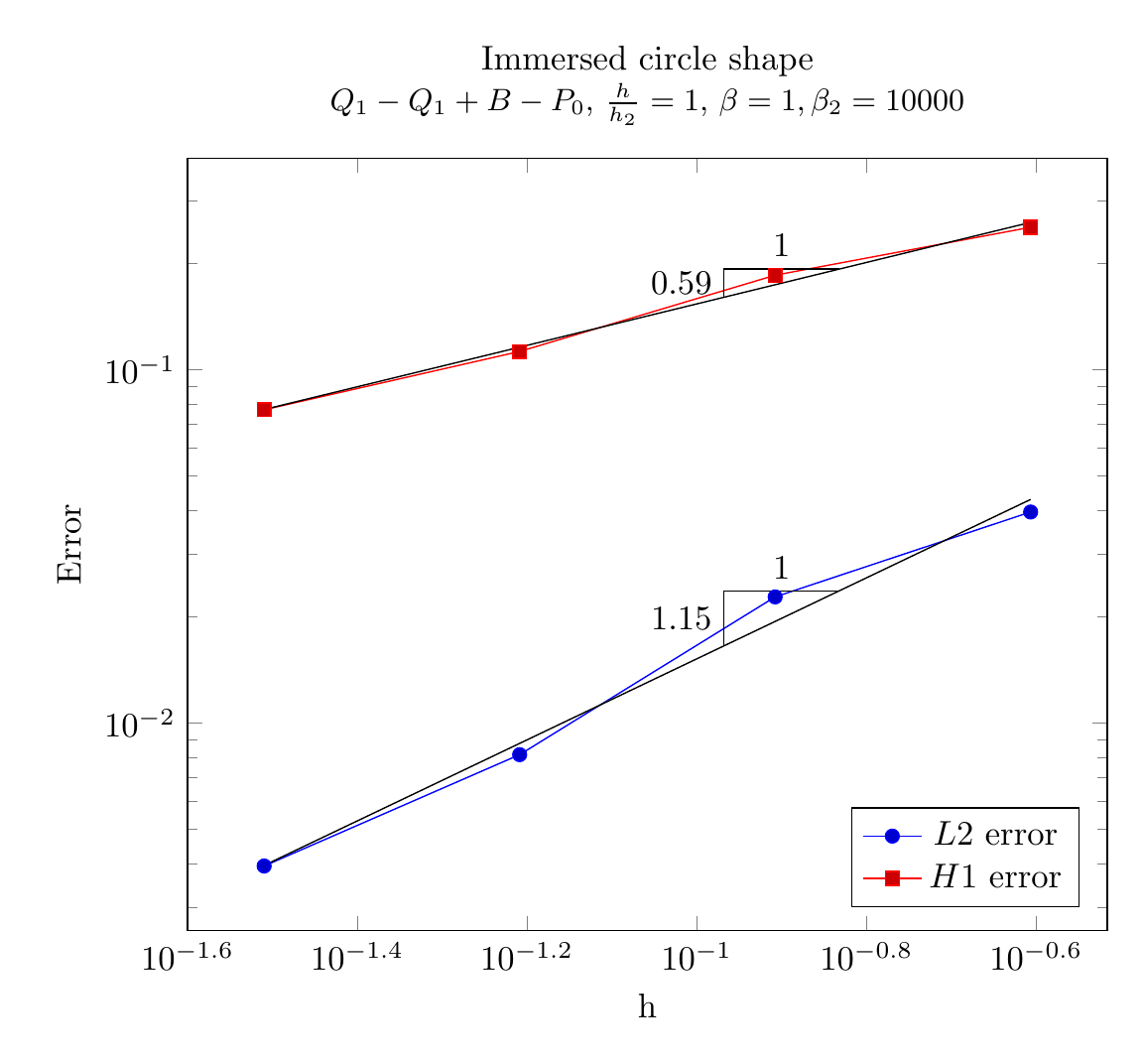} 
\end{minipage}
\begin{minipage}[c]{0.3\linewidth}
		\includegraphics[width=1\linewidth]{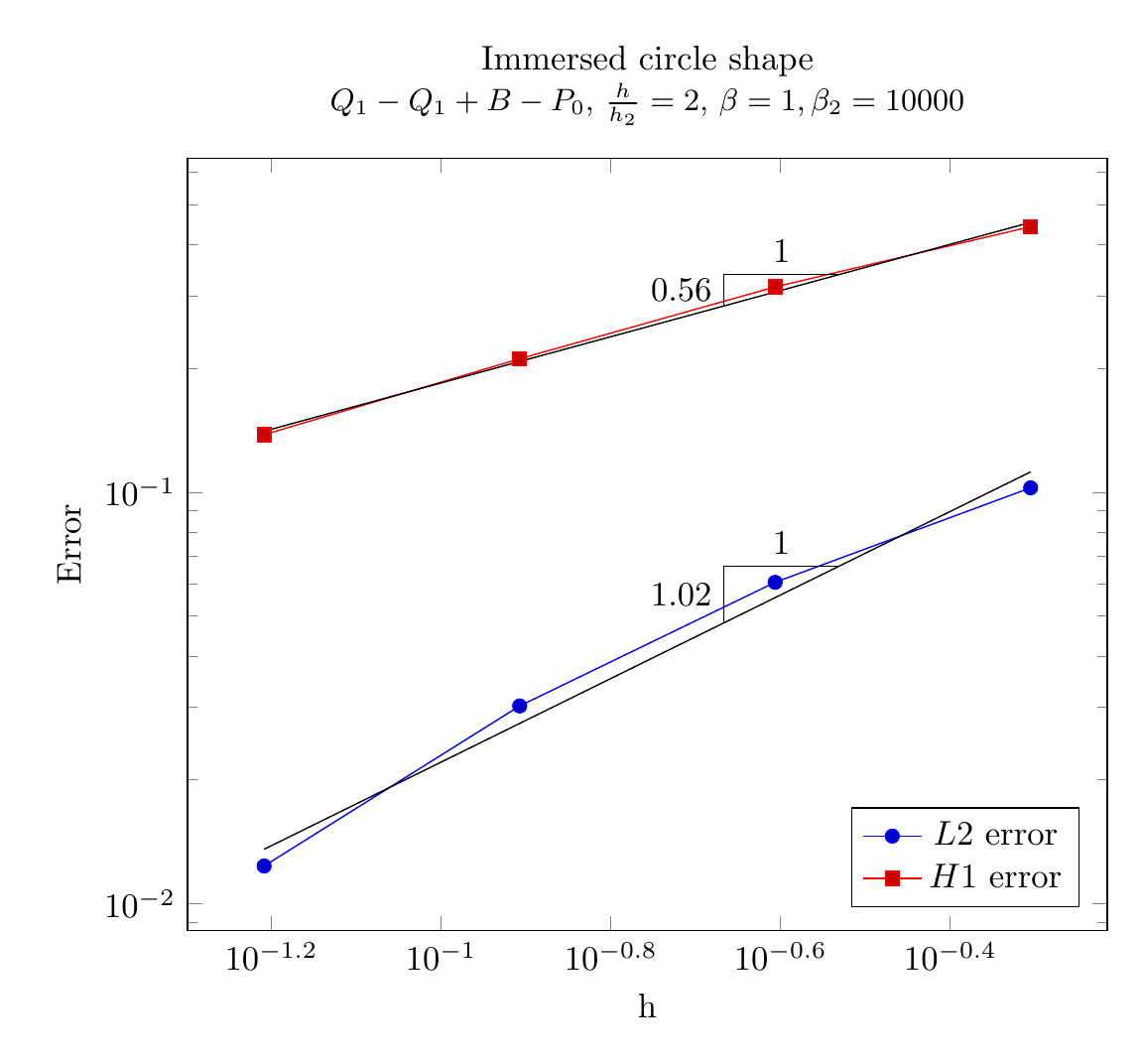} 
\end{minipage}
\captionof{figure}{Convergence of $Q_1-(Q_1+B)-P_0$ FDDLM: $\beta=1~\beta_2=10000$.}
  \label{fig:error1_1_0_110000_c}
\begin{minipage}[c]{0.3\linewidth}
		\includegraphics[width=1\linewidth]{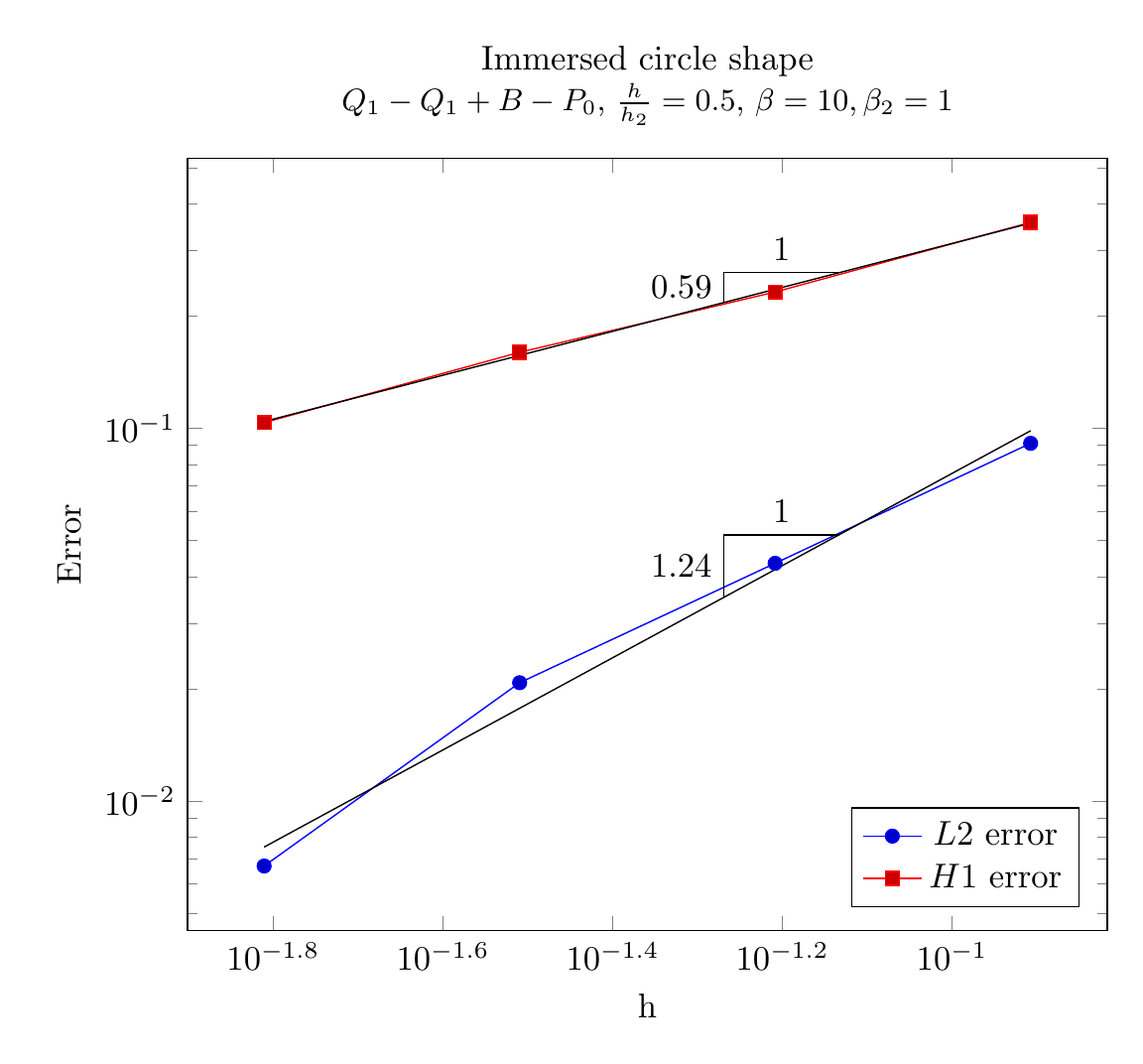} 
\end{minipage}
\begin{minipage}[c]{.3\linewidth}
		\includegraphics[width=1\linewidth]{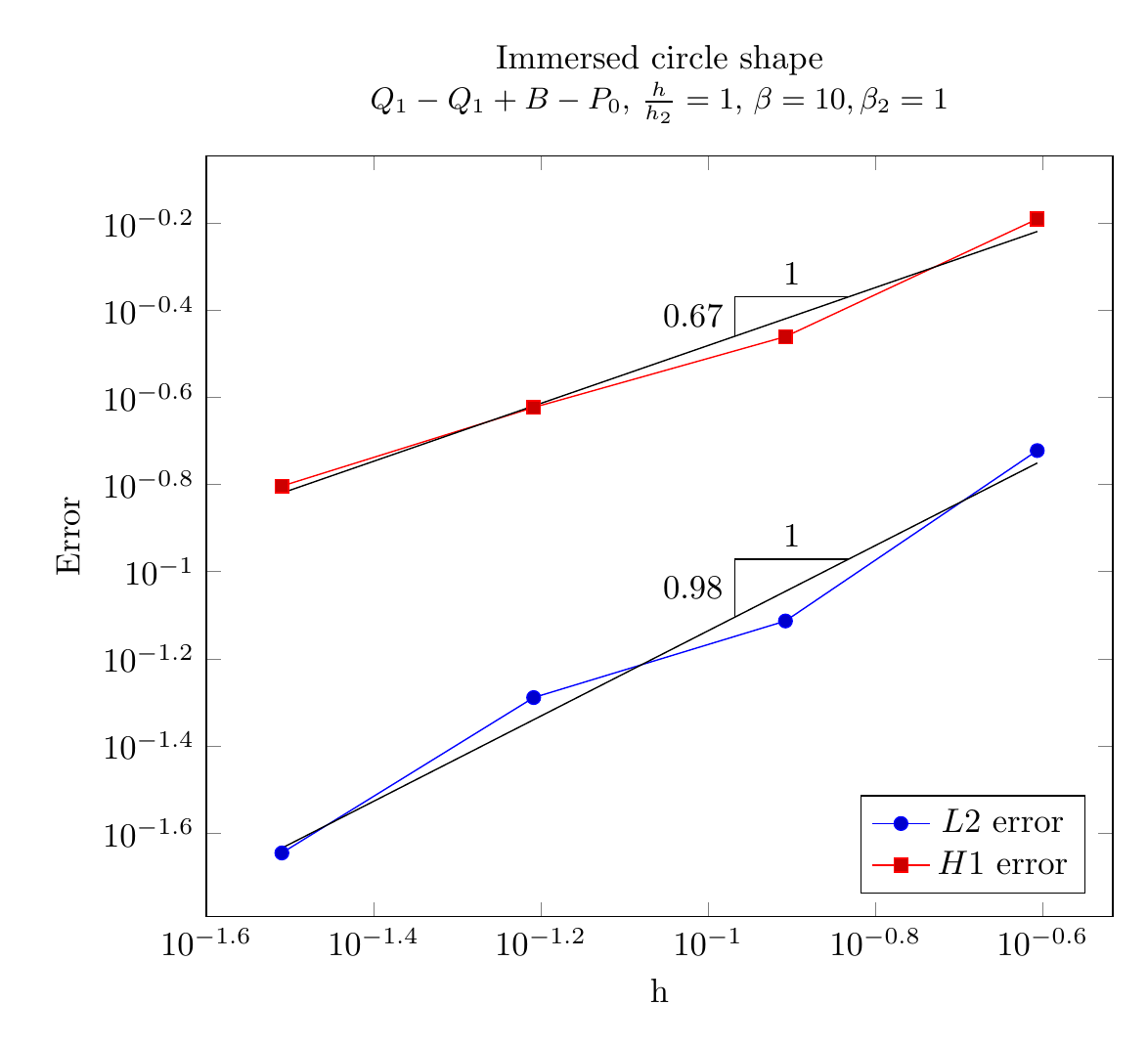} 
\end{minipage}
\begin{minipage}[c]{0.3\linewidth}
		\includegraphics[width=1\linewidth]{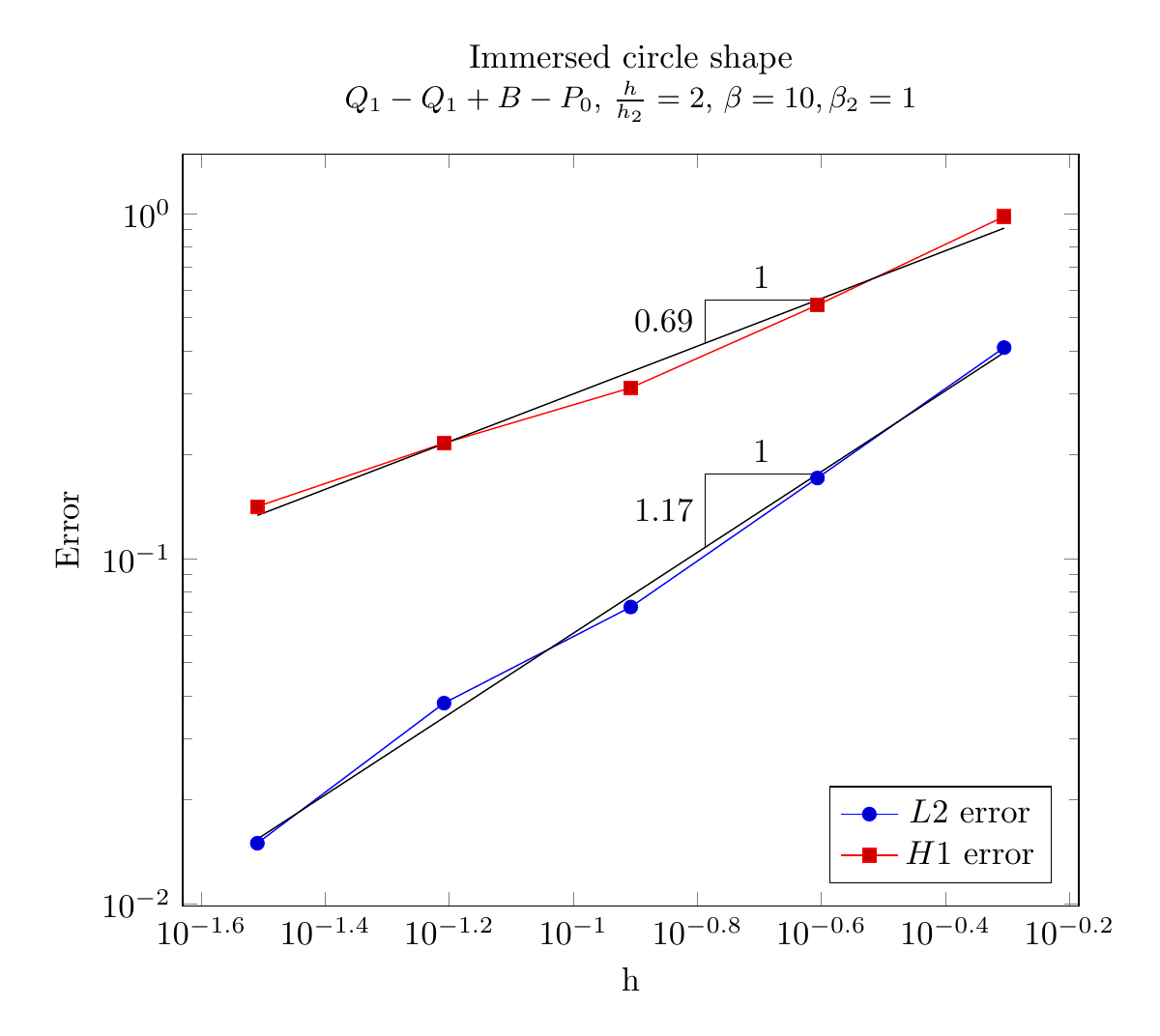} 
\end{minipage}
\captionof{figure}{Convergence of $Q_1-(Q_1+B)-P_0$ FDDLM: $\beta=10~\beta_2=1$.}
  \label{fig:error1_1_0_101_c}

\end{figure}

\begin{figure}[H]
\centering
\begin{minipage}[c]{0.3\linewidth}
		\includegraphics[width=1\linewidth]{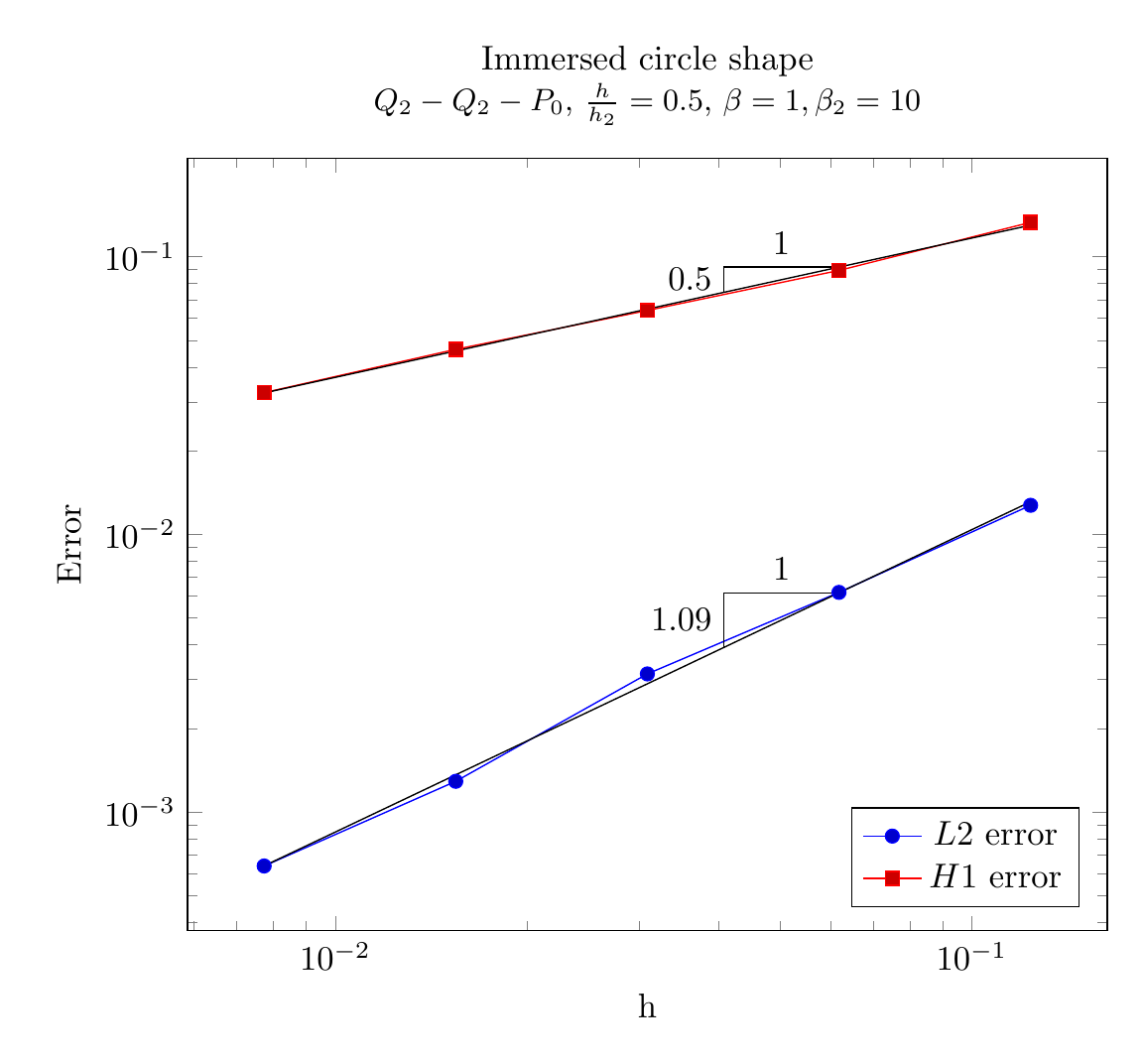} 
\end{minipage}
\begin{minipage}[c]{.3\linewidth}
		\includegraphics[width=1\linewidth]{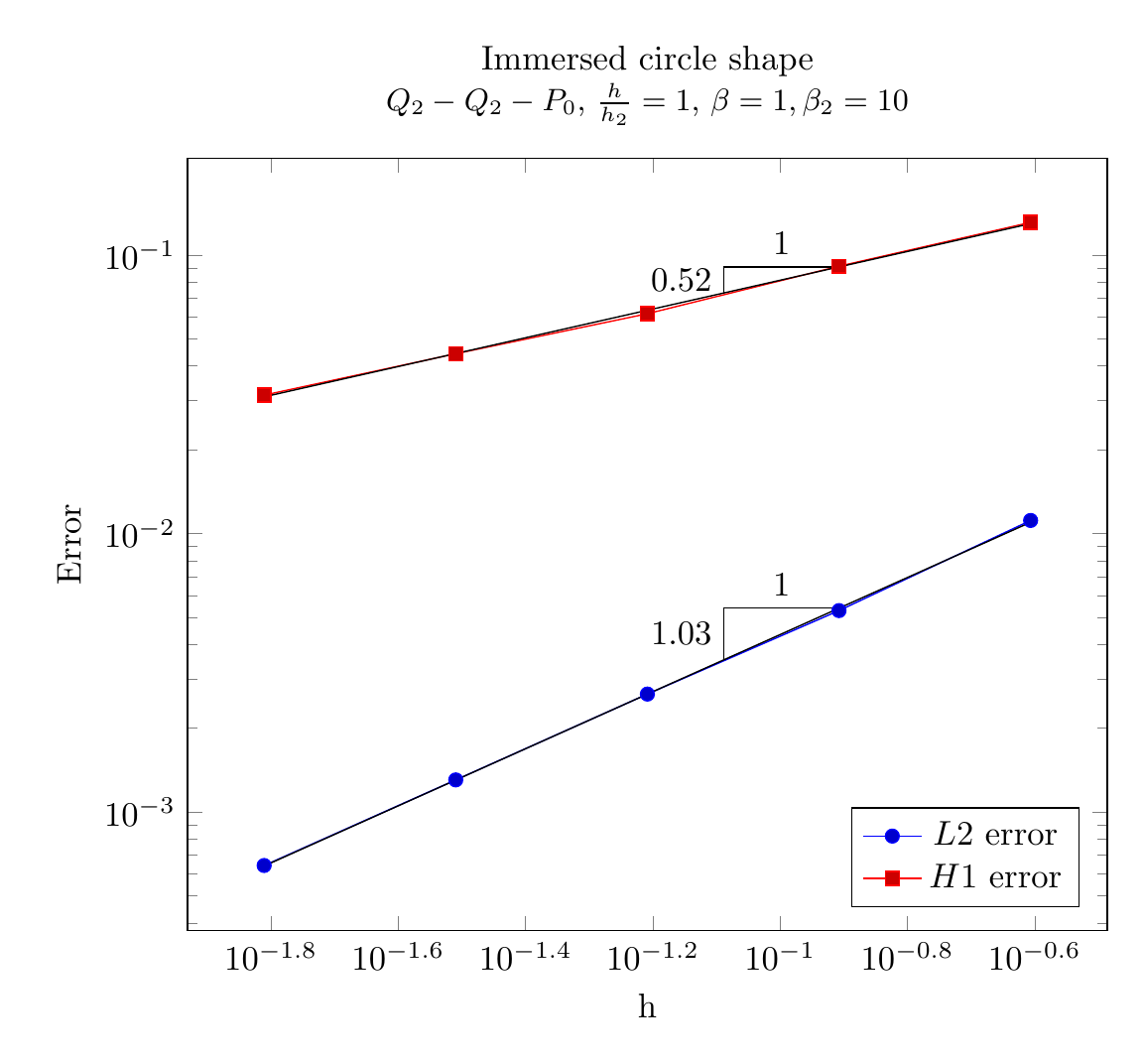} 
\end{minipage}
\begin{minipage}[c]{0.3\linewidth}
		\includegraphics[width=1\linewidth]{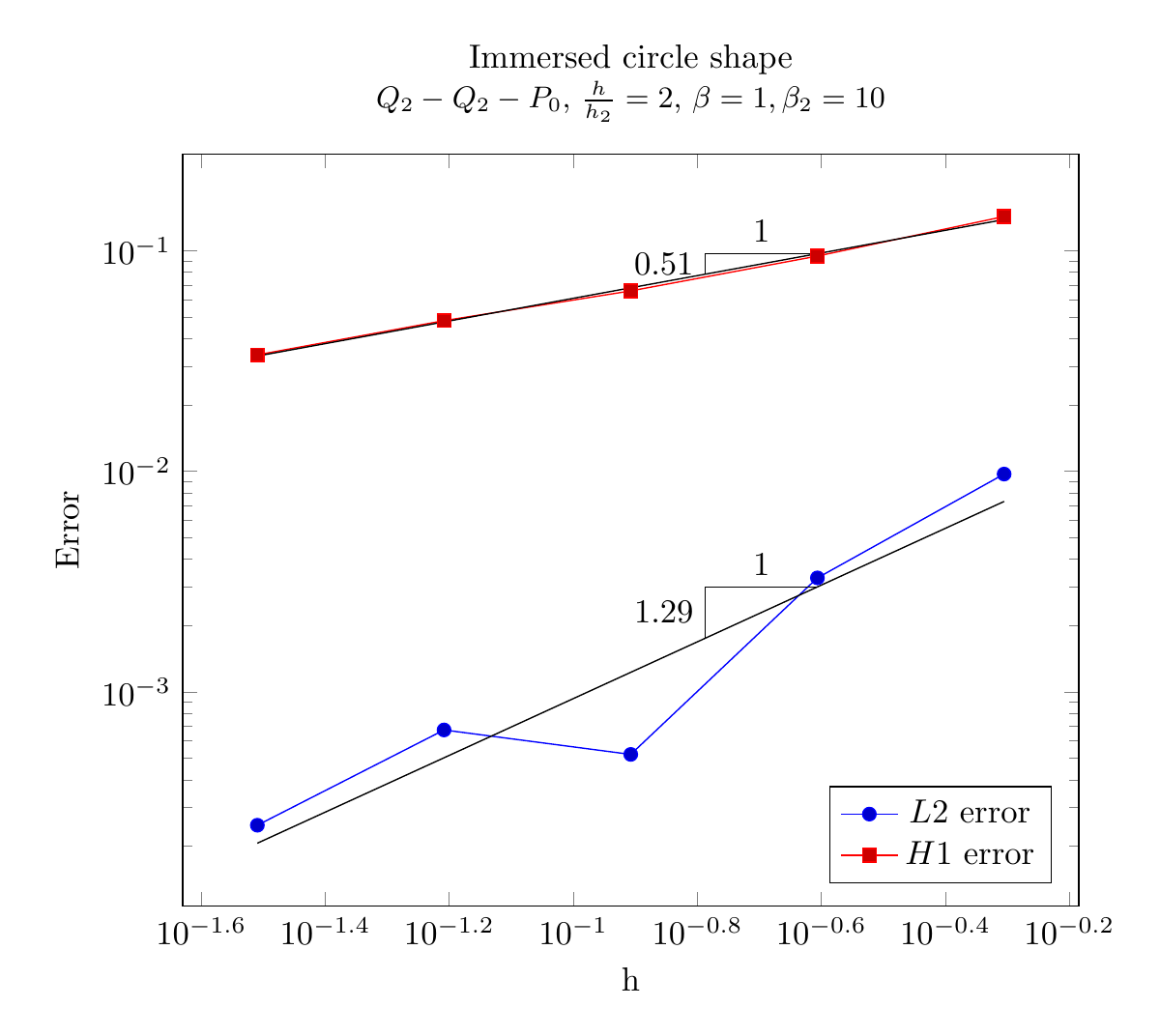} 
\end{minipage}
\captionof{figure}{Convergence of $Q_2-Q_2-P_0$ FDDLM: $\beta=1~\beta_2=10$.}
  \label{fig:error220_110_c}
%
\begin{minipage}[c]{0.3\linewidth}
		\includegraphics[width=1\linewidth]{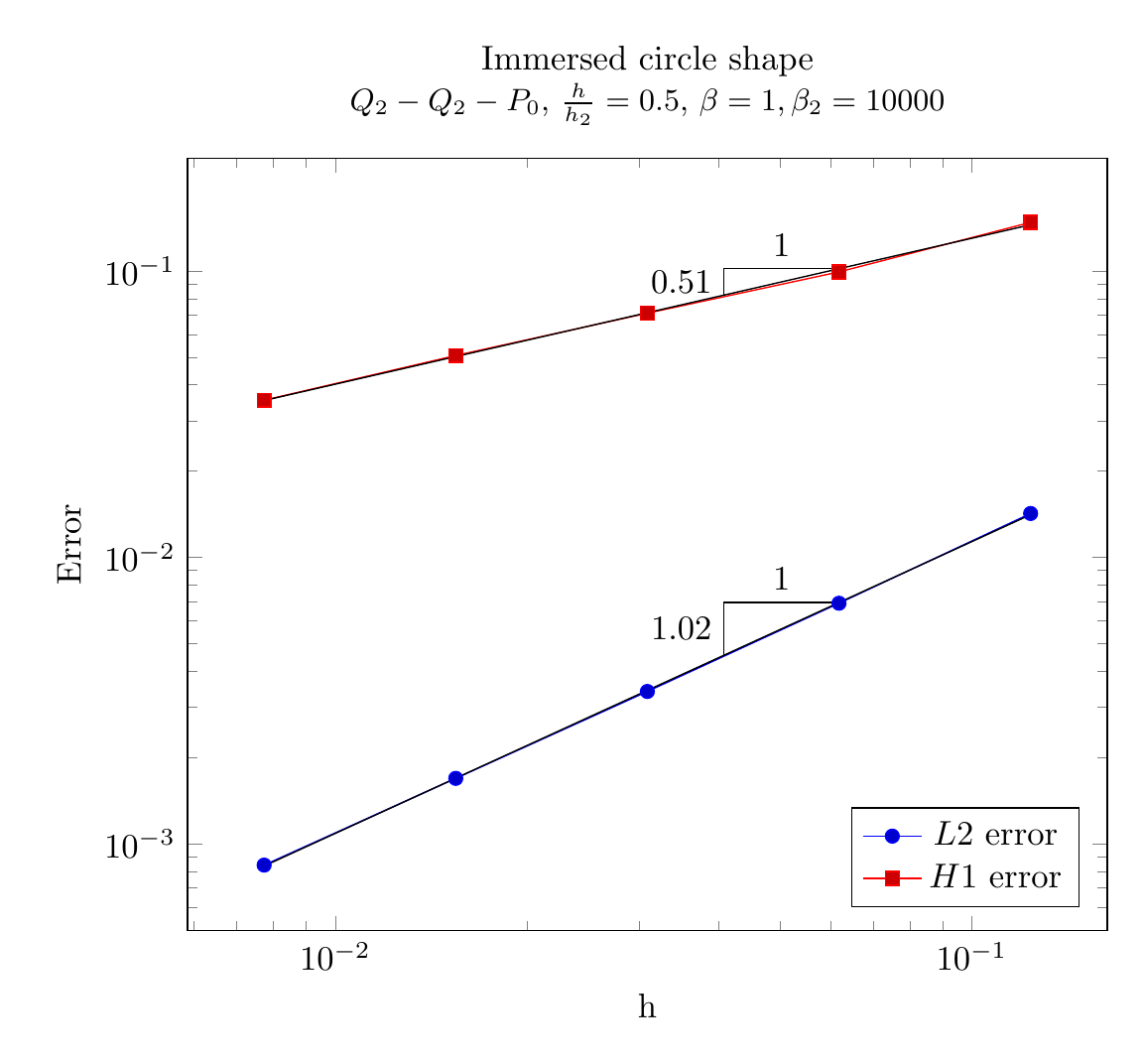} 
\end{minipage}
\begin{minipage}[c]{.3\linewidth}
		\includegraphics[width=1\linewidth]{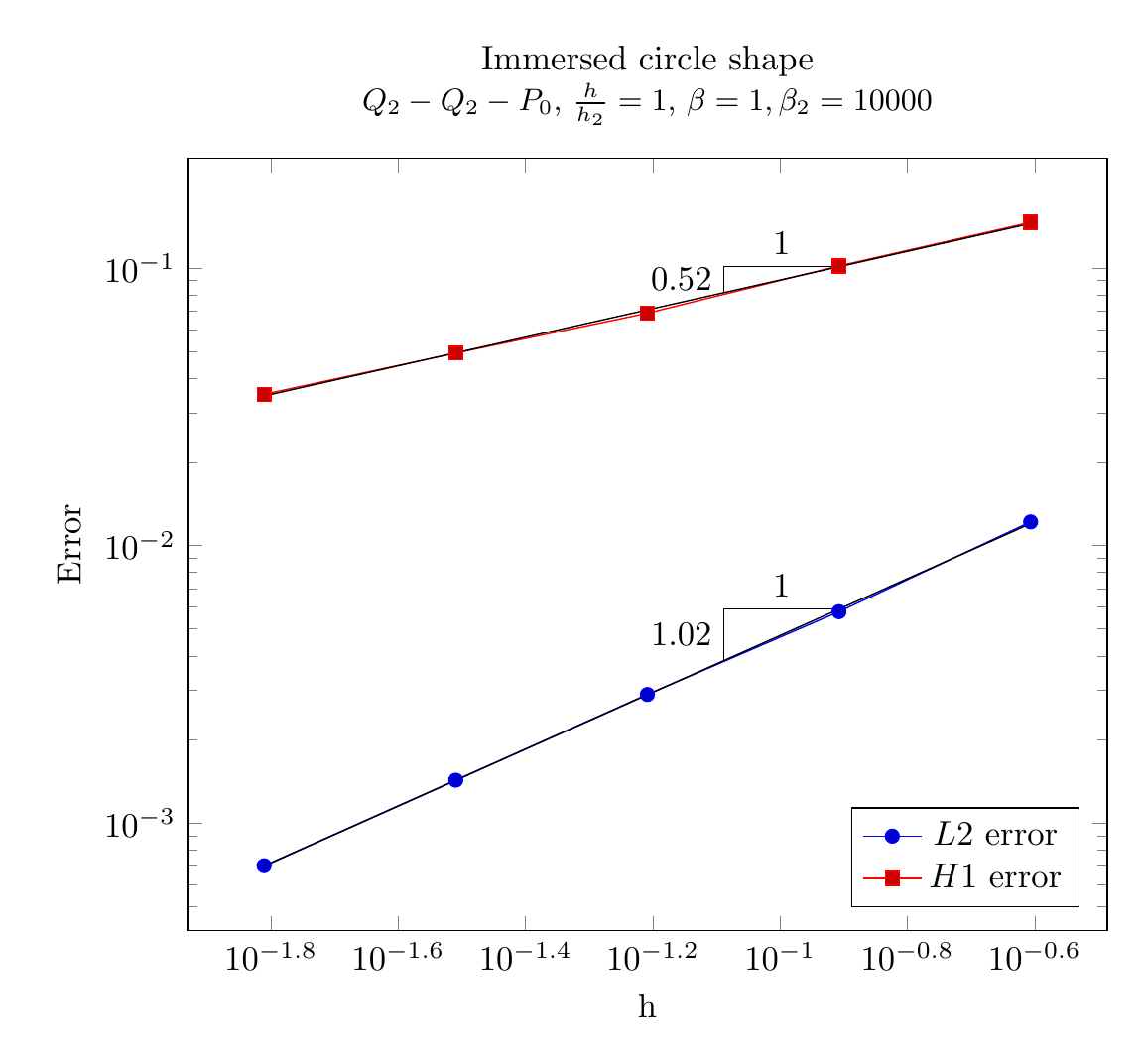} 
\end{minipage}
\begin{minipage}[c]{0.3\linewidth}
		\includegraphics[width=1\linewidth]{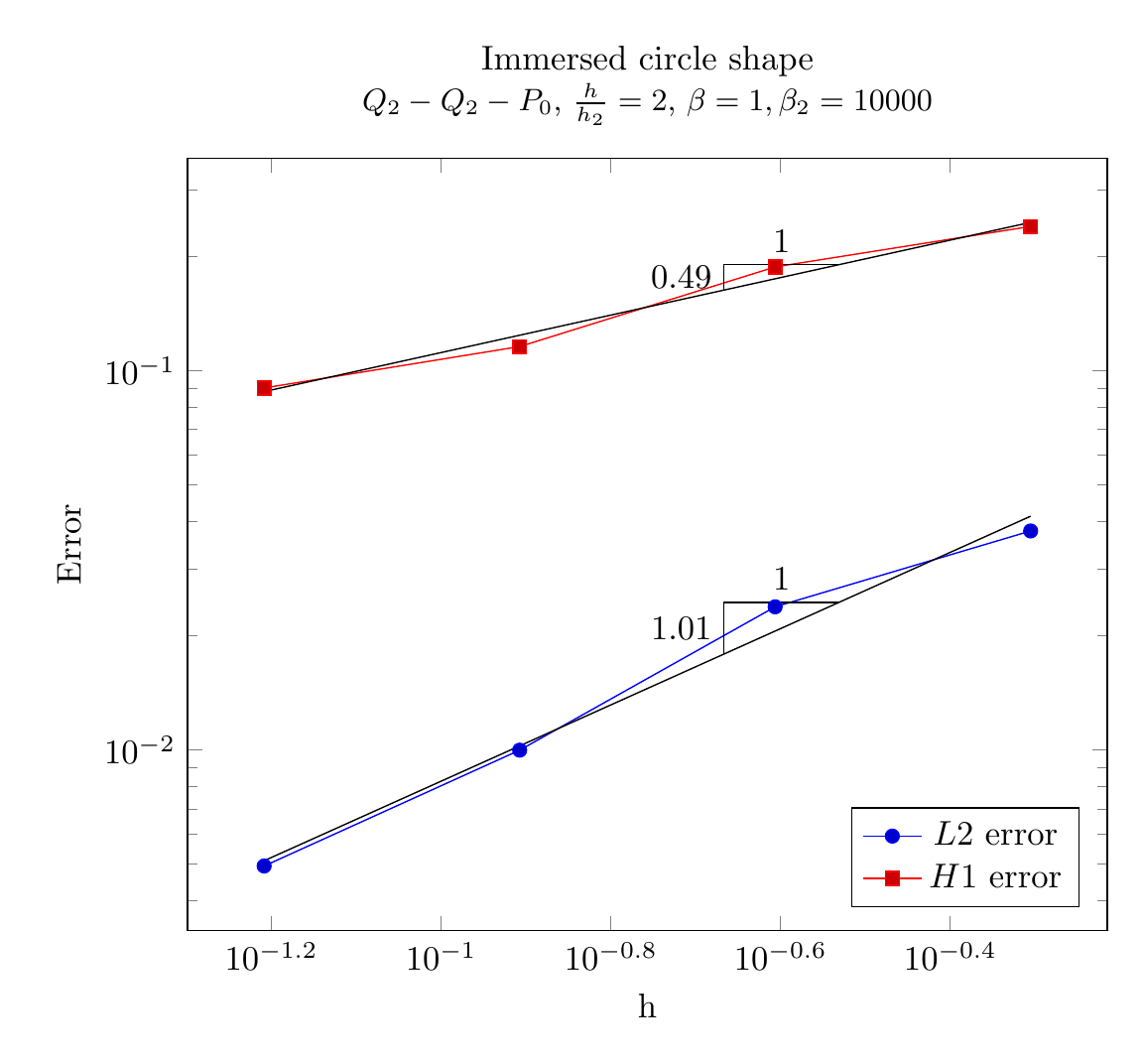} 
\end{minipage}
\captionof{figure}{Convergence of $Q_2-Q_2-P_0$ FDDLM: $\beta=1~\beta_2=10000$.}
  \label{fig:error220_110000_c}
\begin{minipage}[c]{0.3\linewidth}
		\includegraphics[width=1\linewidth]{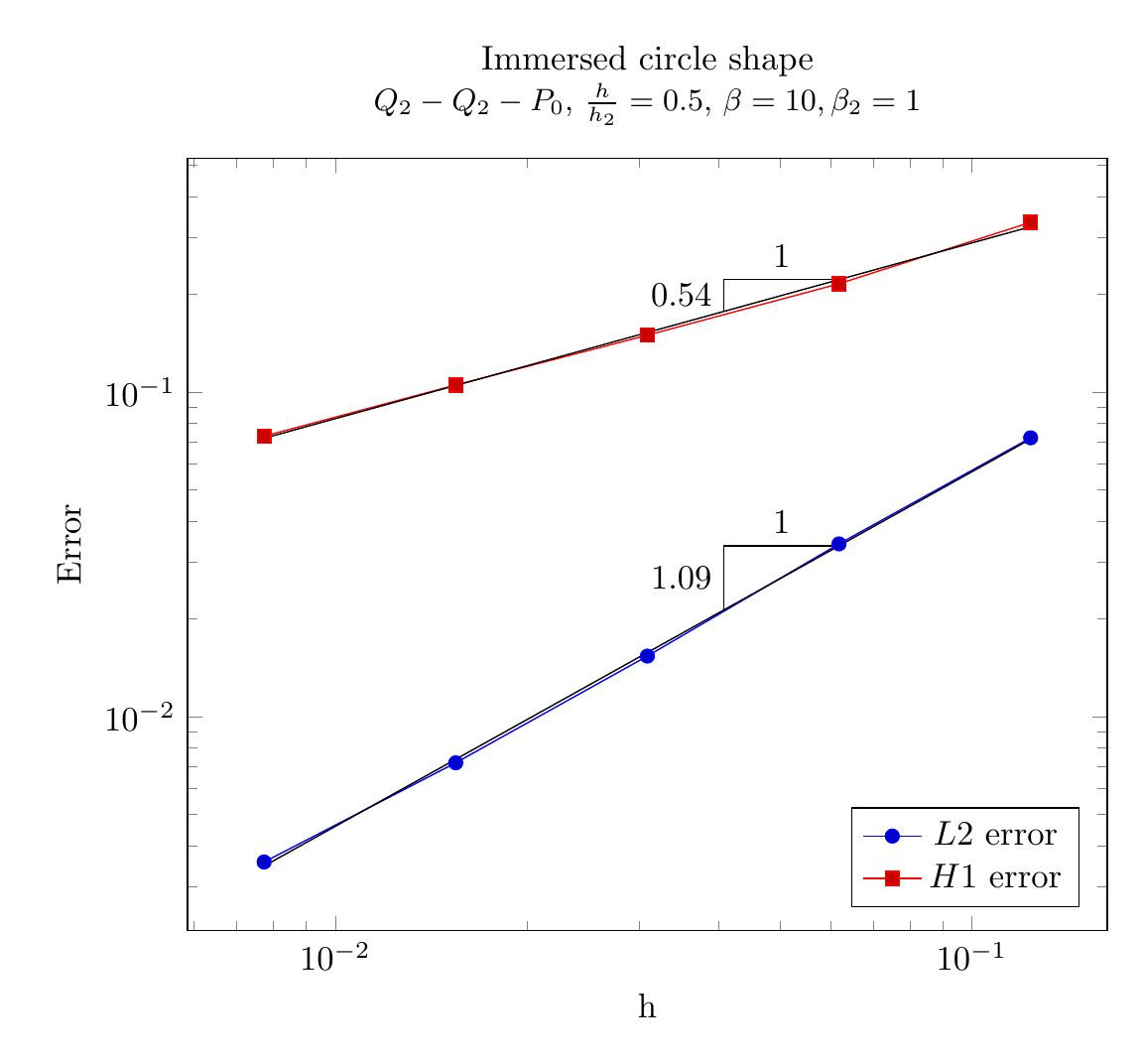} 
\end{minipage}
\begin{minipage}[c]{.3\linewidth}
		\includegraphics[width=1\linewidth]{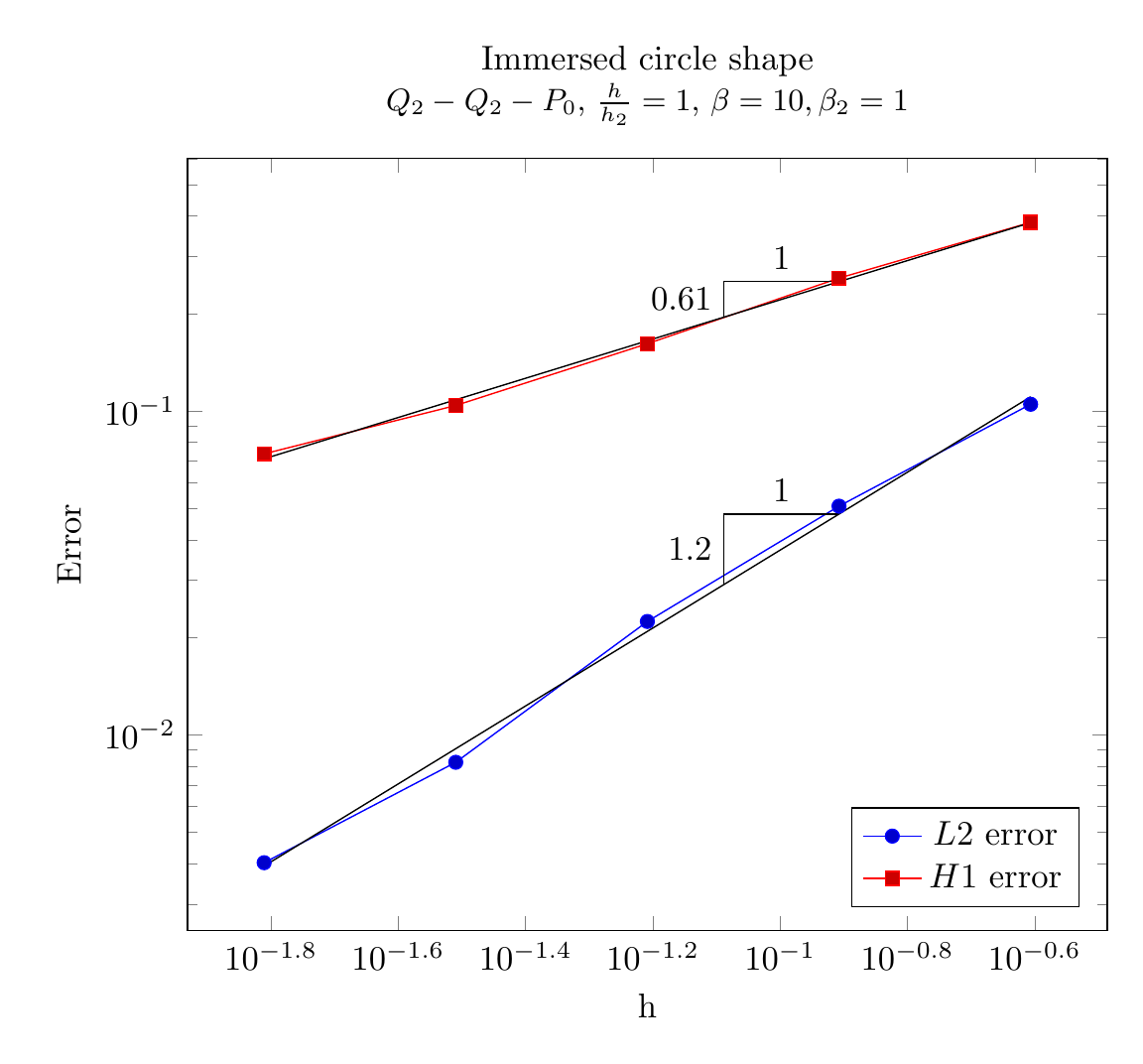} 
\end{minipage}
\begin{minipage}[c]{0.3\linewidth}
		\includegraphics[width=1\linewidth]{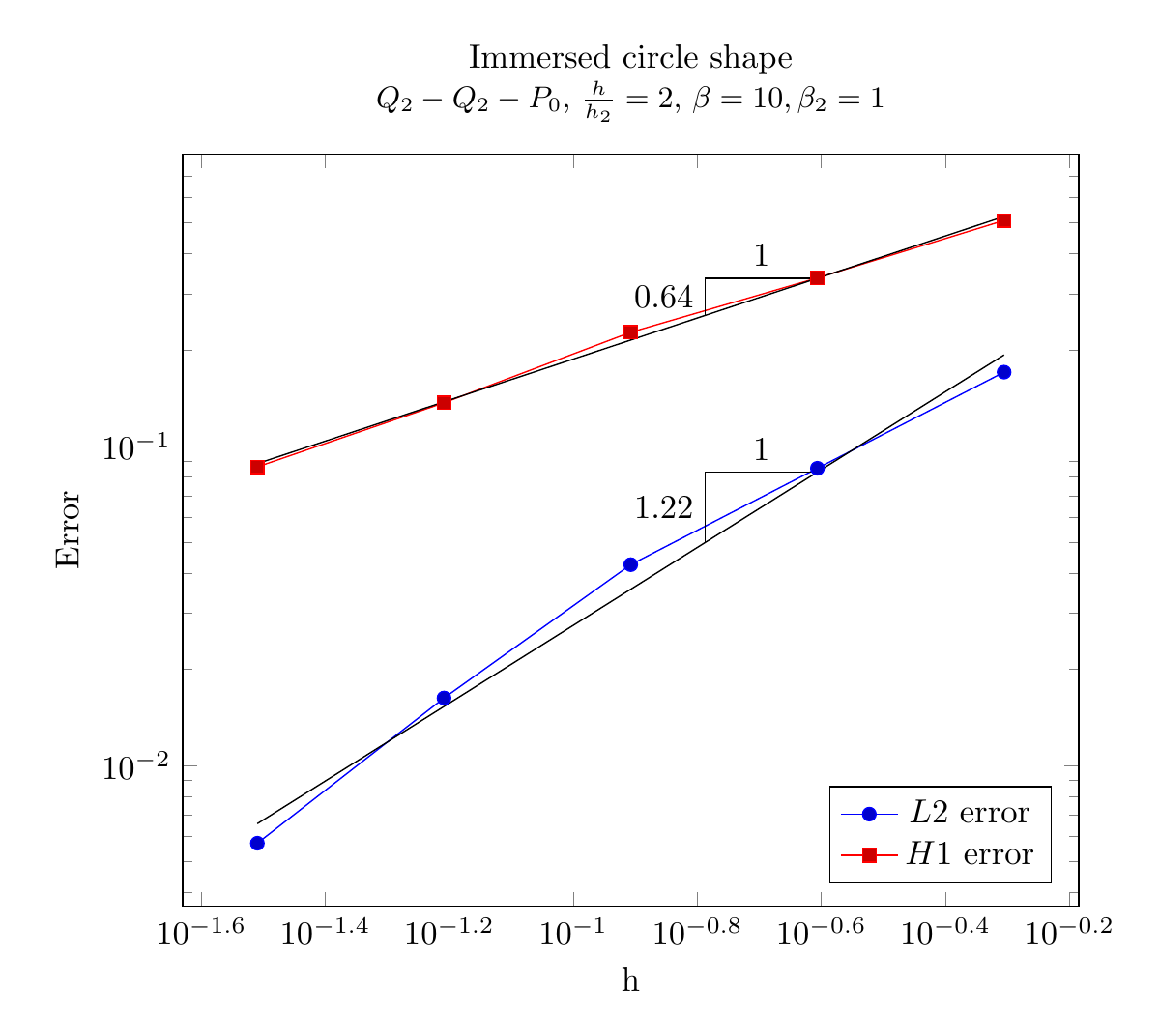} 
\end{minipage}
\captionof{figure}{Convergence of $Q_2-Q_2-P_0$ FDDLM: $\beta=10~\beta_2=1$.}
  \label{fig:error220_101_c}

\end{figure}

\begin{figure}[H]
\centering

\begin{minipage}[c]{.3\linewidth}
		\includegraphics[width=1\linewidth]{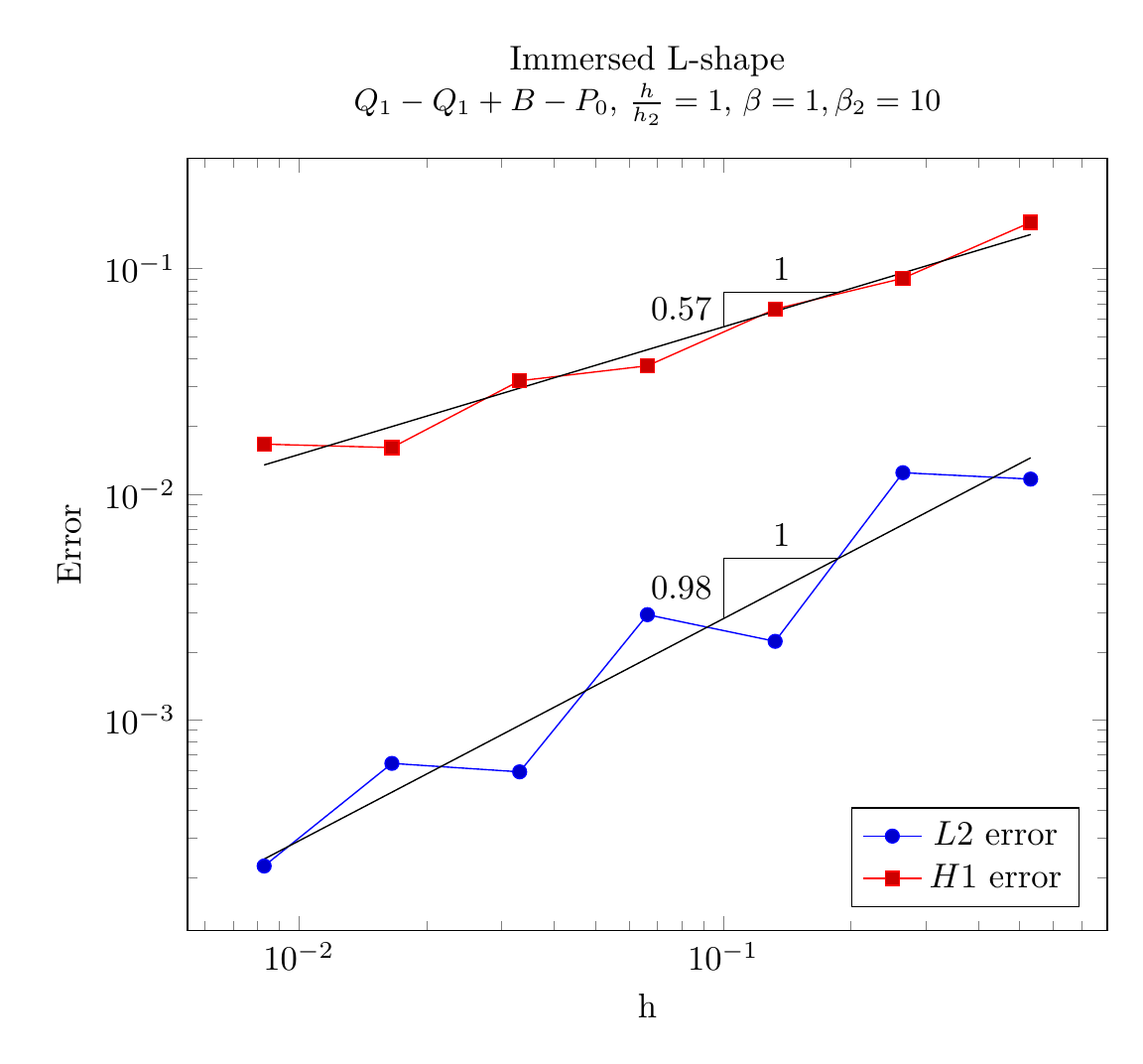} 
\end{minipage}
\begin{minipage}[c]{.3\linewidth}
		\includegraphics[width=1\linewidth]{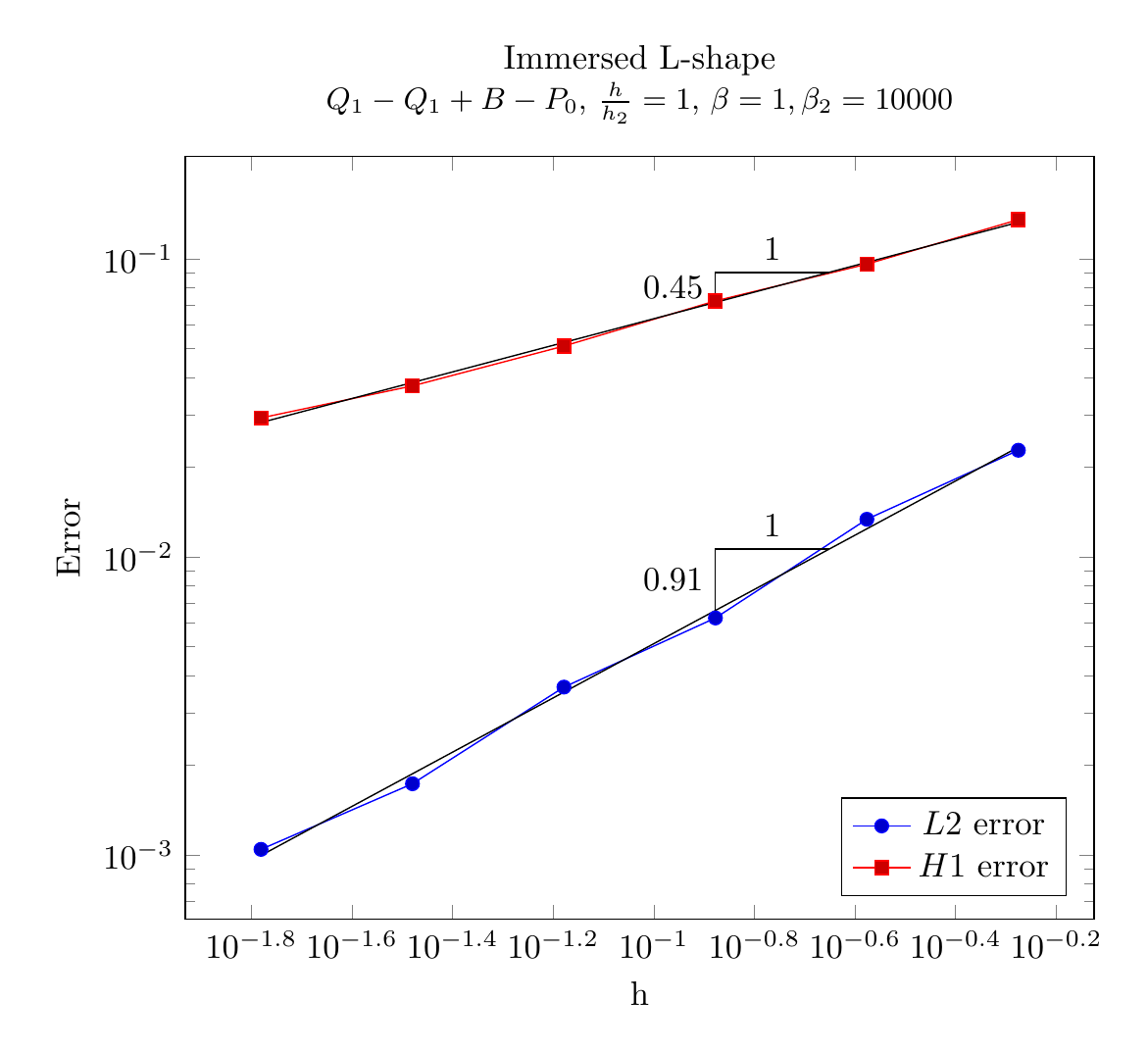} 
\end{minipage}
\begin{minipage}[c]{.3\linewidth}
		\includegraphics[width=1\linewidth]{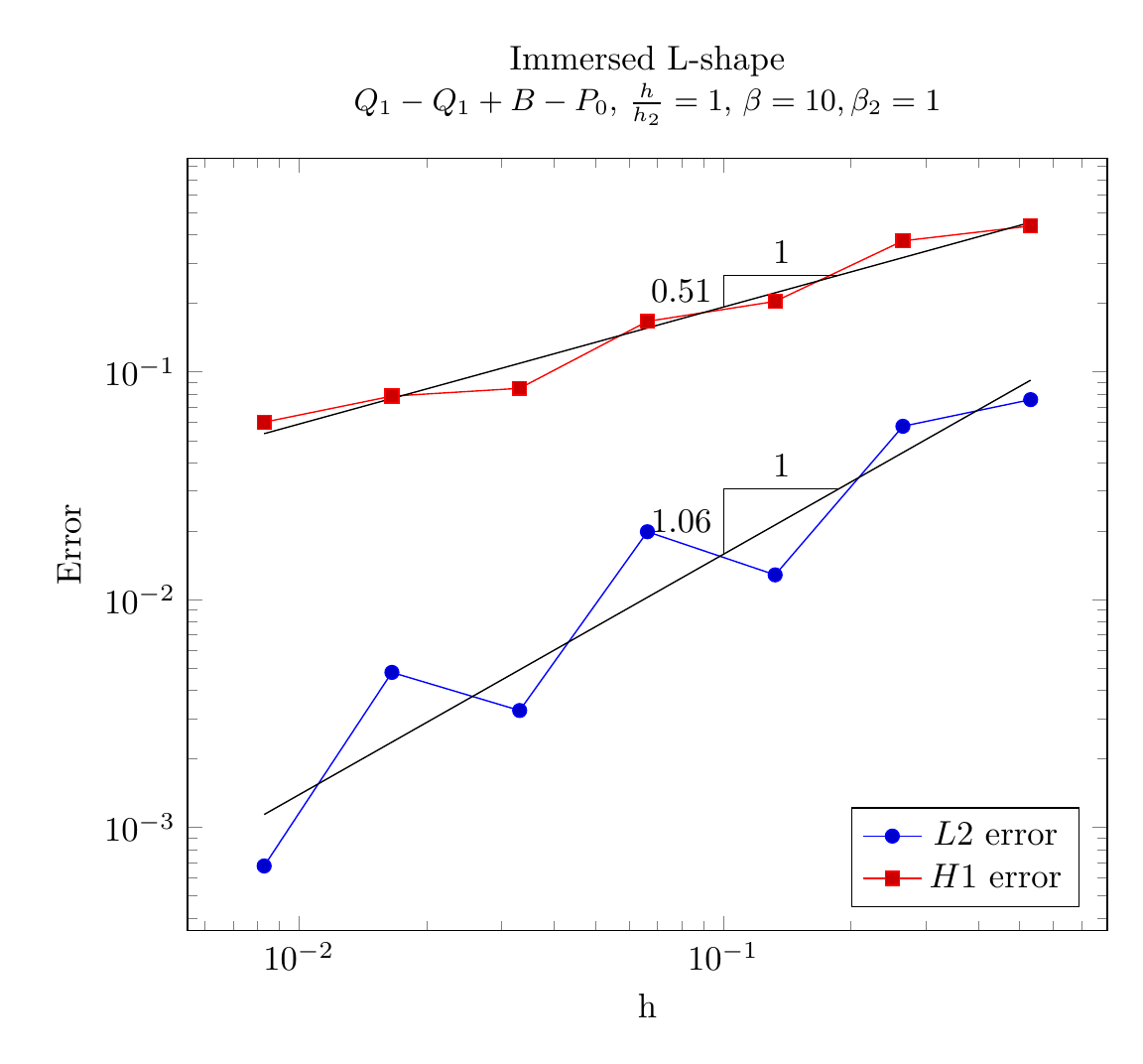} 
\end{minipage}
\captionof{figure}{Convergence of $Q_1-(Q_1+B)-P_0$ FDDLM: immersed L-shape.}
  \label{fig:error1_1_0_L}
  
  \begin{minipage}[c]{.3\linewidth}
		\includegraphics[width=1\linewidth]{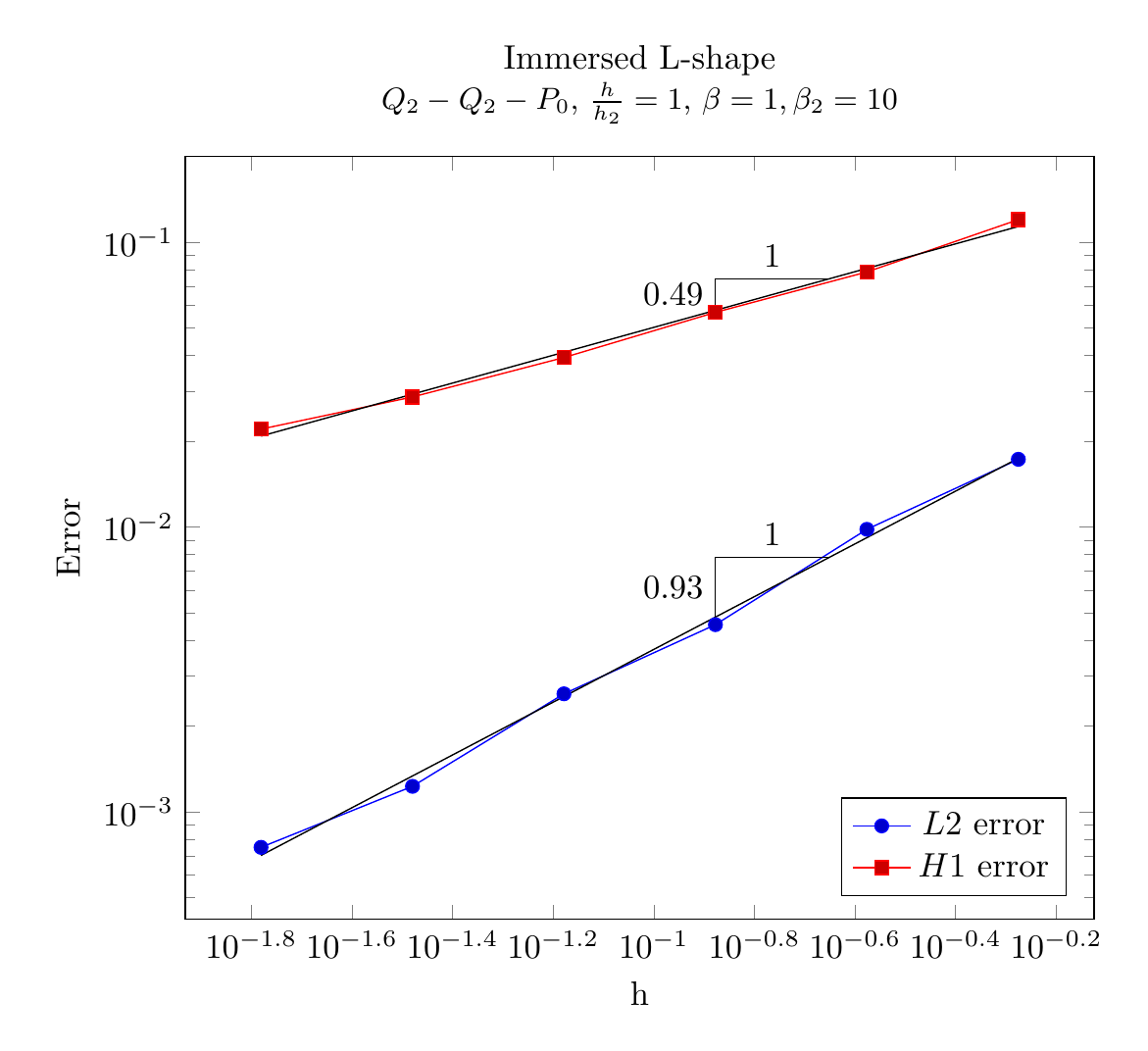} 
\end{minipage}
\begin{minipage}[c]{.3\linewidth}
		\includegraphics[width=1\linewidth]{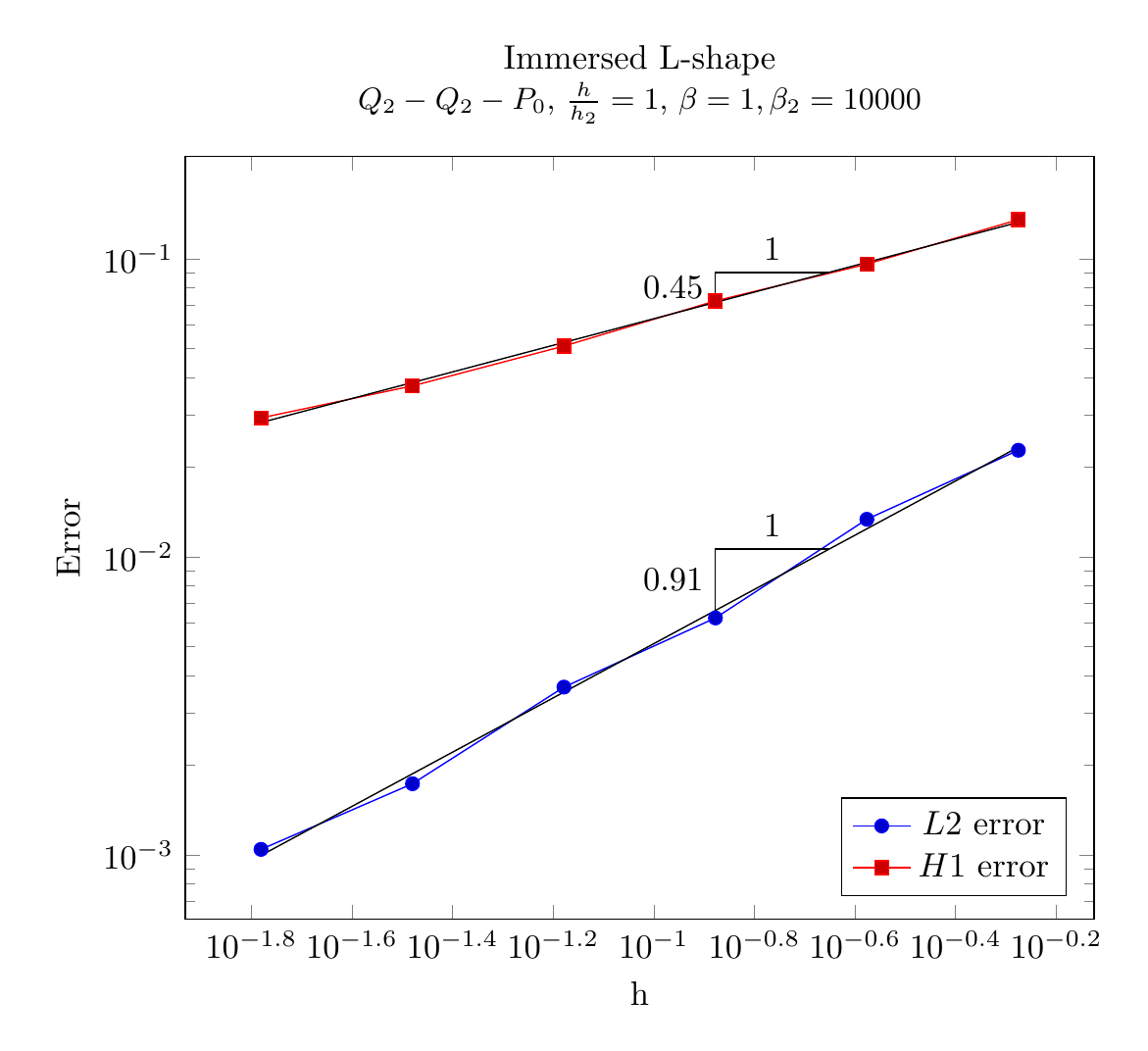} 
\end{minipage}
\begin{minipage}[c]{.3\linewidth}
		\includegraphics[width=1\linewidth]{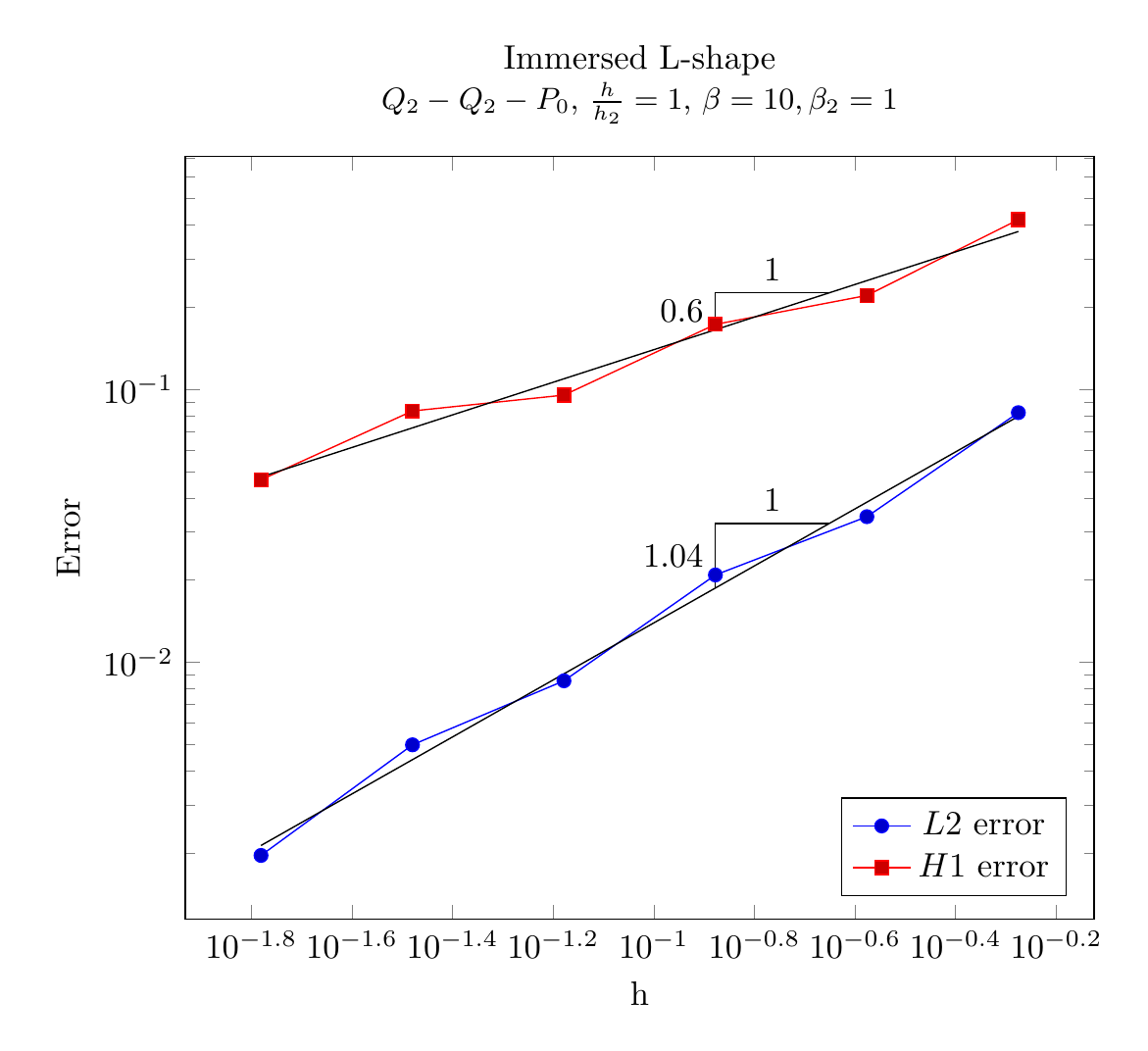} 
\end{minipage}
\captionof{figure}{Convergence of $Q_2-Q_2-P_0$ FDDLM: immersed L-shape.}
  \label{fig:error220_L}
  
   \begin{minipage}[c]{.3\linewidth}
		\includegraphics[width=1\linewidth]{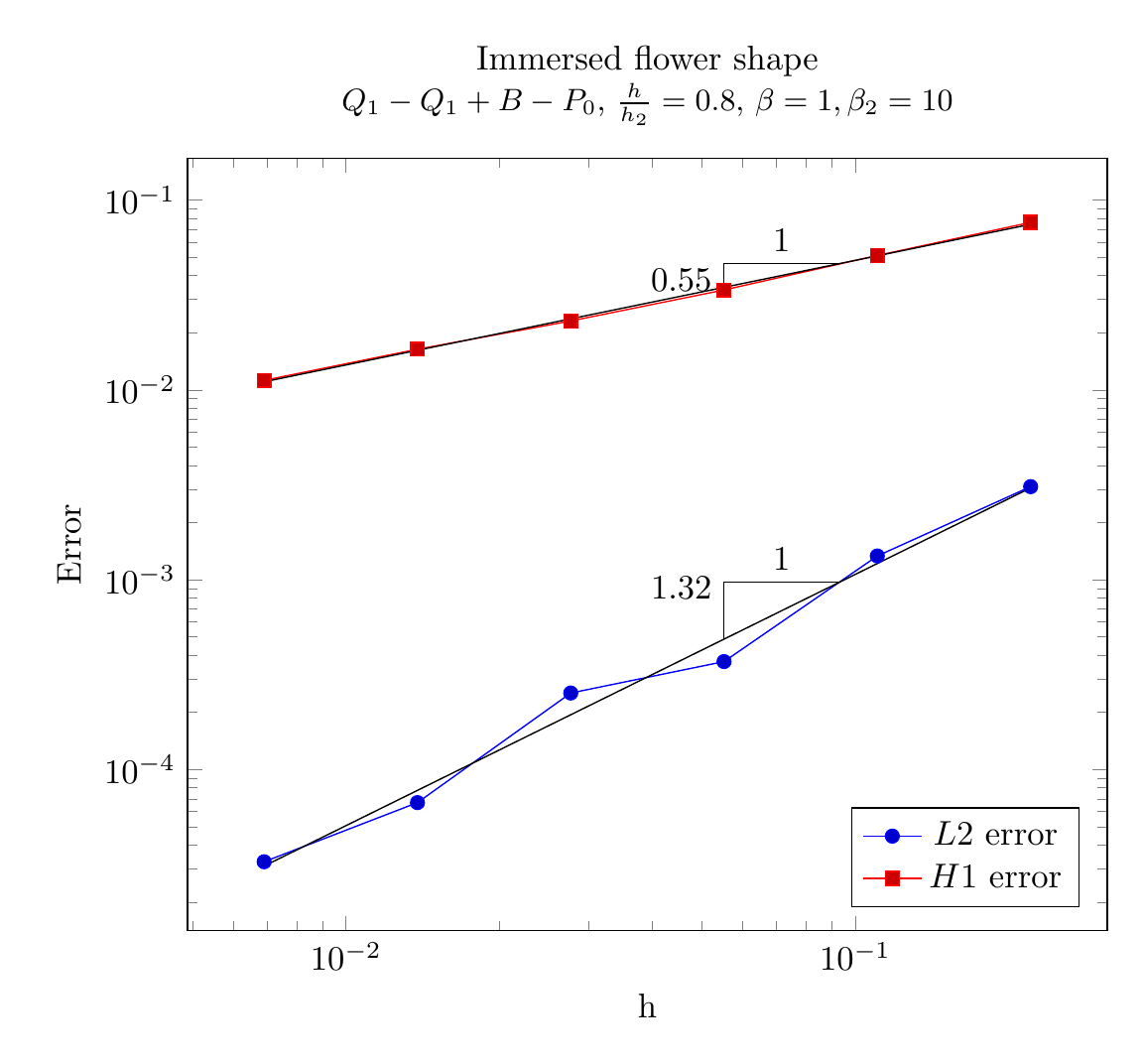} 
\end{minipage}
\begin{minipage}[c]{.3\linewidth}
		\includegraphics[width=1\linewidth]{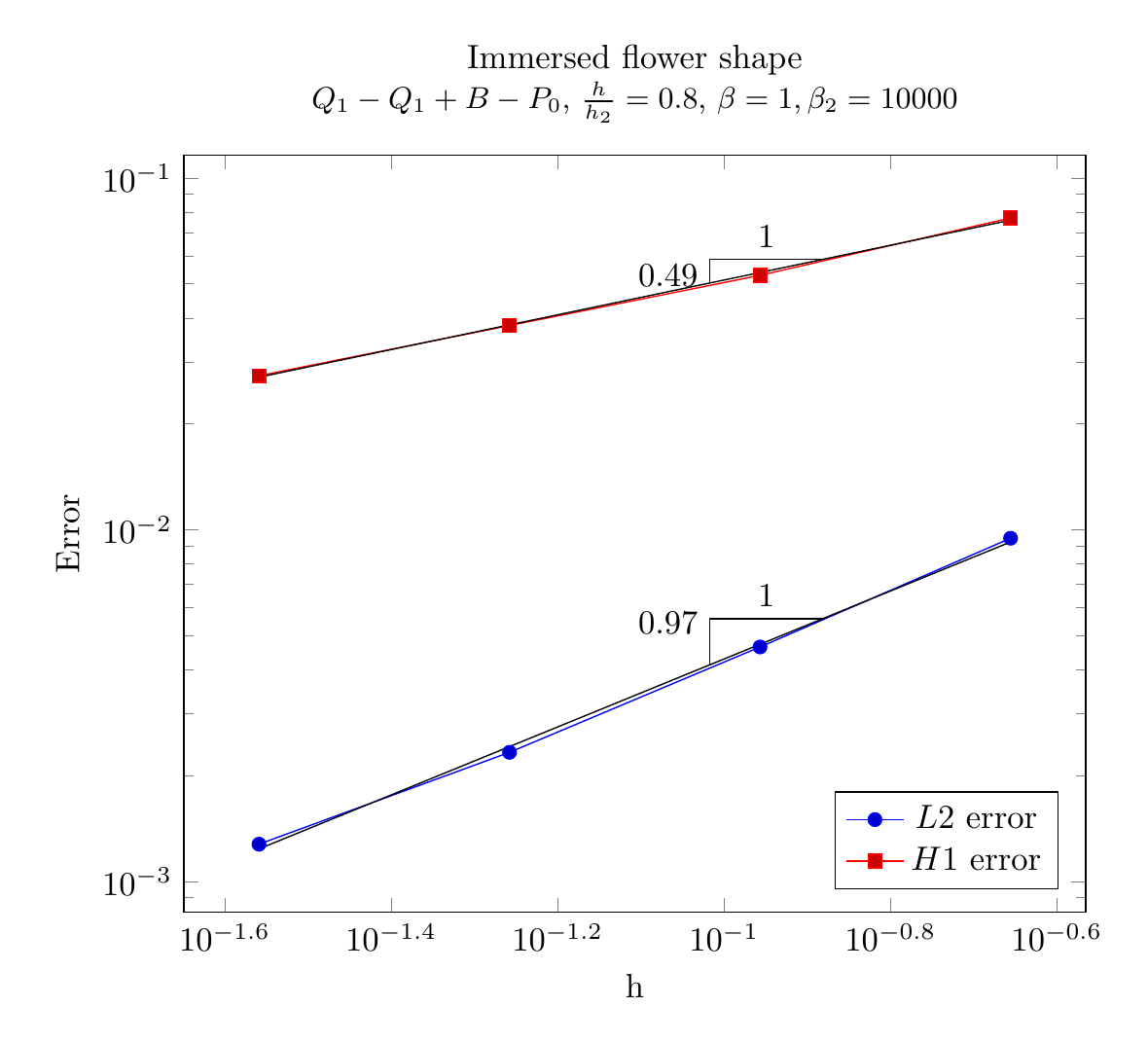} 
\end{minipage}
\begin{minipage}[c]{.3\linewidth}
		\includegraphics[width=1\linewidth]{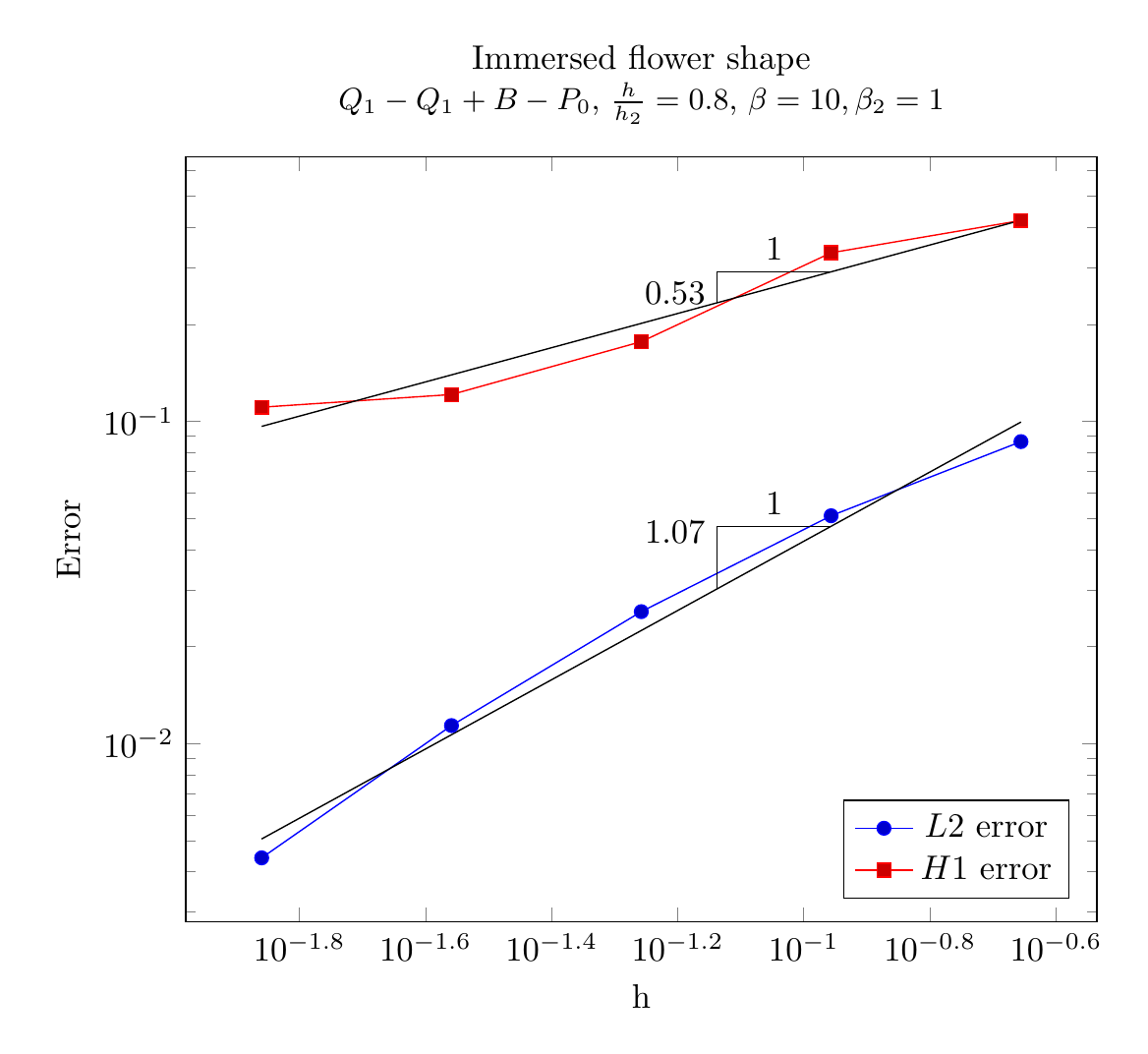} 
\end{minipage}
\captionof{figure}{Convergence of $Q_1-(Q_1+B)-P_0$ FDDLM: immersed flower shape.}
  \label{fig:error1_1_0_f}
  
  \begin{minipage}[c]{.3\linewidth}
		\includegraphics[width=1\linewidth]{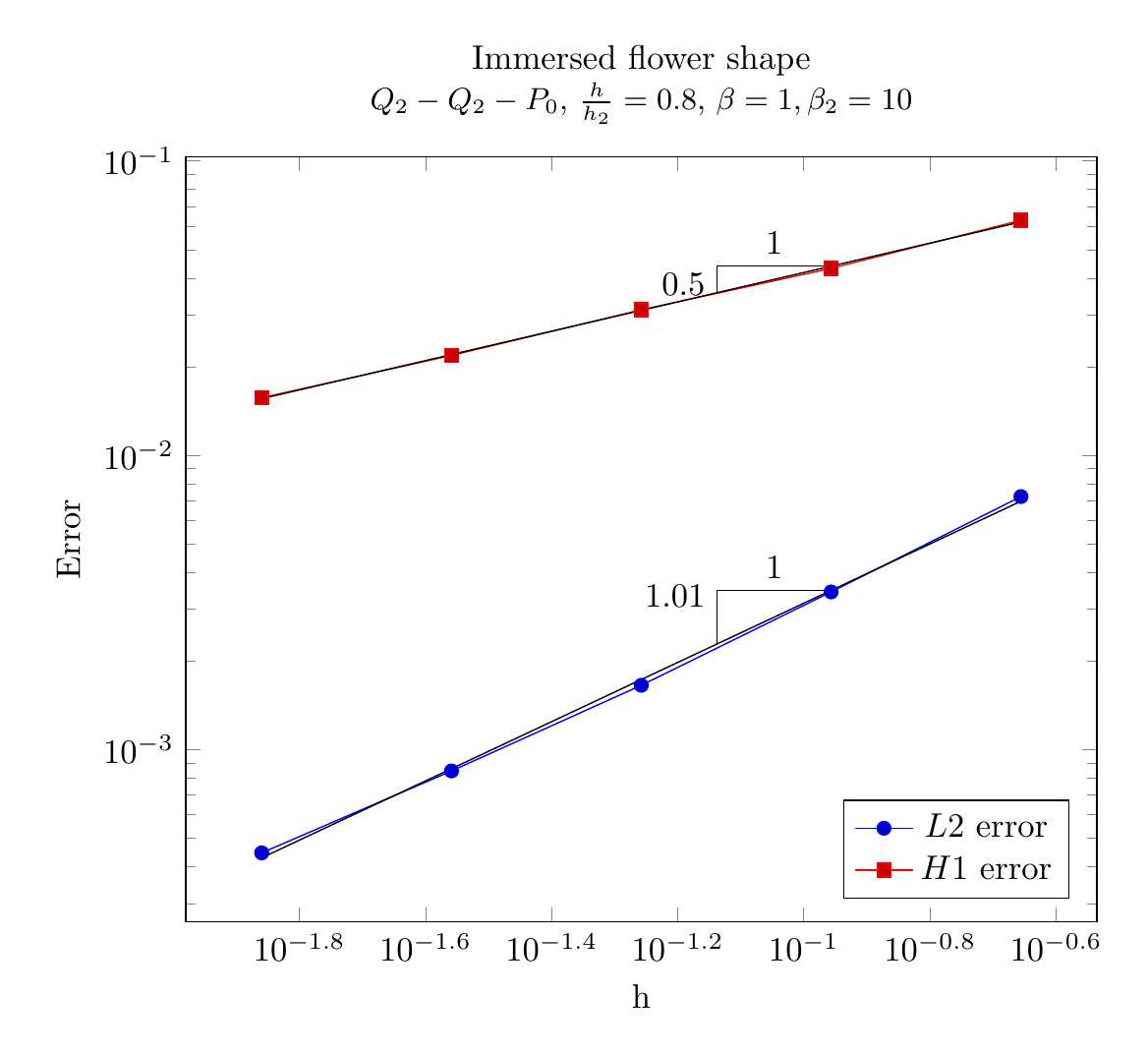} 
\end{minipage}
\begin{minipage}[c]{.3\linewidth}
		\includegraphics[width=1\linewidth]{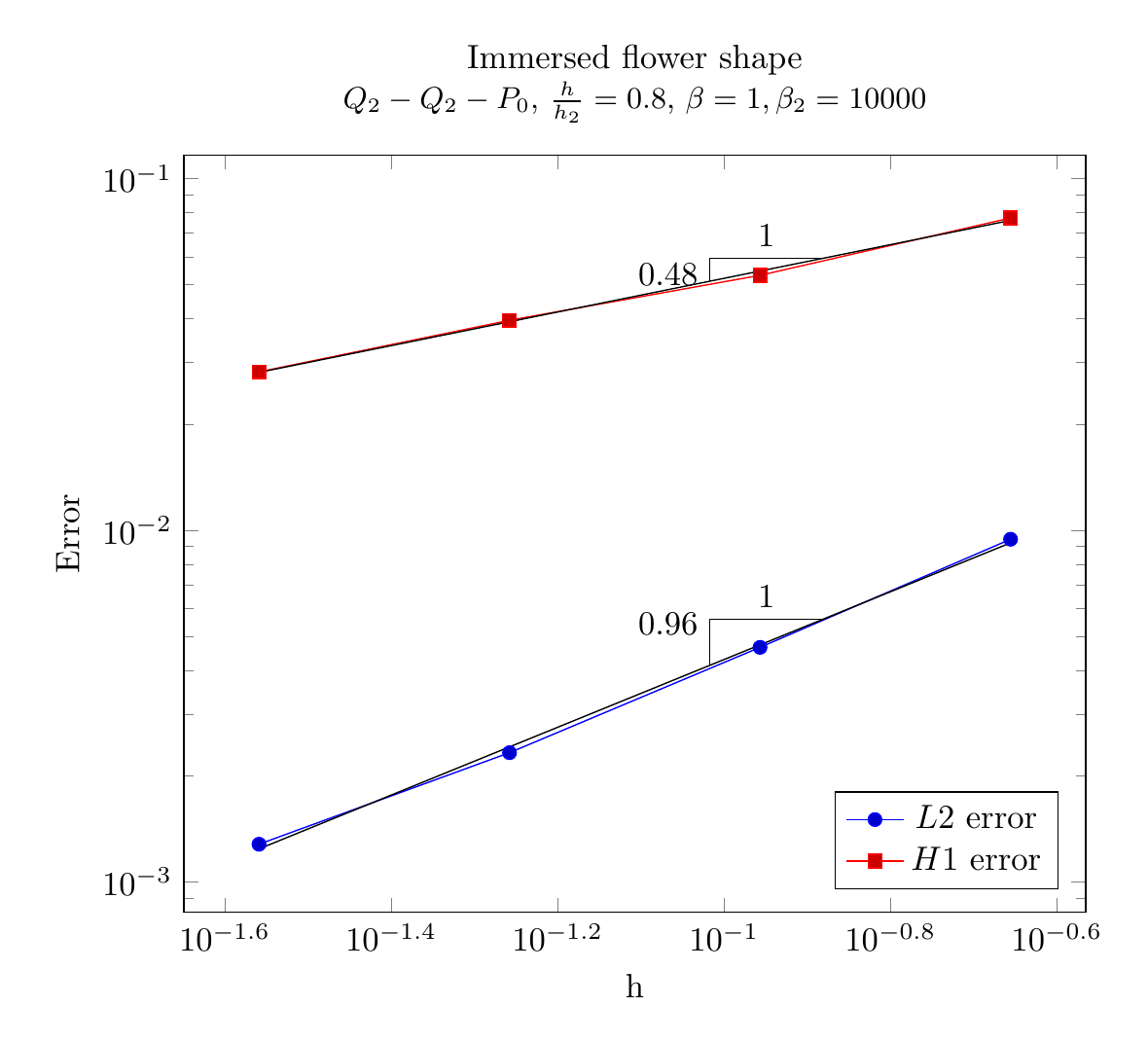} 
\end{minipage}
\begin{minipage}[c]{.3\linewidth}
		\includegraphics[width=1\linewidth]{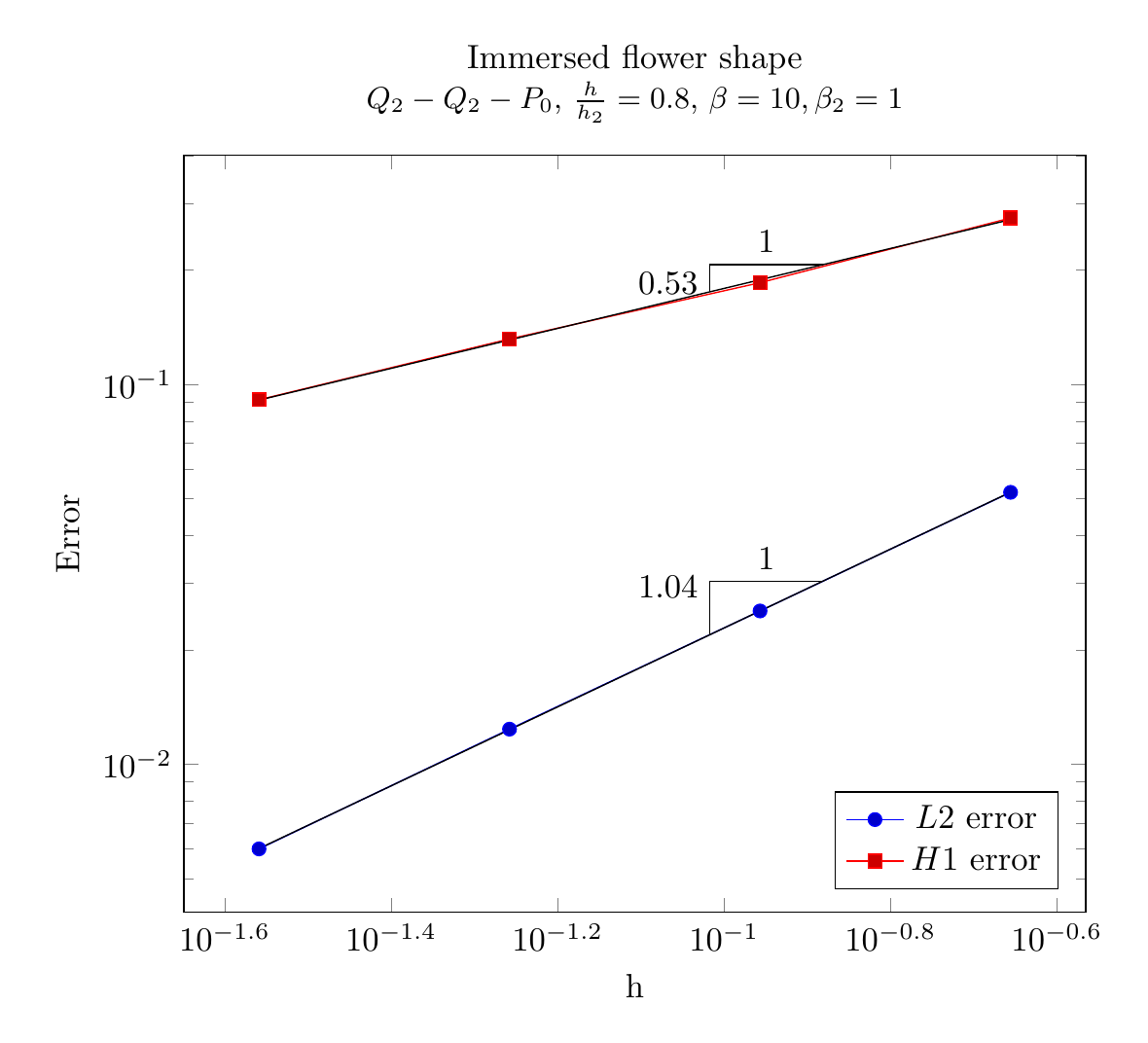} 
\end{minipage}
\captionof{figure}{Convergence of $Q_2-Q_2-P_0$ FDDLM: immersed flower shape.}
  \label{fig:error220_f}
\end{figure}

\begin{figure}[H]
\centering
\begin{minipage}[c]{.3\linewidth}
		\includegraphics[width=1\linewidth]{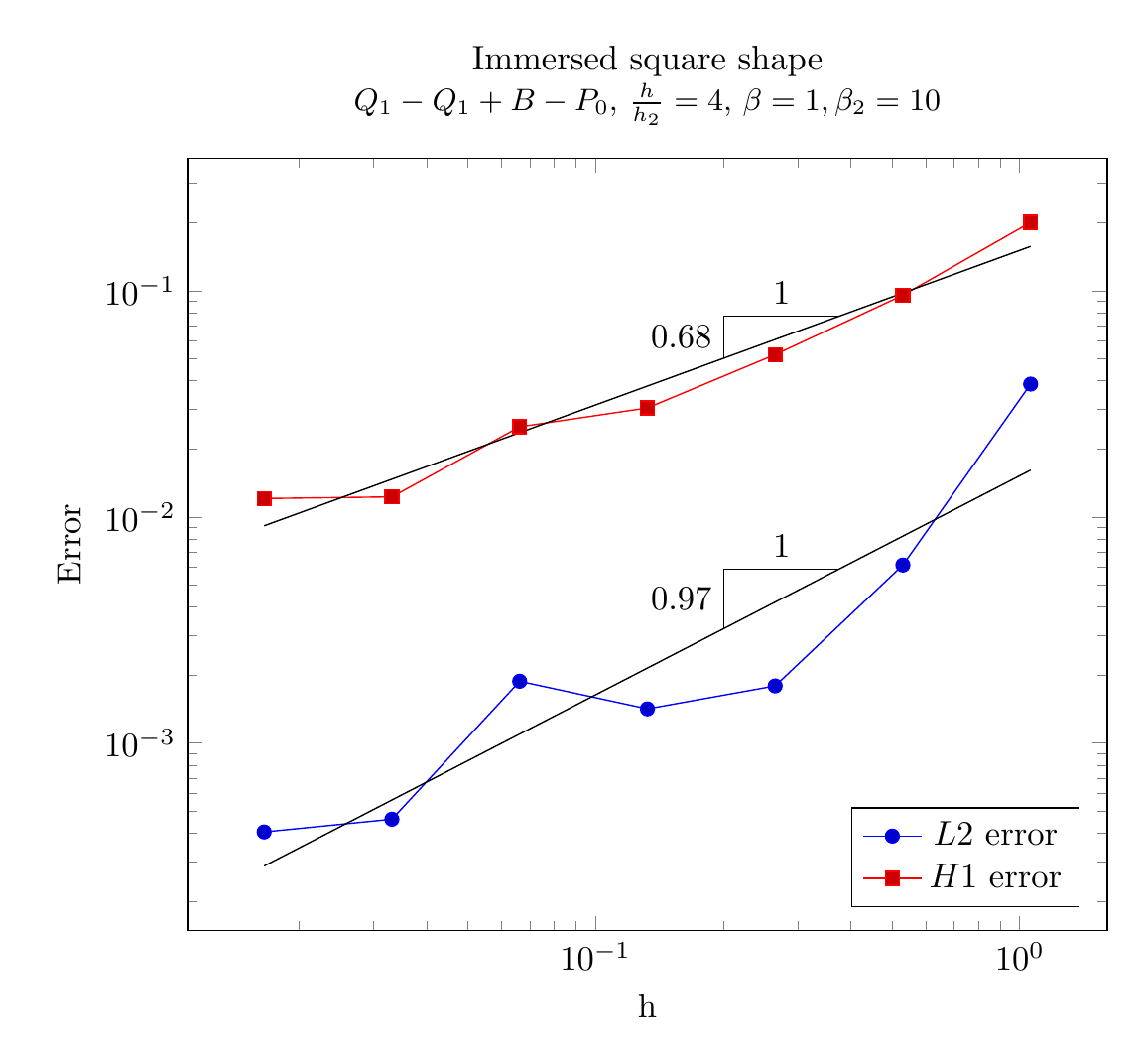} 
\end{minipage}
\begin{minipage}[c]{.3\linewidth}
		\includegraphics[width=1\linewidth]{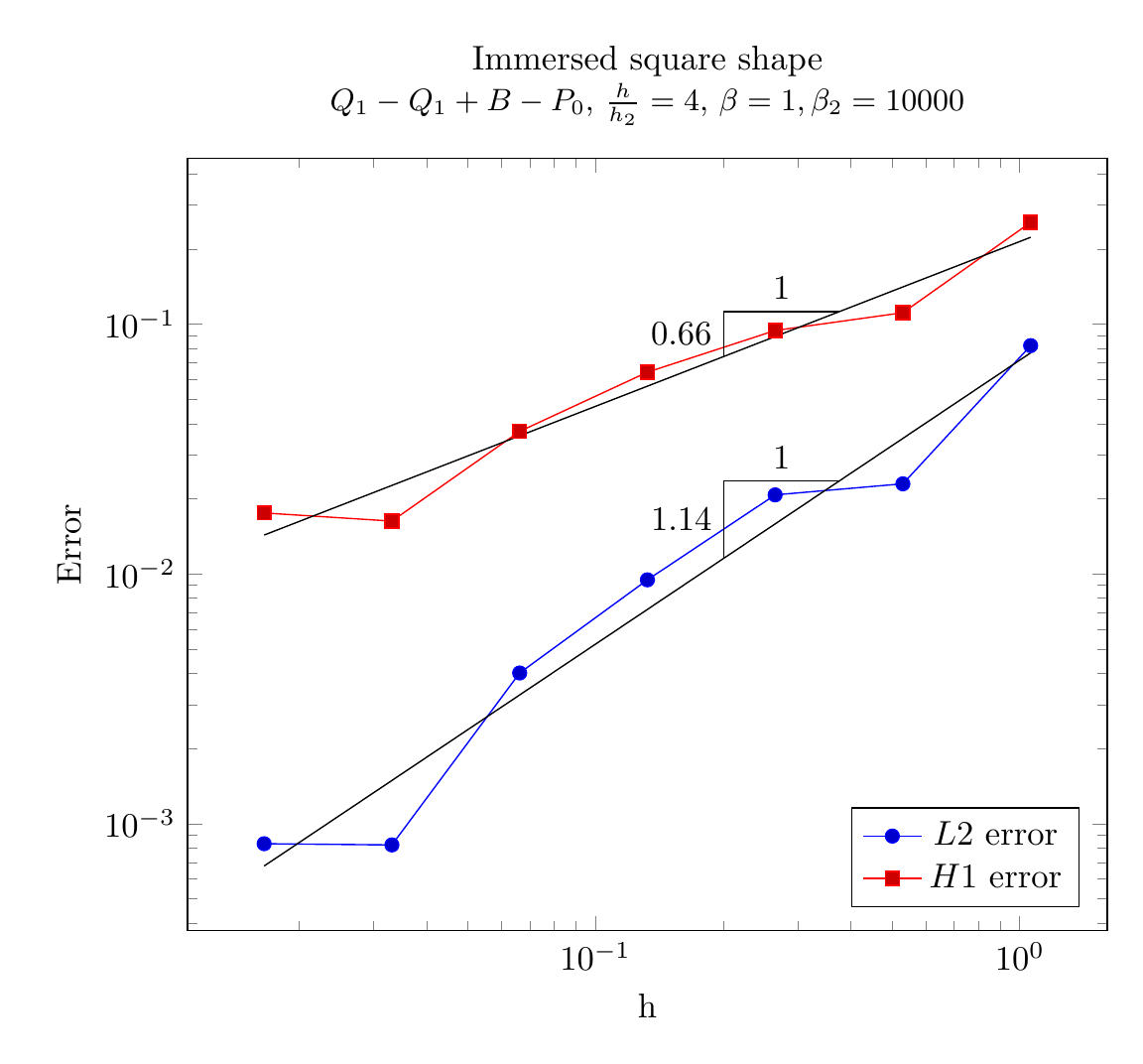} 
\end{minipage}
\begin{minipage}[c]{.3\linewidth}
		\includegraphics[width=1\linewidth]{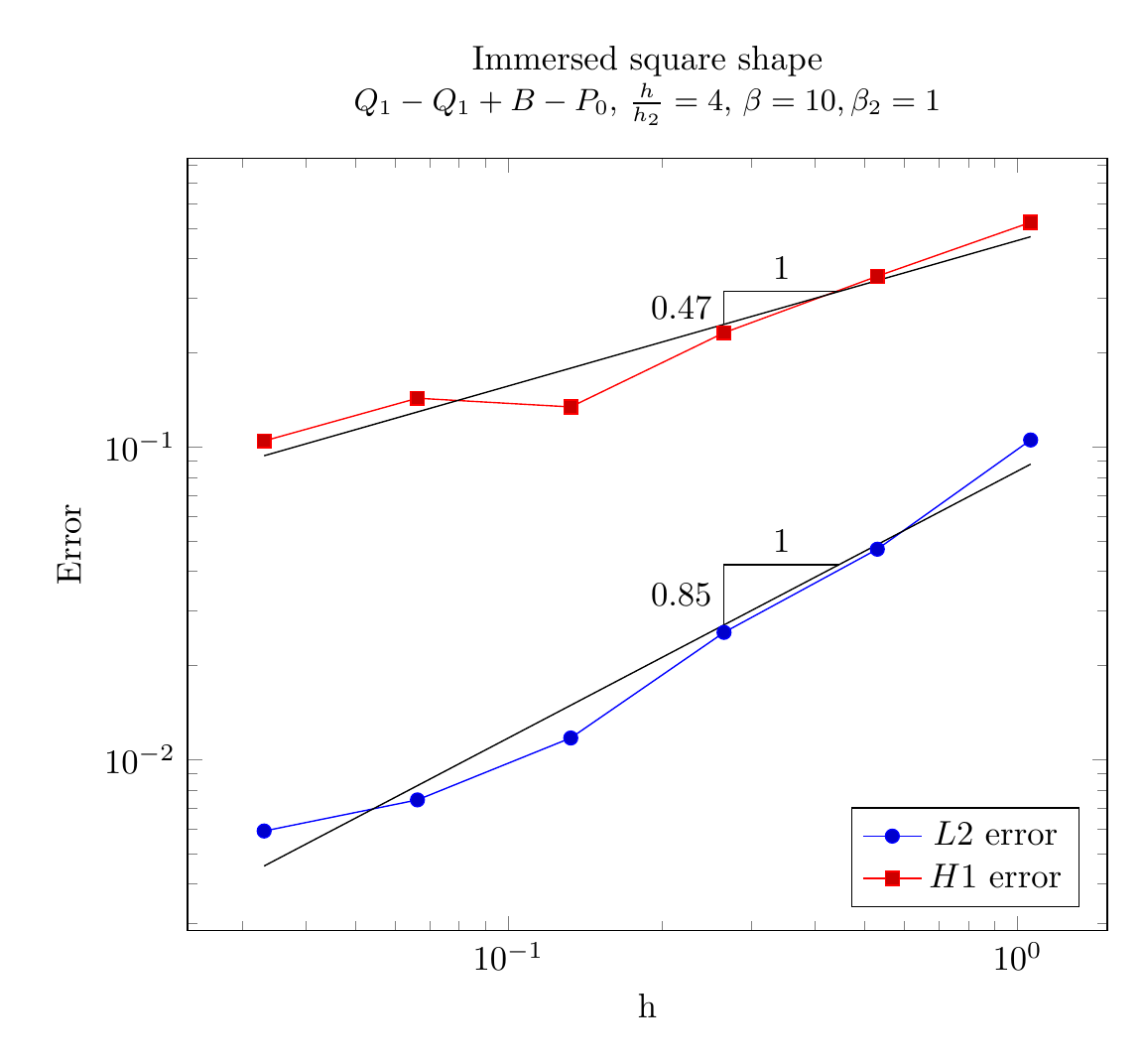} 
\end{minipage}
\captionof{figure}{Convergence of $Q_1-(Q_1+B)-P_0$FDDLM-FEM: immersed square shape.}
  \label{fig:error1_1_0_s}
  
  \begin{minipage}[c]{.3\linewidth}
		\includegraphics[width=1\linewidth]{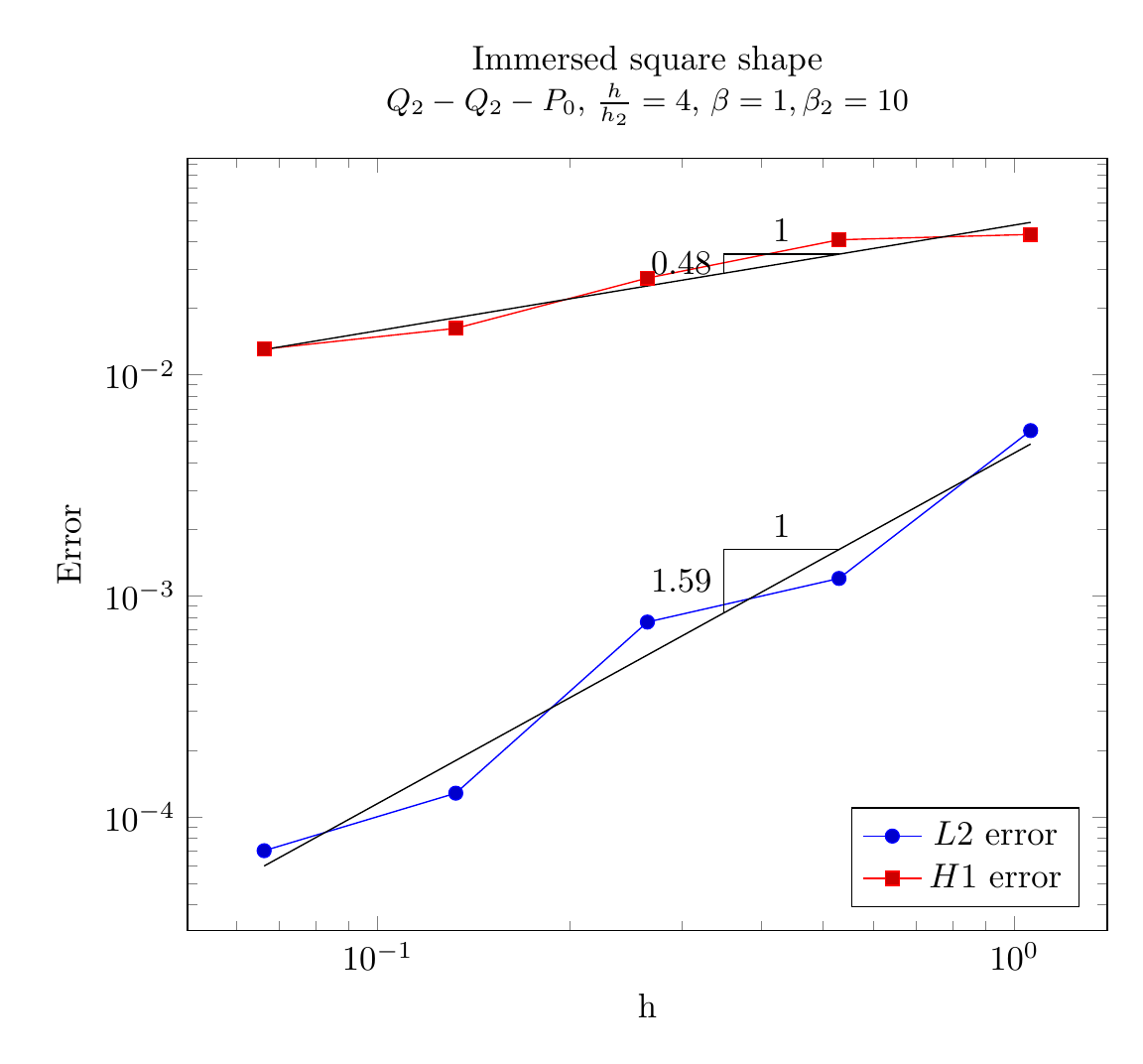} 
\end{minipage}
\begin{minipage}[c]{.3\linewidth}
		\includegraphics[width=1\linewidth]{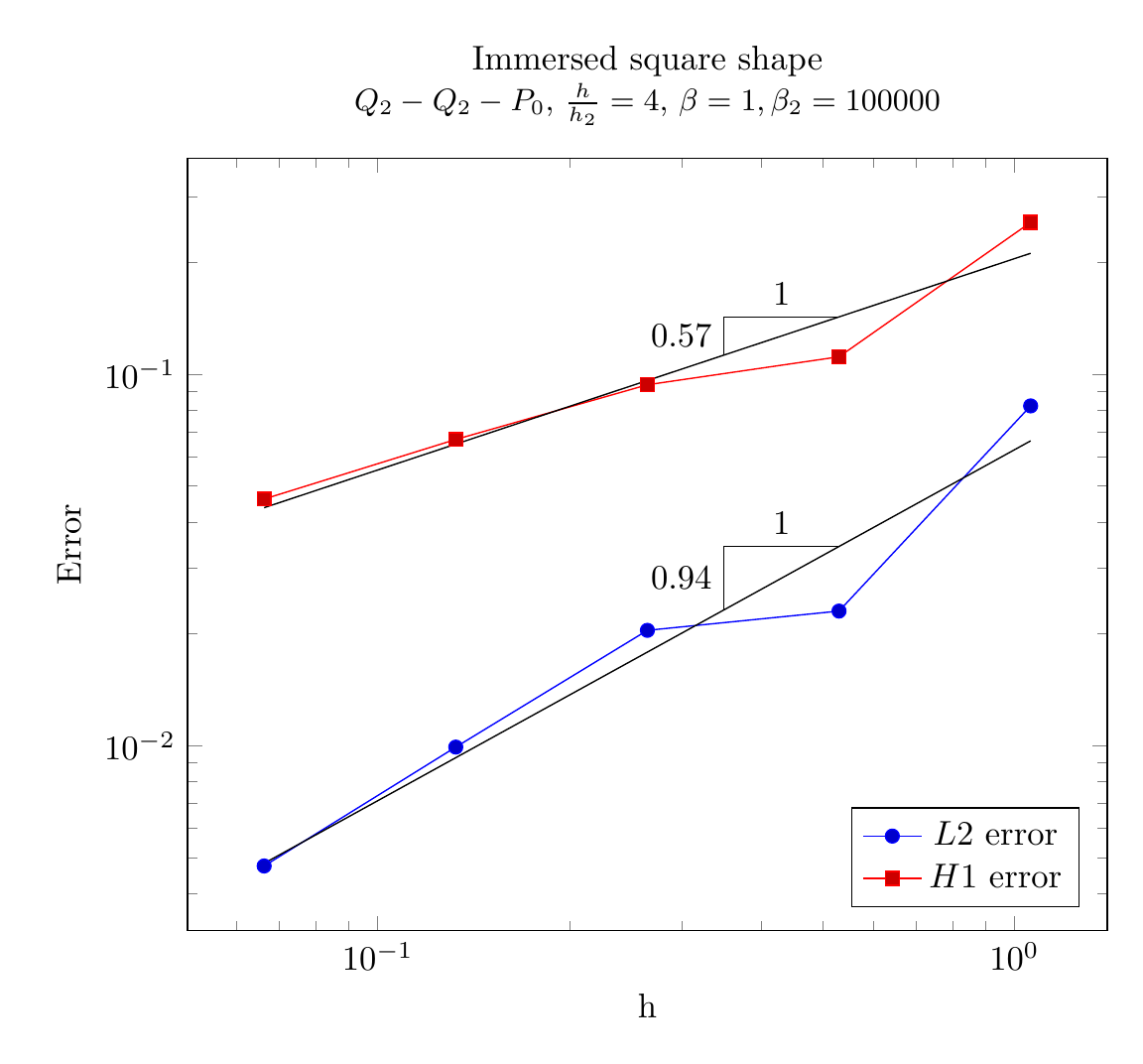} 
\end{minipage}
\begin{minipage}[c]{.3\linewidth}
		\includegraphics[width=1\linewidth]{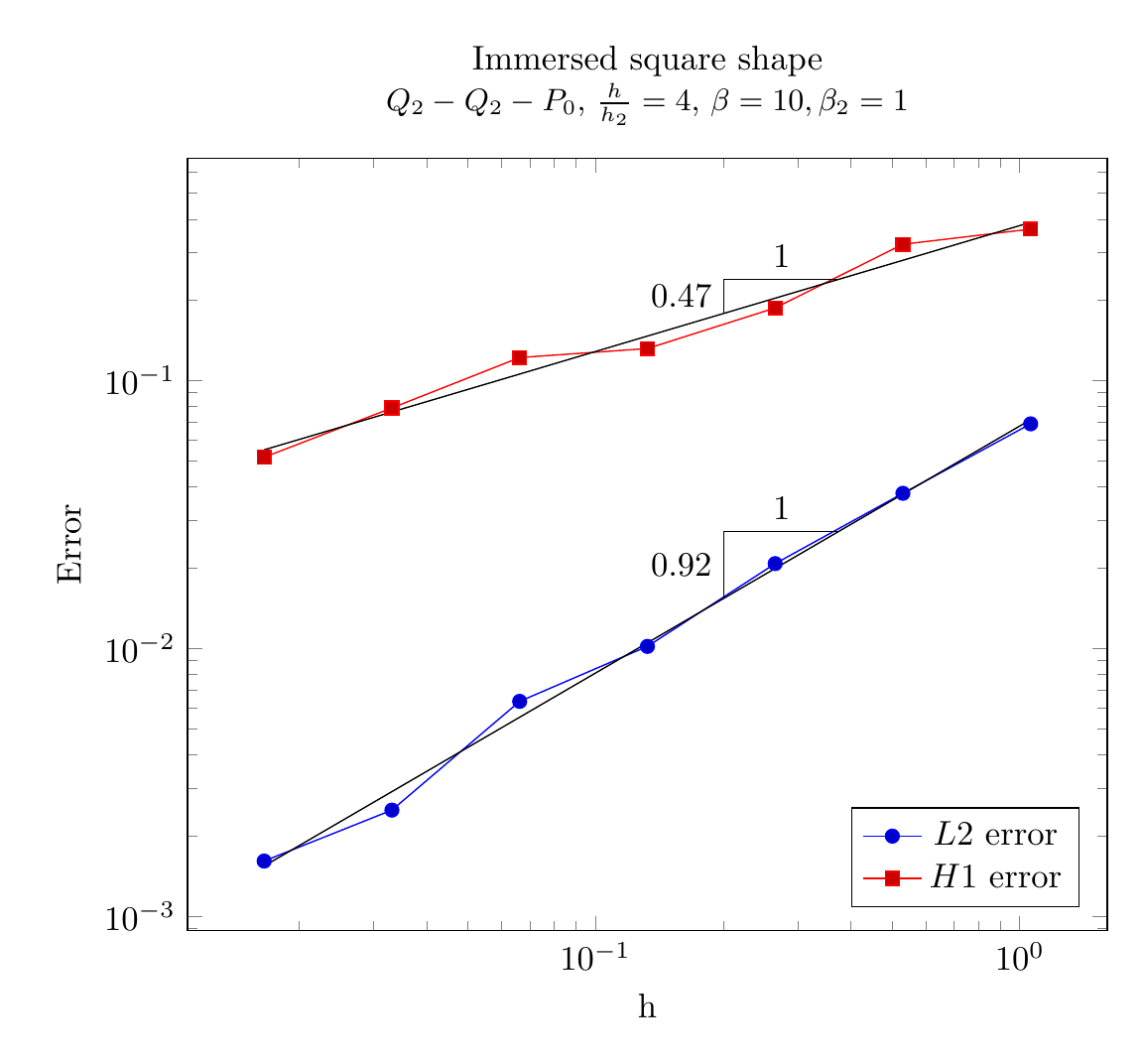} 
\end{minipage}
\captionof{figure}{Convergence of $Q_2-Q_2-P_0$ FDDLM: Immersed Square Shape}
  \label{fig:error220_s}
\begin{minipage}[c]{.3\linewidth}
		\includegraphics[width=1\linewidth]{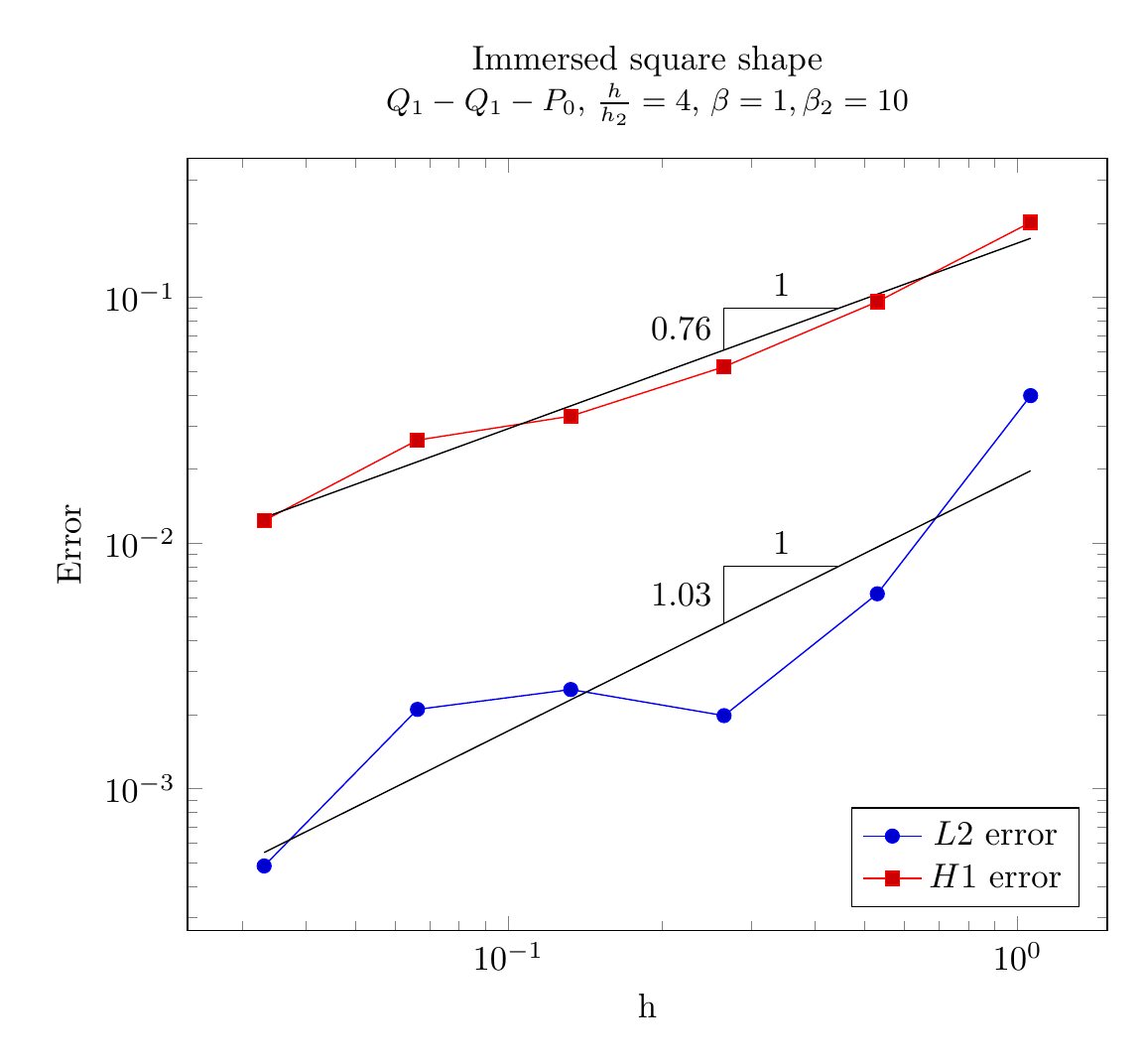} 
\end{minipage}
\begin{minipage}[c]{.3\linewidth}
		\includegraphics[width=1\linewidth]{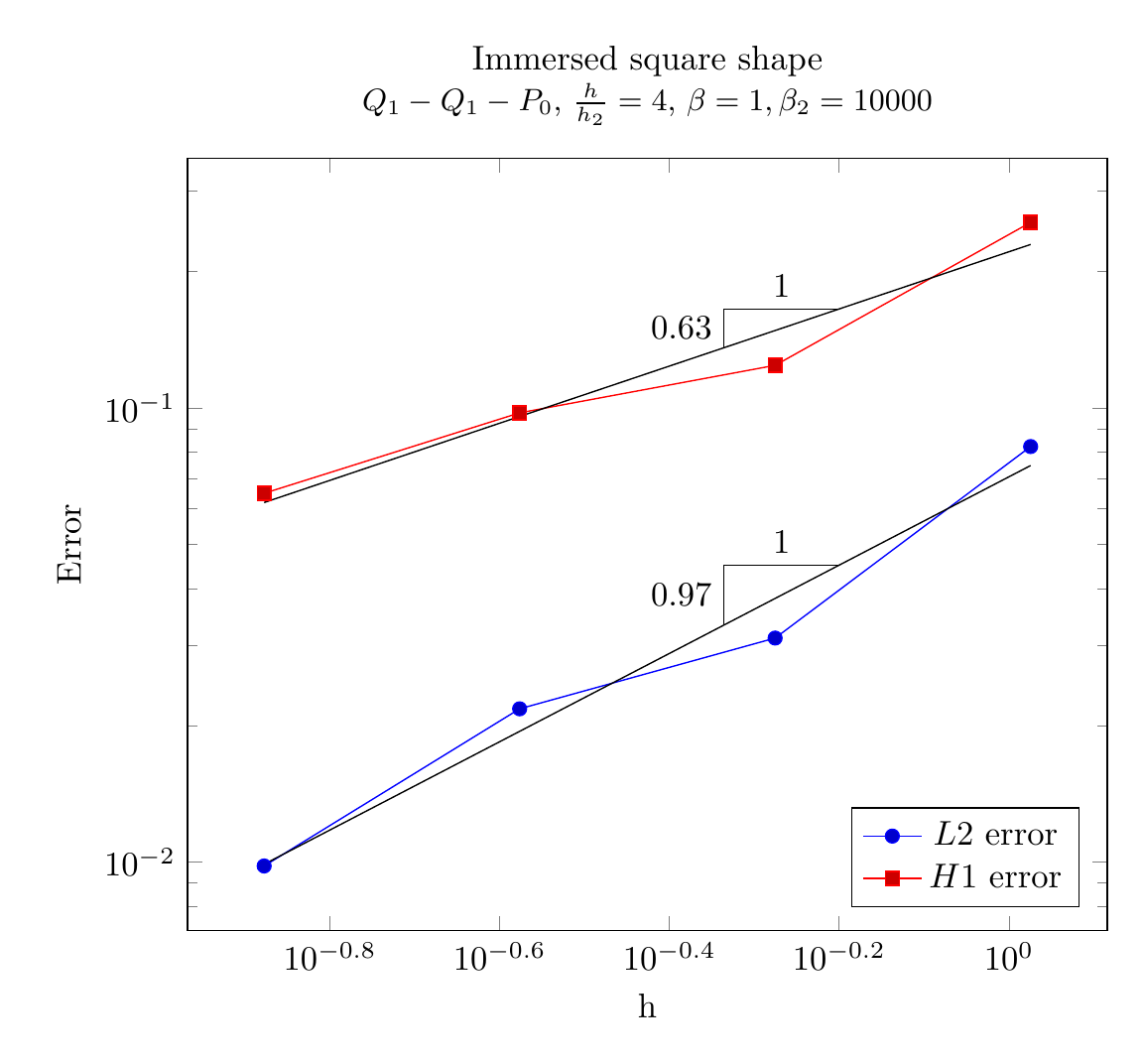} 
\end{minipage}
\begin{minipage}[c]{.3\linewidth}
		\includegraphics[width=1\linewidth]{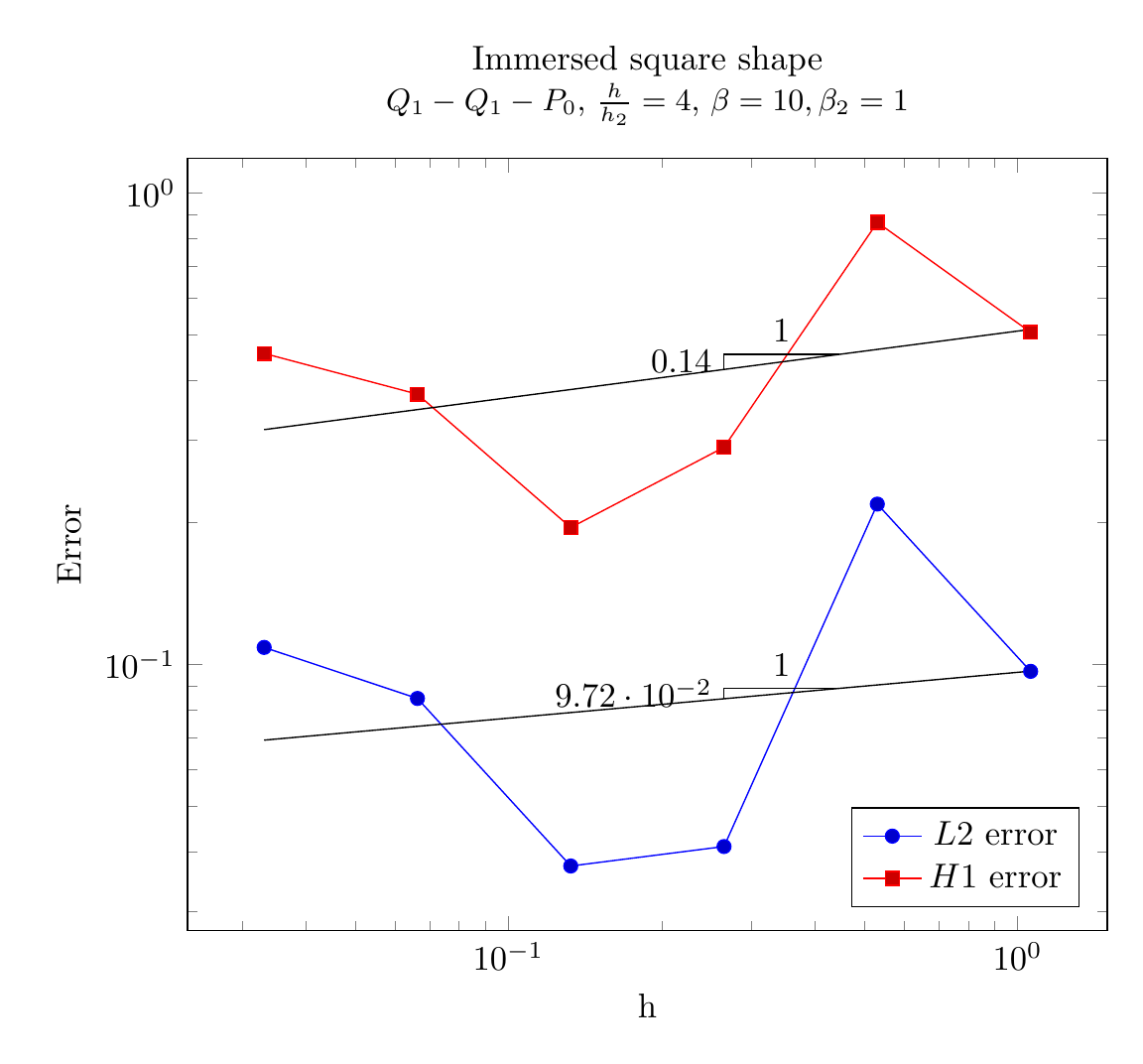} 
\end{minipage}
\captionof{figure}{Non-convergence of $Q_1-Q_1-P_0$ FDDLM: immersed square shape.}
  \label{fig:error110nb_s}
\end{figure}

Given the multitude of possible combinations to test, we have obtained a substantial number of results for each example. While the outcomes are undoubtedly valuable, they tend to be very similar to each other. To ensure a comprehensive coverage of our findings, we will begin by discussing all the different possibilities we have attempted for Example~3 presented in Figures:( \ref{fig:error1_1_0_110_c}-\ref{fig:error220_101_c}). To avoid overwhelming the reader with a deluge of redundant information, we will focus on presenting a representative subset of the most relevant results for each of the other examples. In particular, we present the results using one different ratio between mesh sizes for each as seen in Figures(\ref{fig:error1_1_0_L}-\ref{fig:error220_s}). Lastly, we present three cases using the unstable $(Q_1,Q_1,P_0)$, to show that the rate of convergence does not converge in some situations as illustrated by Figure \ref{fig:error110nb_s}.

The order of convergence of the $L^2$ and $H^1$ norms of the error provide information about how quickly the error decreases as the mesh size is refined. In general, a higher order of convergence indicates that the error decreases more rapidly, which implies a smoother solution. All examples considered encounter jumps in the coefficients that can be quite large. These jumps can pose a challenge for numerical methods, as they can lead to numerical instabilities and affect the accuracy of the solution. Moreover, Examples 1 and 2 exhibit singularities in the geometry due to the re-entrant corners unlike Examples~3 and~4, which can further complicate the problem and possibly lead to numerical instabilities. 

Despite these challenges in the problem of our interest, our proposed numerical scheme is able to achieve a reasonable rate of convergence, demonstrating its stability, robustness, and effectiveness. In general, the results from Cases 1 and 2 for all examples are consistent with our theoretical expectations. Furthermore, the findings from Case 3 indicate that the elker condition constraint may be relaxed to some extent. 

Specifically, we know that the solution belongs to $H^r(\Omega)$ with $1<r<\frac{3}{2}$. Consequently, we would expect the order of convergence to represent a solution that at least lives in $H^1$. The observed orders of convergence of almost $O(h^{1})$ for the $L^2$ norm and almost $O(h^{\frac{1}{2}})$ for the $H^1$ norm of the error are consistent with this expectation.

The convergence rate plots depicted in Figures (\ref{fig:error1_1_0_110_c}-\ref{fig:error220_f}) reveal interesting patterns. Specifically, we observed that in Examples~3 and~4, where the interface is smoother, the convergence plots exhibit more linear behavior. In contrast, Examples~1 and~2, which feature re-entrant corners, exhibit more oscillations in their convergence plots. Interestingly, these oscillations only appear when the ratio between mesh sizes is greater than or equal to 1 in Cases 1 and 2, and the opposite situation in Case 3. Additionally, we found that the $(Q_2,Q_2,P_0)$ element demonstrated linear convergence, whereas the $(Q_1,Q_1+B,P_0)$ element exhibited these oscillations, which were more pronounced in Example~2 than in Example~1.  On the other hand, in Examples~3 and~4, we observed that both elements produced similar results. It is important to note that the rate of convergence is not significantly affected by these oscillations and our scheme is achieving an optimal rate of convergence in all cases.

In the first example, where we have an immersed square shape, we benchmark our results against those obtained in \cite{boffi2014mixed}, when an $L^2$ scalar product is used to evaluate the duality term numerically and found that our results are similar to those reported there.

We conclude this section by going back to the element $Q_1-Q_1-P_0$ due to its appeal coming from its easier implementation.
Recall that, in Section~\ref{sse:weakinfsup}, we showed numerically that this element does not pass the numerical inf-sup test. Unexpectedly, we have an optimal rate of convergence when the
difference $\beta_2-\beta>0$. 
Lastly, when the difference $\beta_2-\beta<0$, this element fails to converge as shown in Figure~\ref{fig:error110nb_s}.
%
\section{Conclusion}
This study proposes two stable elements for solving elliptic interface problems with jump coefficients. Our numerical results confirm that the discrete inf-sup is bounded away from zero for these proposed elements, and the convergence rate is consistent with our theoretical expectations. Specifically, we observe that all cases converge with the optimal convergence rate known in the literature of FEM, but less oscillation appears in the particular ratios between the mesh sizes.
Furthermore, we demonstrate the robustness of our scheme when dealing with different immersed shapes. We observe that Element 2 consistently outperforms Element 1 in all cases. These results extend the findings of \cite{boffi2014mixed} to the case where a discontinuous Lagrange multiplier is used. In future work, it would be interesting to study the a posteriori error estimate for these schemes and use the results to apply local refinements to control the overall error and reduce computational costs by refining only where necessary. Additionally, this work could be extended to a fluid-structure interaction problem, which was the main motivation for our study.


%
\bibliographystyle{apalike}
\bibliography{References}
\end{document}